\documentclass[a4paper,11pt]{article}
\usepackage{a4wide}
\usepackage{theorem}
\usepackage{amsmath}
\usepackage{array}
\usepackage{amssymb}
\usepackage{amsfonts}
\usepackage[french]{babel}
\usepackage{epsf}
\usepackage{epsfig}
\usepackage{graphicx}

\newtheorem{theo}{\indent Th\'eor\`eme \newline}[section]
\newtheorem{defi}[theo]{\indent D\'efinition\newline}
{\theorembodyfont{\rmfamily}%
\newtheorem{rem}[theo]{\noindent Remarque}}
{\theorembodyfont{\rmfamily}%
 \theoremstyle{break}%
}
\newtheorem{prop}[theo]{\indent Proposition\newline}
\newtheorem{lemme}[theo]{\indent Lemme\newline}
\newtheorem{cor}[theo]{\indent Corollaire \newline}

 \def\N{{\mathbb{N}}}
\def\Z{{\mathbb{Z}}}

\def\R{{\mathbb{R}}}
\def\C{{\mathbb{C}}}

\newlength{\indentation}%
\makeatletter
\newcommand\@makefntextsans[1]{%
    \parindent 0em%
    \noindent%
    \hb@xt@0em{\hss}%
    #1}
\def\footnotetextsans{%
     \@ifnextchar [\@xfootnotenextsans%
       {\@footnotetextsans}}
\def\@xfootnotenextsans[#1]{%
  \begingroup%
     \csname c@\@mpfn\endcsname #1\relax%
  \endgroup%
  \@footnotetextsans}
\long\def\@footnotetextsans#1{\insert\footins{%
    \reset@font\footnotesize%
    \interlinepenalty\interfootnotelinepenalty%
    \splittopskip\footnotesep%
    \splitmaxdepth \dp\strutbox \floatingpenalty \@MM%
    \hsize\columnwidth \@parboxrestore%
    \color@begingroup%
      \@makefntextsans{%
        \rule\z@\footnotesep\ignorespaces#1\@finalstrut\strutbox}
    \color@endgroup}}
\makeatother

\begin{document}

\title{Optimalit\'e, congruences et calculs d'invariants des vari\'et\'es symplectiques r\'eelles de dimension quatre}
\author{Jean-Yves Welschinger}
\date{}
\maketitle


{\bf R\'esum\'e:}

Cet article fait suite \`a un pr\'{e}c\'{e}dent dans lequel \'etaient introduits une famille d'invariants par d\'eformation $\chi^d_r$, $d \in H_2 (X ; \Z)$, $r \in \N$, des vari\'et\'es symplectiques 
r\'eelles ferm\'ees de dimension quatre $(X, \omega , c_X)$, invariants qui fournissent des bornes inf\'erieures en g\'eom\'etrie \'enum\'erative r\'eelle. Nous montrons ici par
des m\'ethodes de th\'eorie symplectique des champs que ces bornes inf\'erieures sont optimales lorsque $r \leq 1$ et le lieu r\'eel de la vari\'et\'e contient une sph\`ere, un tore ou
un plan projectif r\'eel (sous des hypoth\`eses plus restrictives dans ce dernier cas). Nous montrons \'egalement qu'une puissance importante de deux divise $\chi^d_r$ pour
des valeurs pas trop grandes de $r$ lorsque le lieu r\'eel contient une sph\`ere ou un plan projectif r\'eel (sous les m\^emes hypoth\`eses plus restrictives dans ce dernier cas)
et proposons enfin quelques calculs explicites dans le cas du plan projectif ou de la quadrique ellipso\"{\i}de ainsi que les formules g\'en\'erales permettant de
les obtenir, lesquelles font intervenir des invariants relatifs pr\'ec\'edemment d\'efinis. 

\section*{Introduction}

Le pr\'esent article fait suite au pr\'{e}c\'{e}dent  \cite{Wels1} dans lequel \'etaient introduits  une famille d'invariants par d\'eformation des vari\'et\'es symplectiques r\'eelles
ferm\'ees de dimension quatre. Une {\it vari\'et\'e symplectique r\'eelle} est une vari\'et\'e symplectique \'equip\'ee 
d'une involution anti-symplectique ; chaque vari\'et\'e projective r\'eelle lisse en fournit un exemple. Ces invariants ont une propri\'et\'e 
imm\'ediate soulign\'ee dans \cite{Wels1}, ils fournissent des bornes inf\'erieures en g\'eom\'etrie \'enum\'erative r\'eelle. Comme son titre l'indique \`a pr\'esent, 
cet article poursuit trois objectifs ; le premier est de montrer l'optimalit\'e de ces bornes inf\'erieures, ce que l'on fera dans plusieurs situations (Th\'eor\`emes \ref{theoopt1} et \ref{theoopt2}), le second 
est de prouver des congruences satisfaites par ces invariants (Th\'eor\`emes \ref{theocong1},  \ref{theocong2} et  \ref{theocong3}) et le dernier de pr\'esenter quelques calculs de 
ces invariants ainsi que
les formules g\'en\'erales permettant de les obtenir (Th\'eor\`emes \ref{theocalcproj},  \ref{theocal2spher} et  \ref{theocal3spher}). Remarquons en passant que les
r\'esultats d'optimalit\'e en g\'eom\'etrie \'enum\'erative r\'eelle se font rares et que les m\'ethodes syst\'ematiques pour y aboutir sont, 
\`a ma connaissance, inexistantes. La m\'ethode syst\'ematique que l'on utilise ici pour aboutir \`a nos r\'esultats vient de la th\'eorie symplectique des champs \cite{EGH}. 

Soit $(X, \omega , c_X)$ une vari\'et\'e symplectique r\'eelle ferm\'ee de dimension quatre et soit $d \in H_2 (X ; \Z)$ une classe d'homologie satisfaisant la relation $(c_X)_* d = -d$.
Choisissons une structure presque complexe auxiliaire $J$ aussi g\'en\'erale que possible parmi les structures $\omega$-positives qui rendent l'involution 
$c_X$ anti-holomorphe. Les {\it courbes $J$-holomorphes rationnelles r\'eelles} homologues \`a $d$, c'est-\`a-dire les sph\`eres $J$-holomorphes invariantes
par $c_X$ homologues \`a $d$, forment alors un espace de dimension $c_1(X)d -1$, o\`u $c_1(X)$ d\'esigne la premi\`ere classe de Chern de la vari\'et\'e $(X, \omega)$.
Nous supposons dans ce travail comme dans  \cite{Wels1} cette dimension positive ou nulle, puisque le cas contraire signifie que l'espace en question est vide, puis faisons
chuter cette dimension \`a z\'ero en imposant quelques contraintes \`a ces courbes, \`a savoir de passer par une collection de $c_1(X)d -1$ points distincts. Ces derniers peuvent
\^etre choisis r\'eels, c'est-\`a-dire fix\'es par $c_X$, ou bien complexes conjugu\'es, c'est-\`a-dire \'echang\'es par $c_X$ ; nous noterons $r$ le nombre de points r\'eels et $r_X$
le nombre de paires de points complexes conjugu\'es, de sorte que $r+2r_X = c_1(X)d -1$. Seul un nombre fini de courbes $J$-holomorphes rationnelles r\'eelles homologues \`a $d$ 
satisfont ces contraintes suppl\'ementaires ; ce nombre d\'epend en g\'en\'eral
des choix auxiliaires de la structure presque complexe et de la configuration de points, essentiellement parce que le corps des r\'eels n'est pas alg\'ebriquement
clos. Toutefois, il ressort de \cite{Wels1} que si l'on compte ces courbes en fonction d'un signe, positif lorsqu'elles ont un nombre pair de points doubles
r\'eels isol\'es et n\'egatif dans le cas contraire, alors l'entier $\chi_r^d$ que l'on obtient est ind\'ependant des choix de la structure presque complexe $J$, de la configuration de points et 
m\^eme de la forme symplectique $\omega$ \`a l'int\'erieur de sa classe de d\'eformation (voir le Th\'eor\`eme $2.1$ de  \cite{Wels1}). Cet entier ne d\'epend que de la classe d'homologie $d$,
du nombre $r$ de points choisis r\'eels et de la r\'epartition de ces points dans les diff\'erentes composantes connexes du lieu r\'eel
$\R X = \text{fix}(c_X)$ de la vari\'et\'e. En fait, la partie r\'eelle d'une sph\`ere holomorphe r\'eelle \'etant connexe, cet invariant $\chi_r^d$ est contraint de s'annuler d\`es que 
ces points ne sont pas tous choisis dans une m\^eme composante $L$ du lieu r\'eel. On adoptera la notation $\chi_r^d (L)$ pour indiquer
que les $r$ points r\'eels sont choisis dans $L$. Le nombre $R_d (\underline{x} , J)$ de courbes $J$-holomorphes rationnelles 
r\'eelles homologues \`a $d$ qui contiennent l'ensemble $\underline{x}$ de points que l'on s'est donn\'e se retrouve ainsi born\'e inf\'erieurement par la valeur absolue de
l'invariant $\chi_r^d (L)$ ; ce sont l\`a les bornes inf\'erieures en g\'eom\'etrie \'enum\'erative r\'eelle que l'on a mentionn\'ees plus haut. Ces bornes s'\'ecrivent
\begin{eqnarray}
\label{bornes}
\vert \chi_r^d (L) \vert  \leq R_d (\underline{x} , J)   \leq N_d,
\end{eqnarray}
comme \'enonc\'ees dans le Corollaire $2.2$ de \cite{Wels1}, le membre $N_d$ d\'esignant le nombre total de courbes $J$-holomorphes rationnelles satisfaisant ces conditions d'incidence 
(c'est un invariant de Gromov-Witten de genre z\'ero de la vari\'et\'e). 

C'est \`a ce stade \`a peu pr\`es que nous a laiss\'e \cite{Wels1} et que l'on reprend ici notre \'etude en appliquant un principe fondamental de la th\'eorie symplectique 
des champs en pr\'esence d'une telle surface lagrangienne $L$ et d'une structure presque complexe $J$ : on allonge le {\it cou} de la structure presque-complexe
au voisinage de $L$ pour lui conf\'erer une longueur arbitrairement grande. Rappelons qu'un voisinage de $L$ dans $X$ est symplectomorphe \`a
un voisinage de la section nulle dans son fibr\'e cotangent $T^*L$, un r\'esultat \'etabli dans \cite{Wein}. \'Etant donn\'ee une m\'etrique riemannienne sur $L$, le fibr\'e unitaire cotangent 
$S^*L = \{ (q,p) \in T^*L \; \vert \, \parallel p \parallel = 1 \}$
muni de la restriction de la forme de Liouville $\lambda$ est une vari\'et\'e de contact de dimension trois. Le compl\'ementaire $T^*L \setminus L$ se trouve \^etre symplectomorphe
\`a la symplectisation $(S^*L \times \R , d(e^t \lambda))$ de cette vari\'et\'e. Ce que l'on appelle cou de longueur arbitrairement grande, c'est une portion
arbitrairement grande $S^*L \times [-n, n]$ de cette symplectisation dans laquelle $J$ envoie le champ de Liouville $\partial / \partial t$ sur le champ de Reeb de $(S^*L , \lambda)$, 
pr\'eserve les plans de contact et est invariante par translation dans le second facteur, voir \cite{EGH} et la strat\'egie g\'en\'erale \'enonc\'ee au \S \ref{subsubsectstrat}. 

Cette technique issue de la th\'eorie symplectique des champs nous permet d'\'etablir les r\'esultats suivants.
Lorsque $L$ est une sph\`ere, un tore ou un plan projectif r\'eel
(mais dans ce dernier cas $(X, \omega , c_X)$  sera elle-m\^eme suppos\'ee symplectomorphe au plan projectif complexe \'eclat\'e en six boules complexes conjugu\'ees au maximum) et lorsque $r_X$ 
est maximal ou en d'autres termes
lorsque $r \leq 1$, les bornes inf\'erieures (\ref{bornes}) sont optimales, atteintes par les structures presque complexes au cou suffisamment long, voir les th\'eor\`emes
\ref{theoopt1} et \ref{theoopt2}. Ainsi, lorsqu'on allonge le cou d'une structure presque complexe en pr\'eservant l'anti-holomorphicit\'e de $c_X$, il arrive
une longueur \`a partir de laquelle toutes les courbes rationnelles r\'eelles sont compt\'ees en fonction d'un m\^eme signe, toutes les \'eliminations possibles entre courbes
s'\'etant r\'ealis\'ees au cours de l'allongement. 
Ce ph\'enom\`ene permet plus g\'en\'eralement d'\'eliminer parfois tous les disques $J$-holomorphes \`a bords dans une lagrangienne, m\^eme en l'absence de structure r\'eelle.
Nous le montrerons dans le cas de sph\`eres  lagrangiennes au paragraphe \ref{subsectmin}
qui fait office de digression, voir les Th\'eor\`emes \ref{theomin3}, \ref{theomin2} et \ref{theoopt4}. Tous ces r\'esultats font l'objet de la premi\`ere partie de cet article. 
Dans la seconde partie,  on d\'emontre qu'une puissance importante de deux divise l'invariant $\chi_r^d$ lorsque $r$ n'est pas trop grand et $L$ est
une sph\`ere ou un plan projectif r\'eel, voir les Th\'eor\`emes \ref{theocong1},  \ref{theocong2} et \ref{theocong3}, le fait que $S^*L$ est un fibr\'e en cercles joue alors un r\^ole important.
Dans la troisi\`eme partie de cet article, on pr\'esente quelques formules permettant le calcul de $\chi_r^d$ dans le plan projectif complexe ou l'ellipso\"{\i}de pour de 
faibles valeurs de $r$, voir les Corollaires \ref{corcalcproj}, \ref{corcalc2spher} et \ref{corcalc3spher}. Ces derni\`eres sont obtenues en brisant la vari\'et\'e en deux morceaux, ce qui brise les 
courbes rationnelles r\'eelles elle-m\^eme
en deux morceaux et permet d'exprimer $\chi_r^d$ en fonction de deux  ingr\'edients, l'un calcul\'e \`a l'aide de courbes r\'eelles dans $T^*L$ qui n'est autre qu'un invariant r\'eel relatif
\`a un diviseur r\'eel sans lieu r\'eel -conique imaginaire pure ou section hyperplane r\'eelle disjointe de l'ellipso\"{\i}de- et l'autre
\`a l'aide de paires de courbes complexes conjugu\'ees dans $X \setminus L$. Les calculs d'invariants relatifs r\'ealis\'es dans \cite{Vak} (voir aussi \cite{IoPar}
et \cite{Katz}) permettent de ma\^{\i}triser
ce deuxi\`eme ingr\'edient. Or, plus $r$ est petit, plus le premier ingr\'edient est simple de sorte que pour les petites valeurs de $r$, on d\'eduit de  \cite{Vak} des
formules de r\'{e}currence g\'en\'erales, voir les Th\'eor\`emes \ref{theocalcproj},  \ref{theocal2spher} et  \ref{theocal3spher}.

Ces r\'esultats d'optimalit\'e, de congruences et de calculs ont \'et\'e annonc\'es dans la note \cite{Wels4} dans le cas du plan  projectif ou de la quadrique de dimension deux ;
ils ont \'et\'e pr\'esent\'es la premi\`ere fois en d\'ecembre $2005$ lors de l'atelier organis\'e en l'honneur de Dusa McDuff, \`a Banff au Canada.\\

{\bf Remerciements :}

Je remercie l'Agence nationale de la recherche pour son soutien ainsi que Y. Eliashberg pour ses encouragements \`a d\'ecouper les vari\'et\'es symplectiques en morceaux.

\tableofcontents

\section{Optimalit\'e}

\subsection{Optimalit\'e des bornes inf\'erieures}
\label{subsectopt}

\subsubsection{\'Enonc\'es des r\'esultats}

Nous \'enon\c{c}ons dans ce premier paragraphe les situations dans lesquelles nous sommes en mesure de montrer l'optimalit\'e des bornes inf\'erieures (\ref{bornes}) en dimension quatre. 
Le paragraphe  \ref{subsectmin}  tiendra lieu de digression en dimension sup\'erieure.

\begin{theo}
\label{theoopt1}
Soit $(X, \omega , c_X)$ une vari\'et\'e symplectique r\'eelle ferm\'ee de dimension quatre et soit $d \in H_2 (X ; \Z)$ une classe d'homologie satisfaisant $(c_X)_* d = -d$.
Supposons que le lieu r\'eel de cette vari\'et\'e poss\`ede une sph\`ere ou un plan projectif r\'eel $L$. Dans ce dernier cas, supposons que $(X, \omega , c_X)$ est elle-m\^eme
symplectomorphe au plan projectif complexe \'eclat\'e en six points complexes conjugu\'es au maximum. Les bornes inf\'erieures (\ref{bornes}) sont sous ces hypoth\`eses
optimales d\`es que $0 \leq r \leq 1$, atteintes par les structures presque-complexes g\'en\'erales ayant un long cou au voisinage de $L$. Le signe de l'invariant $\chi_r^d (L)$ est
en outre dans ce cas d\'etermin\'e par l'in\'egalit\'e $(-1)^{\frac{1}{2}(d^2 - c_1(X)d + 2)} \chi^d_r (L) \geq 0$.
\end{theo}

\begin{rem}
La derni\`ere partie du Th\'eor\`eme \ref{theoopt1} signifie que le signe du coefficient de plus bas degr\'e du polyn\^ome $\chi^d (T)$ introduit dans \cite{Wels1} s'interpr\`ete
comme la parit\'e du genre lisse de la classe $d$. Le fait que ce signe puisse \^etre n\'egatif en degr\'es congrus \`a trois ou quatre modulo quatre dans le plan projectif complexe
met en d\'efaut la Conjecture $6$ de \cite{IKS1}. Nous montrerons en effet au \S \ref{sectcalculs} que cet invariant ne s'annule pas en degr\'es sup\'erieurs \`a cinq, voir le 
Th\'eor\`eme \ref{theocalcproj}
\end{rem}

\begin{cor}
Soit $d$ une classe d'homologie de dimension deux du plan projectif complexe ou de la quadrique ellipso\"{\i}de et $0 \leq r \leq 1$.
Les bornes inf\'erieures (\ref{bornes}) sont
atteintes pour la structure complexe standard lorsque les points complexes conjugu\'es sont choisis tr\`es proches d'une conique imaginaire pure dans le premier cas
et d'une section hyperplane r\'eelle disjointe de $L$ dans le second.
\end{cor}

{\bf D\'emonstration :}

Dans ces deux cas, la structure complexe standard de la vari\'et\'e poss\`ede un cou infiniment long au voisinage de $L$. Il s'agit d'un voisinage fibr\'e en disques de la conique imaginaire
pure ou de la section hyperplane r\'eelle priv\'e de la conique ou de la section elle m\^eme. Comme par ailleurs le plan projectif et la quadrique sont des surfaces convexes, l'hypoth\`ese
de g\'en\'ericit\'e de la structure presque-complexe du Th\'eor\`eme  \ref{theoopt1} est satisfaite (voir les Th\'eor\`emes \ref{theocalcproj} et \ref{theocal2spher} pour un r\'esultat plus g\'en\'eral). 
Le Th\'eor\`eme  \ref{theoopt1} s'applique donc et fournit  le r\'esultat. $\square$

\begin{theo}
\label{theoopt2}
Soit $(X, \omega , c_X)$ une vari\'et\'e symplectique r\'eelle ferm\'ee de dimension quatre dont le lieu r\'eel poss\`ede un tore $L$ et soit $d \in H_2 (X ; \Z)$ une classe d'homologie 
satisfaisant $(c_X)_* d = -d$. Les bornes inf\'erieures (\ref{bornes}) sont
optimales lorsque $r=1$, atteintes par les structures presque-complexes g\'en\'erales ayant un long cou au voisinage de $L$. Lorsque le lieu r\'eel
est connexe -r\'eduit au tore $L$-, l'invariant $\chi_1^d (L)$ est en outre positif. Dans le cas g\'en\'eral, le signe de l'invariant $\chi_1^d (L)$ est
d\'etermin\'e par l'in\'egalit\'e $(-1)^{\frac{1}{2}(d^2 - c_1(X)d + 2)} \chi^d_1 (L) \geq 0$ lorsque le lieu r\'eel des courbes rationnelles ne s'annule pas dans $H_1 (L ; \Z /2\Z)$,
tandis qu'il est d\'etermin\'e par l'in\'egalit\'e $(-1)^{\frac{1}{2}(d^2 - c_1(X)d + 2)} \chi^d_1 (L) \leq 0$ lorsque ce dernier s'annule.
\end{theo}

\begin{rem}
Dans le cas particulier de la quadrique hyperbolo\"{\i}de, la positivit\'e de $\chi_1^d (L)$ a d\'ej\`a \'et\'e observ\'ee dans \cite{IKS1} par d'autres m\'ethodes.
\end{rem}

\subsubsection{Strat\'egie g\'en\'erale}
\label{subsubsectstrat}

\begin{figure}[h]
\label{figarbre}
\begin{center}
\includegraphics{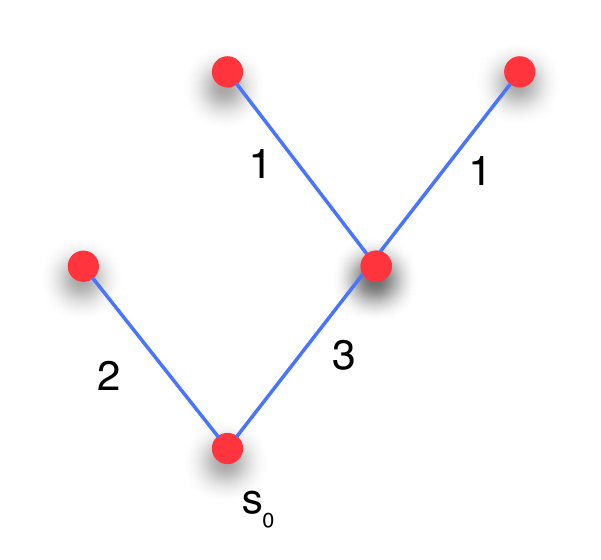}
\end{center}
\caption{Exemple d'arbre $A_C$}
\end{figure}

On allonge le cou d'une structure presque-complexe g\'en\'erique jusqu'\`a briser la vari\'et\'e $(X, \omega , c_X)$ en deux,
le fibr\'e cotangent \`a $L$ d'une part et le compl\'ementaire $X \setminus L$ de l'autre. Chacune de ces deux parties se retrouve munie d'une structure presque-complexe, not\'ee $J_L$ et 
$J_X$ respectivement, qui rendent respectivement $c_L$ et $c_X$ antiholomorphes, o\`u $c_L : (q,p) \in T^* L \mapsto (q, -p)  \in T^* L$. De plus,
en dehors d'un compact, ces structures sont cylindriques sur une structure $CR$ de la vari\'et\'e de contact $(S^* L , \lambda)$. Nous avons ici not\'e $S^* L$ le fibr\'e unitaire cotangent de $L$
pour une m\'etrique \`a courbure constante, de sorte que les orbites p\'eriodiques du 
flot du champ de vecteurs de Reeb $R_\lambda$ associ\'e \`a la forme de Liouville $\lambda$, c'est-\`a-dire du flot g\'eod\'esique, viennent en familles.
Rappelons qu'une fois identifi\'e le compl\'ementaire d'un compact de $T^* L$ ou $X \setminus L$ avec une partie de la symplectisation
$S^* L \times \R$ de $S^* L$, la structure presque complexe  $J_L$ ou $J_X$  est d\'efinie en dehors de ce compact par la structure $CR$ 
et la relation $J \partial / \partial t =  R_\lambda$. Nous r\'ealisons  cette scission de sorte que les $r_X$ paires de points complexes
conjugu\'ees que l'on s'est donn\'ees se retrouvent dans  $X \setminus L$. Le th\'eor\`eme de compacit\'e de 
th\'eorie symplectique des champs \cite{BEHWZ} permet de
comprendre le devenir des courbes rationnelles r\'eelles homologues \`a $d$ et qui passent par $\underline{x}$. Ces courbes se brisent en des courbes \`a deux \'etages, $J_L$-holomorphes
(resp.  $J_X$-holomorphes) pour celles habitant l'\'etage $T^* L$ (resp. $X \setminus L$), et asymptotes \`a des orbites p\'eriodiques de $R_\lambda$, la p\'eriode pouvant \^etre multiple
de la p\'eriode fondamentale. La r\'eunion de ce nombre fini de composantes est invariante par l'involution $c_X$, de sorte que ces composantes sont organis\'ees en paires de
composantes complexes conjugu\'ees de $T^* L$ ou $X \setminus L$ et d'une composante de $T^* L$ laiss\'ee invariante par $c_L$.
Chaque courbe \`a deux \'etages limite $C$ peut donc \^etre cod\'ee par un arbre $A_C$ ayant une racine $s_0$ et ses ar\^etes \'equip\'ees de multiplicit\'es enti\`eres strictement positives.
Chaque sommet de cet arbre repr\'esente une composante du quotient de la courbe limite par l'action de $c_X$, composante qui se trouve dans l'\'etage $T^* L$ si ce sommet est \`a distance paire 
de $s_0$ et dans l'\'etage $X \setminus L$ sinon. Le sommet $s_0$ quant \`a lui repr\'esente l'unique composante
laiss\'ee invariante par l'involution $c_L$ de $T^* L$. Le quotient de cette composante est une h\'emisph\`ere point\'ee \`a bord dans $L$. Chaque ar\^ete
adjacente \`a un sommet donn\'e repr\'esente une asymptote de la composante correspondante \`a ce sommet et la multiplicit\'e de l'ar\^ete n'est autre que la multiplicit\'e de l'orbite de Reeb
limite correspondante. Par exemple, l'arbre repr\'esent\'e par la figure \ref{figarbre} repr\'esente une courbe rationnelle r\'eelle \`a deux \'etages et neuf composantes. La composante
racine est une sph\`ere r\'eelle dans $T^* L$ ayant deux paires de pointes complexes conjugu\'ees asymptotes \`a deux paires d'orbites de Reeb, l'une de multiplicit\'e deux, l'autre de
multiplicit\'e trois.
L'\'etage $X \setminus L$ contient une paire de plans $J_X$-holomorphes complexes conjugu\'es asymptotes \`a la paire d'orbites de Reeb de multiplicit\'e deux pr\'ec\'edente, cette paire est
cod\'ee par la feuille de l'arbre adjacente \`a l'ar\^ete de multiplicit\'e deux. Cet \'etage $X \setminus L$ contient \'egalement une paire de sph\`eres $J_X$-holomorphes complexes conjugu\'ees
ayant trois pointes dont deux sont asymptotes \`a des orbites de Reeb simples et la troisi\`eme asymptote \`a la paire d'orbites de Reeb de multiplicit\'e trois d\'efinie plus haut, cette paire
de sph\`eres est cod\'ee par le sommet trivalent. Enfin, l'\'etage $T^* L$ contient \'egalement deux paires de plans $J_L$-holomorphes complexes conjugu\'es asymptotes aux paires d'orbites de Reeb simples pr\'ec\'edentes, ces plans sont cod\'es par les deux feuilles restantes de l'arbre.
L'arbre $A_C$ vient de plus avec une fonction qui associe \`a chaque sommet \`a distance impaire de $s_0$ les classes d'homologies relatives de la paire de courbes correspondantes ainsi que
les paires de points complexes conjugu\'es de $\underline{x} $ que ces courbes contiennent.\\

\subsubsection{D\'emonstration des Th\'eor\`emes \ref{theoopt1} et \ref{theoopt2}}
\label{subsectoptdem}

Cette strat\'egie g\'en\'erale \'etant  pos\'ee, remarquons que dans chacun des cas qui concernent le Th\'eor\`eme  \ref{theoopt1}, $S^* L$ est un fibr\'e en cercles. Par suite, le fibr\'e normal de chaque 
courbe $C_s$ associ\'e \`a un sommet $s$ de l'arbre $A_C$ est canoniquement trivialis\'e le long des orbites de Reeb asymptotes par le flot de Reeb. Notons $\mu_s$ le double de l'obstruction 
\`a \'etendre cette trivialisation sur $C_s$ tout entier. La dimension de l'espace des modules dans lequel habite $C_s$ s'exprime en fonction de cet indice de Maslov $\mu_s $, voir la 
Proposition \ref{propmaslov} de notre formulaire donn\'e au paragraphe \ref{subsectformulaire}. 
Cette dimension vaut, pour peu que $C_s$ soit une courbe simple, c'est-\`a-dire ne soit pas un rev\^etement ramifi\'e non-trivial d'une autre, $\mu_s + 2$ lorsque $s \neq s_0$ et 
$\frac{1}{2} \mu_s + 1$ lorsque $s = s_0$ puisque la courbe $C_s$ est alors contrainte d'\^etre pr\'eserv\'ee par l'involution $c_X$ ce qui a pour effet de diviser la dimension  par deux.
Notons pour chaque sommet $s$ de l'arbre, sa valence par $v_s$ et la somme des multiplicit\'es des ar\^etes adjacentes par $k_s$. L'indice de Maslov $\mu_s $ s'exprime
pour les sommets \`a distance paire de $s_0$ en fonction de $v_s$ et $k_s$, voir la Proposition \ref{propcotangent} de notre formulaire donn\'e au paragraphe \ref{subsectformulaire}.

Supposons pour commencer que $L$ est une sph\`ere et notons $S_1$ (resp. $S_2$) l'ensemble des sommets \`a distance impaire (resp. paire) de $s_0$. Lorsque 
$s \in S_2 \setminus \{ s_0 \}$,
$\mu_s = 2k_s + 2v_s - 4$ tandis que l'indice de Maslov de l'h\'emisph\`ere associ\'ee \`a $s_0$ vaut $\mu_{s_0} = 2k_{s_0}  + 2v_{s_0}  - 2$. Par suite,
$$\sum_{s \in S_2} \mu_s = 2k + 2v - 4 \# S_2 + 2,$$
o\`u $v$ d\'esigne le nombre total d'ar\^etes de l'arbre et $k$ la somme de leurs multiplicit\'es. Si l'on suppose que toutes les courbes de l'\'etage $X \setminus L$ sont simples, la 
g\'en\'ericit\'e de $J_X$ impose la positivit\'e de toutes les dimensions des espaces de modules intervenant, soit $\mu_s +  2 \geq 0$ pour tout $s \in S_1$. Lorsque la courbe 
$C_s$ contient $f_s$ points de notre configuration, cette condition d'incidence impose l'in\'egalit\'e plus fine $\mu_s +  2 \geq 2f_s$. On d\'eduit au total la minoration
$$\sum_{s \in S_1} \mu_s \geq - 2 \# S_1 + 2r_X.$$
Le nombre d'ar\^etes d'un arbre diff\`ere du nombre de sommets par un, soit
$v =  \# S_1 +  \# S_2 - 1$, de sorte que l'indice de Maslov total de la courbe $C$ satisfait $\mu \geq 2k - 2 \# S_2 + 2r_X \geq 2r_X$.
Or cet indice de Maslov total est par ailleurs major\'e par $c_1 (X) d -2$, le degr\'e du fibr\'e normal d'une courbe rationnelle immerg\'ee homologue \`a $d$.  Par hypoth\`ese, 
ce degr\'e vaut ici $2r_X$ puisque l'orientabilit\'e de $L$ impose l'imparit\'e de $r$. Les minorations pr\'ec\'edentes sont par cons\'equent des \'egalit\'es, de sorte que $k = \# S_2$. 
En particulier, toutes les orbites de Reeb
intervenant sont simplement rev\^etues et tous les sommets de $S_2$ sont des feuilles. La courbe r\'eelle cod\'ee par $s_0$ n'est autre qu'un cylindre r\'eel sur 
une orbite de Reeb simple. Un tel cylindre est n\'{e}cessairement plong\'e, voir le Lemme \ref{lemmepointsdoubles} de notre formulaire. Le r\'esultat en d\'ecoule ; peu 
avant la brisure de la vari\'et\'e, toutes les courbes rationnelles r\'eelles ont leurs parties r\'eelles plong\'ees, de sorte que les points doubles r\'eels \'eventuels de ces courbes sont 
tous isol\'es. Ce nombre de points doubles est de m\^eme parit\'e que le genre lisse de la courbe. 

Il s'agit \`a pr\'esent d'aboutir \`a la m\^eme conclusion sans supposer que les courbes $C_s$ soient simples. L'indice de Maslov $\mu^l$ d'un rev\^etement de degr\'e $l$ d'une courbe
simple d'indice $\mu$ s'\'ecrit $\mu^l = l \mu + 2R$ o\`u $R$ est l'indice de ramification. Cet indice de Maslov peut \^etre strictement plus petit que $\mu$ uniquement lorsque 
$\mu$ est n\'egatif, donc \'egal \`a $-2$ et encore
faut-il que la courbe rev\^etue ne soit pas plane. Cela ne concerne donc ni les courbes de l'\'etage $T^* L$, ni les courbes de $X \setminus L$ soumises \`a des conditions d'incidence. Notons
$s_1 , \dots , s_j$ les sommets de $A_C$ correspondant \`a ces derni\`eres et calculons la contribution \`a l'indice de Maslov total de chaque composante connexe de l'arbre priv\'e des sommets
$s_1 , \dots , s_j$. Pour ce faire, notons $S_1'$ (resp. $S_2'$) l'ensemble des sommets \`a distance impaire (resp. paire) de $s_0$ d'une telle composante connexe de $A_C \setminus \{ s_1 , \dots , s_j \}$.
Comme pr\'ec\'edemment, $\sum_{s \in S_2'} \mu_s = 2k + 2v - 4 \# S_2 + 2\delta$, o\`u $\delta$ vaut un si la composante en question contient $s_0$ et z\'ero sinon, tandis que
\begin{eqnarray}
\sum_{s \in S_1'} \mu_s & = & \sum_{s \in S_1'} (l_s \tilde{\mu}_s + 2R_s) \text{  o\`u } l_s \text{  d\'esigne le degr\'e du rev\^etement, } R_s \text{  l'indice de ramification} \nonumber \\
&& \text{et } \tilde{\mu}_s \text{  l'indice de Maslov de la courbe simple sous-jacente} \nonumber \\
& \geq &  -2 \sum_{s \in S_1'} l_s + 2 \sum_{s \in S_1'} (l_s \tilde{\chi}_s - \chi_s) \text{  o\`u } \chi \text{  d\'esigne la caract\'eristique d'Euler,}  \nonumber \\
& \geq &  2 \sum_{s \in S_1'} (l_s -  l_s \tilde{v}_s  +  v_s) - 4 \# S_1' \text{  o\`u } v_s \text{  d\'esigne le nombre de pointes}  \label{equ3} \\
& \geq &  2 \sum_{s \in S_1'} (l_s -  k_s  + v_s) - 4 \# S_1'. \label{equ1}
\end{eqnarray}
Apr\`es sommation, on d\'eduit $\sum_{s \in S_1' \cup S_2'} \mu_s  \geq 2 \sum_{s \in S_1'} l_s + 2k' + 2v' -  4+ 2\delta$, o\`u $v'$ et $k'$ d\'esignent respectivement le nombre 
d'ar\^etes attach\'ees \`a $s_1 , \dots , s_j$ et leur multiplicit\'e totale. Les minorations  $\sum_{s \in S_1' \cup S_2'} \mu_s  \geq  2k'$ si $\delta$ vaut un et 
$\sum_{s \in S_1' \cup S_2'} \mu_s  \geq  2(k' - 1)$ sinon en r\'esultent. La contribution totale \`a l'indice de Maslov des sommets autres que $s_1 , \dots , s_j$ se trouve donc 
minor\'ee par $2j$. La contribution des sommets $s_1 , \dots , s_j$
est quant \`a elle minor\'ee par $2r_X - 2j$, de sorte qu'on aboutit \`a nouveau \`a la minoration $\mu \geq  2r_X$. On conclut donc comme pr\'ec\'edemment.

Supposons \`a pr\'esent que $L$ est un plan projectif r\'eel. D'apr\`es la Proposition  \ref{propcotangent} de notre formulaire donn\'e au paragraphe \ref{subsectformulaire}, l'indice 
de Maslov d'un sommet  $s \in S_2 \setminus \{ s_0 \}$ vaut $\mu_s = k_s + 2v_s - 4$ tandis que l'indice de Maslov 
de l'h\'emisph\`ere associ\'ee \`a $s_0$ vaut $\mu_{s_0} = k_{s_0}  + 2v_{s_0}  - 2$ de sorte que
$\sum_{s \in S_2} \mu_s = k + 2v - 4 \# S_2 + 2$. Les hypoth\`eses faites sur la vari\'et\'e garantissent l'absence de courbes simples d'indices de Maslov strictement n\'egatifs autres
que des plans dans l'\'etage $X \setminus L$. En effet, cet \'etage est isomorphe au fibr\'e en droites complexes de degr\'e quatre sur la conique imaginaire pure \'eclat\'e en 
six points complexes
conjugu\'es au maximum. La classe d'homologie relative d'une courbe dans cet espace s'\'ecrit $ne + kf - \sum_i \alpha_i E_i$, o\`u $e$ est la section nulle du fibr\'e, $f$ une fibre et
$E_i$ les \'eventuels diviseurs exceptionnels. L'irr\'eductibilit\'e de la courbe $C_s$ impose les in\'egalit\'es $\alpha_i \leq n$ d\`es que $n \geq 1$, ce que l'on obtient comme 
cons\'equence
de la positivit\'e d'intersection avec les courbes exceptionnelles $J_X$-holomorphes $E_i$ et $f - E_i$. L'indice de Maslov d'une telle courbe vaut $2(6n + 2k -  \sum_i \alpha_i -2)$,
il est positif d\`es que $k,n$ sont non nuls. Pour chaque sommet $s \in S_1$ l'in\'egalit\'e $\mu_s +  2 \geq 0$ s'en d\'eduit. Lorsque la courbe $C_s$ contient $f_s$ points de 
notre configuration, cette condition d'incidence impose l'in\'egalit\'e plus fine $\mu_s +  2 \geq 2f_s$. De l\`a la minoration $\sum_{s \in S_1} \mu_s \geq - 2 \# S_1 + 2r_X$ et 
finalement apr\`es sommation l'estimation de l'indice de Maslov total $\mu \geq k - 2 \# S_2 + 2r_X $. Il reste \`a remarquer que pour chaque $s \in S_2 \setminus \{ s_0 \}$, l'entier
$k_s$ doit \^etre pair puisque le noyau du morphisme $H_1 (S^* L ; \Z) \to H_1 ( L ; \Z)$ est engendr\'e par une orbite de Reeb double. L'in\'egalit\'e pr\'ec\'edente se r\'e\'ecrit 
donc  \`a pr\'esent
$\mu \geq 2r_X $ si les parties r\'eelles des courbes rationnelles que l'on consid\`ere sont non nulles dans $H_1 ( L ; \Z)$  et $\mu \geq 2r_X - 1$ sinon.
Or cet indice de Maslov est par ailleurs major\'e par le degr\'e $c_1 (X) d -2$ du fibr\'e normal \`a une courbe rationnelle immerg\'ee homologue \`a $d$, degr\'e qui par hypoth\`ese 
vaut ici $2r_X + r - 1$.
Ainsi, toutes les minorations pr\'ec\'edentes sont des \'egalit\'es, de sorte que les sommets \`a distances paires de $s_0$ sont soit des cylindres sur des orbites simples, soit des plans
sur des orbites de Reeb doubles. Le sommet $s_0$ quant \`a lui code un cylindre r\'eel sur une orbite simple lorsque $r$ est nul, et soit un cylindre r\'eel sur une orbite double,
soit une sph\`ere r\'eelle ayant deux paires de pointes complexes conjugu\'ees asymptotes \`a des orbites simples lorsque $r$ vaut un. Dans tous ces cas, une telle courbe est plong\'ee,
de par le Lemme \ref{lemmepointsdoubles}. On conclut comme pr\'ec\'edemment, ce qui ach\`eve la d\'emonstration du Th\'eor\`eme \ref{theoopt1}.

Supposons enfin que $L$ soit un tore et munissons-le d'une m\'etrique plate de sorte que son fibr\'e unitaire cotangent $(S^* L , \lambda)$ soit un tore standard de dimension trois. 
Le flot de Reeb fournit \`a nouveau une
trivialisation canonique du fibr\'e normal aux courbes $C_s$ le long de leurs orbites de Reeb asymptotes. L'obstruction $\mu_s$  \`a \'etendre cette trivialisation sur $C_s$ tout entier
vaut cette fois-ci $2v_s -4$ si $s \neq s_0$ est \`a distance paire de $s_0$ et $2v_s -2$ si $s = s_0$, voir la Proposition  \ref{propcotangent} de notre formulaire donn\'e au 
paragraphe \ref{subsectformulaire}. Si $s$ est au contraire \`a distance impaire de $s_0$, la dimension de l'espace des modules dans lequel habite $C_s$ s'\'ecrit $\mu_s + 2 - v_s$, 
d'apr\`es la Proposition \ref{propmaslov}. Contrairement aux 
cas pr\'ec\'edents, la passage \`a un
rev\^etement ramifi\'e ne peut faire qu'augmenter cette dimension. On d\'eduit donc de la parit\'e de $\mu_s$ l'in\'egalit\'e $\mu_s \geq 0$,  in\'egalit\'e stricte lorsque $v_s > 2$. Si la courbe
est contrainte de passer par $f_s$ points de notre configuration, cette in\'egalit\'e se trouve renforc\'ee en $\mu_s \geq 2f_s$. En sommant les contributions de tous les sommets de
l'arbre, on s'aper\c{c}oit donc que l'indice de Maslov total $\mu$ de la courbe est minor\'e par $2r_X$. Comme cet indice de Maslov est par ailleurs major\'e par $c_1 (X) d -2$ et comme 
par hypoth\`ese
$r=1$, la minoration pr\'ec\'edente est une \'egalit\'e. Il en est par suite de m\^eme pour toutes les minorations faites, de sorte que toutes les composantes de l'\'etage $T^* L$ sont
des cylindres. Les cylindres autres que celui associ\'e \`a $s_0$ sont disjoints de $L$ pour un choix g\'en\'erique de $J_X$. Le cylindre r\'eel associ\'e \`a $s_0$ est un rev\^etement
d'un cylindre plong\'e sur une orbite de Reeb simple. En effet, quitte \`a passer \`a un rev\^etement du fibr\'e cotangent \`a $L$, on peut supposer le cylindre asymptote \`a une orbite
simple. Un tel cylindre est, une fois l'orbite fix\'ee, unique et plong\'e, ce qui est imm\'ediat pour la structure complexe standard de $T^* L$ et est une propri\'et\'e invariante par 
d\'eformation de la structure presque-complexe. On en d\'eduit que peu avant la brisure de la vari\'et\'e, toutes les courbes
rationnelles r\'eelles que l'on consid\`ere poss\'edaient un nombre de points doubles r\'eels non-isol\'es pair si le degr\'e du rev\^etement
est impair et impair sinon. En effet, la perturbation du rev\^etement $k$-uple d'une courbe simple du tore produit $k - 1$ points d'auto-intersection modulo deux.
Le nombre de points doubles r\'eels isol\'es de ces courbes rationnelles se trouve donc \^etre de la m\^eme parit\'e que le genre lisse de la courbe lorsque le lieu r\'eel des courbes
rationnelles est non-nul dans $H_1 (L ; \Z / 2\Z)$ et de la parit\'e oppos\'ee lorsque celui-ci s'annule. Le Th\'eor\`eme \ref{theoopt2} est d\'emontr\'e. $\square$

\subsection{Minimisation du nombre de membranes $J$-holomorphes}
\label{subsectmin}

Soit $C$ une membrane $J$-holomorphe \`a bord dans une sous-vari\'et\'e lagrangienne $L$ d'une vari\'et\'e symplectique ferm\'ee $(X, \omega)$. Notons $\chi$ la caract\'eristique d'Euler de
cette membrane, $d \in H_2 (X , L ; \Z)$ sa classe d'homologie relative et $\mu_{TX} \in H^2 (X , L ; \Z)$ la classe de Maslov de la paire $(X , L)$. La dimension attendue de l'espace des d\'eformations
de $C$ s'\'ecrit $\langle \mu_{TX} , d \rangle + (n-3)\chi$. Cette dimension chute lorsque l'on impose \`a $C$ des contraintes suppl\'ementaires. Si l'on impose par exemple
\`a cette membrane de rencontrer $p$ cycles de codimensions $2 + q_1, \dots , 2 + q_p$, cette dimension attendue chute de la somme $q = q_1+ \dots + q_p$. Deux probl\`emes g\'en\'eraux sous-tendent
nos r\'esultats. Il s'agit d'une part de compter le nombre de membranes $J$-holomorphes homologues \`a $d$ soumises \`a de telles conditions d'incidence de sorte que ce comptage ne d\'epende pas
de $J$ et ne d\'epende des conditions d'incidence qu'\`a homologie pr\`es. Il s'agit d'autre part de minimiser ce nombre de membranes. Si nous ne pouvons r\'epondre au premier probl\`eme
dans ce degr\'e de g\'en\'eralit\'e, il nous est par contre parfois possible de r\'epondre au second sans m\^eme supposer l'\'egalit\'e $q = \langle \mu_{TX} , d \rangle + (n-3)\chi$, lorsque le minimum en question est nul. Le pr\'esent paragraphe est consacr\'e aux  r\'esultats que l'on a pu obtenir
dans cette direction. Ici encore le minimum est atteint en allongeant le cou d'une structure presque complexe g\'en\'erale.

\subsubsection{En dimension sup\'erieure}

\begin{theo}
\label{theomin3} 
Soit $L$ une sph\`ere lagrangienne dans une vari\'et\'e symplectique ferm\'ee $(X, \omega)$ satisfaisant $c_1 (X) = \lambda \omega$, $\lambda \leq 0$
et soit $E > 0$. Supposons la dimension de $X$ sup\'erieure \`a cinq. Pour toute structure presque-complexe $J$ g\'en\'erale ayant un cou suffisamment long au voisinage de $L$,
cette vari\'et\'e ne poss\`ede ni membrane $J$-holomorphe reposant sur $L$ ni courbe $J$-holomorphe rencontrant $L$ qui soit d'\'energie inf\'erieure \`a $E$. 
Ce r\'esultat reste valable en dimension quatre pour les courbes ou membranes de genre nul. 
\end{theo}
Rappelons que l'\'energie d'une courbe $C$ est par d\'efinition l'int\'egrale de la forme $\omega$ sur cette courbe. Les vari\'et\'es projectives \`a fibr\'e canonique nul ou ample, par 
exemple les intersections compl\`etes de multidegr\'es $(d_1, \dots , d_k)$
de l'espace projectif de dimension $N$ d\`es lors que $\sum_{i=1}^k d_i \geq N+1$, satisfont les hypoth\`eses du Th\'eor\`eme \ref{theomin3}. Remarquons qu'une modification de ce dernier
s'applique \'egalement aux vari\'et\'es dont le fibr\'e canonique est le produit d'un fibr\'e ample et d'un fibr\'e port\'e par un diviseur effectif disjoint de $L$. Le Th\'eor\`eme \ref{theomin3} permet
de d\'efinir l'homologie de Floer de deux sph\`eres lagrangiennes proches dans les vari\'et\'es symplectiques dont la premi\`ere classe de Chern s'annule, j'esp\`ere d\'evelopper ce r\'esultat 
prochainement.

\begin{theo}
\label{theomin2} 
Soit $L$ une sph\`ere lagrangienne dans une vari\'et\'e symplectique ferm\'ee semipositive $(X, \omega)$ de dimension $2n \geq 6$ et soit $d \in H_2 (X , L ; \Z)$. \'Ecrivons 
$\langle \mu_{TX} , d \rangle + (n-3)\chi = q + r$ avec $q \in \Z$, $0 \leq r < 2 + (n-3)\chi$ et $\chi \leq 2$. Lorsque $q \geq 0$, choisissons $p$ cycles de $X \setminus L$ de codimensions 
$2 + q_1, \dots , 2 + q_p$ de sorte que $q = q_1+ \dots + q_p$.  D\`es que la structure presque complexe g\'en\'erale $J$ poss\`ede un cou suffisamment long au voisinage 
 de $L$, cette vari\'et\'e ne contient aucune membrane $J$-holomorphe homologue \`a $d$, de caract\'eristique d'Euler $\chi$ qui rencontre ces $p$ cycles et repose sur $L$.
Ce r\'esultat reste valable pour des membranes de genre nul lorsque $n = 2$.
\end{theo}

{\bf Exemple : la quadrique ellipso\"{\i}de.}

Soit $X$ la quadrique ellipso\"{\i}de de dimension complexe $n \geq 3$ et $H$ une section hyperplane disjointe de $L$. Le groupe $H_2 (X , L ; \Z)$ est monog\`ene, engendr\'e par
la classe $d_0$ satisfaisant $\langle H , d_0 \rangle = +1$. La premi\`ere classe de Chern de $X$ vaut $n H$, d'o\`u l'on d\'eduit le calcul $\langle \mu_{TX} , l d_0 \rangle
= 2ln$ quel que soit l'entier $l$. \'Ecrivons $l = (n-1)a + b$, le Th\'eor\`eme \ref{theomin2} s'applique par exemple lorsque $n+1 \leq 2b < 2n$, les membranes sont des disques
et lorsque toutes les conditions d'incidence sont ponctuelles.

\begin{theo}
\label{theoopt4} 
Soit $(X, c_X)$ une vari\'et\'e alg\'ebrique r\'eelle convexe de dimension trois dont le lieu r\'eel poss\`ede une sph\`ere $L$. Supposons l'existence d'une classe d'homologie 
 $d \in H_2 (X ; \Z)$ satisfaisant $c_1 (X) d = 2 \mod (4)$. L'invariant $\chi_1^d (L)$ est alors n\'egatif et les bornes inf\'erieures (\ref{bornes}) sont optimales. Dans le cas de
 l'ellipso\"{\i}de, ces bornes sont atteintes pour la structure complexe alg\'ebrique lorsque les conditions d'incidence non r\'eelles sont choisies suffisamment proches d'une section 
 hyperplane r\'eelle disjointe de $L$.
\end{theo}

Il se peut que l'ellipso\"{\i}de de dimension trois soit en fait le seul exemple de vari\'et\'e satisfaisant les hypoth\`eses du Th\'eor\`eme \ref{theoopt4}. L'invariant
$\chi_1^d (L)$ qui appara\^{\i}t dans ce  th\'eor\`eme a \'et\'e construit dans \cite{Wels2}. Remarquons qu'en reprenant les notations du Th\'eor\`eme \ref{theomin2}, 
ce Th\'eor\`eme \ref{theoopt4} traite du cas $r=n-1$ et montre ainsi en un sens l'optimalit\'e des hypoth\`eses faites dans ce Th\'eor\`eme \ref{theomin2}.
Nous montrerons en effet dans la troisi\`eme partie de cet article la non-trivialit\'e de l'invariant $\chi_1^d (L)$  pour l'ellipso\"{\i}de de dimension trois et calculerons ce dernier,
voir le \S \ref{subsec3spher}. \\

{\bf D\'emonstration des Th\'eor\`emes  \ref{theomin3}, \ref{theomin2} et \ref{theoopt4} :}

Nous suivons la strat\'egie g\'en\'erale \'enonc\'ee au paragraphe \ref{subsubsectstrat} en \'equipant la sph\`ere lagrangienne d'une m\'etrique \`a courbure constante. 
Les \'eventuelles membranes qui survivraient \`a l'allongement du cou de $J$ jusqu'\`a la brisure de
la vari\'et\'e seraient cette fois-ci cod\'es par des graphes $A_C$ ayant $b+1$ sommets marqu\'es $s_0, \dots , s_b$ correspondant aux 
composantes ayant un bord dans $L$. Les sommets \`a distances paires de $s_0, \dots , s_b$ codent \`a nouveau les composantes de l'\'etage $T^* L$ et les sommets \`a distances impaires
les composantes  de l'\'etage $X \setminus L$. Le flot de Reeb trivialise le fibr\'e normal de chaque composante $C_s$ associ\'ee au sommet $s$
d'un graphe $A_C$ le long de ses orbites de Reeb limites. Notons $\mu_s$ le double de l'obstruction \`a \'etendre cette trivialisation sur $C_s$ toute enti\`ere.
Notons \'egalement, pour chaque sommet $s$ du graphe, sa valence par $v_s$, la somme des multiplicit\'es des ar\^etes adjacentes par $k_s$ et la caract\'eristique d'Euler 
de la courbe qu'il code par $\chi_s$.
L'indice de Maslov des composantes cod\'ees par les sommets \`a distances paires de $s_0, \dots , s_b$, c'est-\`a-dire des courbes $C_s$ de l'\'etage $T^* L$, s'exprime 
d'apr\`es  la Proposition \ref{propcotangent} de notre formulaire par la relation $\mu_s = 2(n-1)k_s - 2\chi_s$.
Pour calculer la contribution totale des sommets \`a distances impaires de $s_0$,
il faut tenir compte du fait que certaines des composantes associ\'ees peuvent rev\^etir des courbes simples.
Notons pour chacun de ces sommets $l_s$ le degr\'e du rev\^etement, $\tilde{\mu}_s$ l'indice de Maslov de la courbe simple sous-jacente et $\tilde{\chi}_s$ sa caract\'eristique d'Euler.
D'apr\`es la Proposition \ref{propmaslov} de notre formulaire, la  dimension de l'espace des modules dans lequel habite cette courbe simple vaut 
$\tilde{\mu}_s + (n-1)(\tilde{\chi}_s + \tilde{v}_s)$. La g\'en\'ericit\'e de la structure presque complexe assure donc la minoration $\tilde{\mu}_s \geq - (n-1)(\tilde{\chi}_s + \tilde{v}_s)$. Les courbes simples
sous-jacente \'etant soumises \`a nos $p$ conditions d'incidence, cette derni\`ere minoration peut apr\`es sommation \^etre am\'elior\'ee de $q$. Par cons\'equent,
\begin{eqnarray*}
\sum_{s \in S_1} \mu_s & = & \sum_{s \in S_1} \big( l_s (\tilde{\mu}_s + 2 \tilde{\chi}_s) - 2 \chi_s \big) \\
& \geq &  q - (n-3)\sum_{s \in S_1} l_s (\tilde{\chi}_s + \tilde{v}_s) -  2\sum_{s \in S_1} (k_s + \chi_s) \; \text{Ê puisque }  l_s \tilde{v}_s \leq k_s. \\
\end{eqnarray*}
Nous en d\'eduisons
$$2 \chi  + \sum_{s \in S_1 \cup S_2} \mu_s   \geq  q +  2(n-2)k -(n-3)\sum_{s \in S_1} l_s (\tilde{\chi}_s + \tilde{v}_s),$$ 
o\`u $k = \sum_{s \in S_1} k_s$.
Lorsque $n \geq 3$, utilisant les majorations $\tilde{\chi}_s + \tilde{v}_s \leq 2$ et $l_s \leq k_s$, nous aboutissons \`a $\sum_{s \in S_1 \cup S_2} \mu_s  + 2 \chi  \geq  q +  2$. Lorsque $n=2$,
nos hypoth\`eses imposent $\tilde{\chi}_s + \tilde{v}_s = 2$ de sorte qu'\`a nouveau $\sum_{s \in S_1 \cup S_2} \mu_s  + 2 \chi  \geq  q +  2$. 
Le Th\'eor\`eme \ref{theomin2} suppose la vari\'et\'e semipositive, les \'eventuelles composantes compactes de l'\'etage $X \setminus L$ ont donc un indice de Maslov positif. 
Par cons\'equent, l'indice de Maslov total
satisfait la majoration $2 \chi  + \sum_{s \in S_1 \cup S_2} \mu_s  \leq \langle \mu_{TX} , d \rangle  \leq   q + r - (n - 3)\chi < q + 2$. Ces minoration et majoration \'etant incompatibles, aucune membrane
ne peut survivre jusqu'\`a la brisure de la vari\'et\'e. Le Th\'eor\`eme \ref{theomin2} est d\'emontr\'e. Dans le cas du Th\'eor\`eme \ref{theomin3}, $q=0$ et nous d\'eduisons par recollement des
composantes cod\'ees par le graphe $A_C$ une membrane symplectique $C$ de $(X,L)$ d'indice de Maslov $ \langle \mu_{TX} , [C] \rangle \geq 2$. Or par hypoth\`ese,
$\langle \mu_{TX} , [C] \rangle = 2 \langle c_1(X) , c \rangle = 2 \lambda \langle \omega , c \rangle \leq 0$, o\`u $c \in H_2 (X ; \Z)$ rel\`eve $[C] \in H_2 (X , L ; \Z)$. Cette impossibilit\'e
d\'emontre le Th\'eor\`eme \ref{theomin3}.

Le Th\'eor\`eme \ref{theoopt4} correspond au cas o\`u $r = n-1$. Dans, ce cas, les minoration et majoration pr\'ec\'edentes co\"{\i}ncident, de sorte que toutes les in\'egalit\'es sont des \'egalit\'es.
En particulier, $k=1$ de sorte que chaque graphe $A_C$ se trouve r\'eduit \`a deux sommets reli\'es par une ar\^ete simple. La courbe r\'eelle cod\'ee par $s_0$ est
un cylindre sur une orbite de Reeb simple. L'\'etat spinoriel de ces courbes se calcule comme suit. En perturbant le point r\'eel dans toutes les directions dans $L$, on s'aper\c{c}oit que
toutes ces courbes ont le m\^eme \'etat spinoriel qu'une conique obtenue comme section plane r\'eelle de la quadrique ellipso\"{\i}de r\'eelle. Ce dernier vaut $-1$ comme on le
v\'erifie en d\'eformant l'\'equateur vers un parall\`ele proche d'un p\^ole de $L$. $\square$

\subsubsection{En dimension quatre}

Nous noterons ${\cal M}_{g,b}$ l'espace des modules des structures complexes de la surface compacte connexe orient\'e de genre $g$ ayant $b$ composantes de bord.

\begin{prop}
\label{propminsphere}
Soit $L$ une sph\`ere lagrangienne dans une vari\'et\'e symplectique ferm\'ee de dimension quatre $(X, \omega)$. On suppose que cette derni\`ere ne poss\`ede pas de sph\`ere
symplectique $S$ satisfaisant $\langle c_1 (X) , [S] \rangle > 0$. Soit $(d, g, b) \in H_2 (X , L ; \Z) \times \N \times \N^*$ et $K$ un compact de ${\cal M}_{g,b}$. Alors, pour toute
structure presque-complexe g\'en\'erale ayant un cou suffisamment long au voisinage de $L$, la vari\'et\'e ne poss\`ede pas de membrane $J$-holomorphe homologue \`a $d$ \`a 
bord dans $L$ et conforme \`a un \'el\'ement de $K$.
\end{prop}

{\bf D\'emonstration de la Proposition \ref{propminsphere} :}

On poursuit la strat\'egie g\'en\'erale d\'ecrite au paragraphe \ref{subsubsectstrat} pr\'ec\'edent en \'equipant $L$ d'une m\'etrique \`a courbure constante
et en allongeant le cou d'une structure presque complexe g\'en\'erique jusqu'\`a briser la vari\'et\'e en deux morceaux.
D'apr\`es le th\'eor\`eme de compacit\'e de th\'eorie symplectique des champs \cite{BEHWZ}, les membranes que l'on consid\`ere se brisent en courbes \`a deux \'etages qui sont cette fois-ci cod\'ees 
par des graphes $A_C$ ayant $b$ sommets marqu\'es $s_1, \dots , s_b$ correspondant aux 
$b$ composantes de bord. Les sommets \`a distances paires de $s_1, \dots , s_b$ codent \`a nouveau les composantes de l'\'etage $T^* L$ et les sommets \`a distances impaires
les composantes  de l'\'etage $X \setminus L$. Par hypoth\`ese, l'\'etage $X \setminus L$ ne poss\`ede pas de courbe $J$-holomorphe rationnelle asymptote \`a des orbites
de Reeb du fibr\'e unitaire cotangent $S^* L$. En effet, une telle courbe $J$-holomorphe rationnelle simple $C$ aurait d'apr\`es la Proposition \ref{propmaslov} un indice de Maslov
$\mu \geq -2$. Notons $v \geq 1$ le nombre de pointes asymptotes \`a des orbites de Reeb de $S^* L$ et $\chi (C)$ la caract\'eristique d'Euler de $C$. En recollant \`a $C$ en chacune
de ses pointes un plan $J$-holomorphe de $T^* L$, on obtient une sph\`ere symplectique $S$ de $X$. Le fibr\'e tangent \`a $X$ est trivialis\'e le long des pointes de $C$ par le flot de Reeb. 
D'apr\`es ce qui pr\'ec\`ede, le double de l'obstruction \`a  \'etendre cette trivialisation le long de $C$ vaut $\mu + 2 \chi (C)$ alors qu'elle vaut deux le long de chaque plan de $T^* L$ d'apr\`es
la Proposition \ref{propcotangent}. Finalement, l'indice de Maslov de $S$ vaudrait $\mu + 4 \geq 2$, ce qui est exclu par les hypoth\`eses. Remarquons \`a pr\'esent que 
chaque composante des courbes \`a deux \'etages est asymptote \`a une r\'eunion de cylindres $J$-holomorphes sur les orbites de Reeb limites. Ces cylindres ont un module infini.
On en d\'eduit que peu avant la brisure de la vari\'et\'e $X$, lorsque $J$ poss\`ede un cou extr\^emement long, les membranes $J$-holomorphes poss\`edent \'egalement des anneaux de grands 
modules dont les \^ames sont homotopes aux orbites de Reeb cod\'ees par les ar\^etes de l'arbre $A_C$. Au moins un de ces anneaux ne borde pas de disque, lequel proviendrait n\'ecessairement
d'un plan de $T^* L$, puisque les membranes ont un bord dans $L$. Par suite, lorsque le cou de la structure  presque-complexe $J$ est suffisamment allong\'e, les membranes $J$-holomorphes
qui survivent \`a cet allongement ont une structure conforme n'appartenant pas au compact $K$. $\square$

\begin{prop}
\label{propmin1} 
Soit $L$ une surface lagrangienne orientable hyperbolique dans une vari\'et\'e symplectique ferm\'ee de dimension quatre $(X, \omega)$ et  soit $d \in H_2 (X , L ; \Z)$. On note
$N_d^g (\underline{x} , J)$ le nombre de courbes $J$-holomorphes homologues \`a $d$ \`a bords dans $L$, de topologie et de structure conforme donn\'ees et qui passent par
une configuration $\underline{x}$ de points distincts de $(X, \omega)$ de cardinal ad\'equat, pour $J \in {\cal J}_\omega$ g\'en\'erique. Ce nombre
$N_d^g (\underline{x} , J)$ s'annule pour toute structure presque-complexe g\'en\'erale ayant un cou suffisamment long au voisinage de $L$.
\end{prop}

{\bf D\'emonstration de la Proposition \ref{propmin1} :}

On \'equipe \`a nouveau $L$ d'une m\'etrique \`a courbure constante et on allonge le cou d'une structure presque complexe g\'en\'erique au voisinage de $L$ jusqu'\`a briser la vari\'et\'e en deux,
ceci de mani\`ere \`a ce que les points de la configuration 
$\underline{x}$ disjoints de $L$ se retrouvent dans l'\'etage  $X \setminus L$. Notons $b+1$ le nombre de composantes connexes du bord des courbes que l'on consid\`ere. 
D'apr\`es le th\'eor\`eme de compacit\'e de th\'eorie symplectique des champs \cite{BEHWZ}, ces derni\`eres se brisent en courbes \`a deux \'etages qui sont cette fois-ci cod\'ees 
par des graphes $A_C$ ayant $b+1$ sommets marqu\'es $s_0, \dots , s_b$ correspondant aux 
$b+1$ composantes de bord. Les sommets \`a distances paires de $s_0, \dots , s_b$ codent \`a nouveau les composantes de l'\'etage $T^* L$ et les sommets \`a distances impaires
les composantes  de l'\'etage $X \setminus L$. Les orbites de Reeb du fibr\'e unitaire cotangent $S^* L$ sont cette fois-ci non-d\'eg\'en\'er\'ees, on fixe la trivialisation standard de  $S^* L$ le long 
de ces orbites de Reeb, de sorte que leur indice de Conley-Zehnder soit nul, voir la Proposition $1.7.3$ de \cite{EGH}. La dimension de l'espace des modules d'une composante simple $C_s$ 
de l'\'etage $X \setminus L$ est donn\'ee par le Th\'eor\`eme $2.8$ de \cite{HWZ}, elle vaut
$\mu_s^{CZ} + \chi_s$ o\`u $\mu_s^{CZ}$ est l'indice de Conley-Zehnder total de la composante et $\chi_s$ sa caract\'eristique d'Euler. Un rev\^etement ramifi\'e
d'une courbe simple ne peut en particulier qu'augmenter cette dimension puisque les indices de Conley-Zehnder des orbites de Reeb ont ici la propri\'{e}t\'{e} de s'additionner sous de tels rev\^etements. 
Par suite, l'in\'egalit\'e $\mu_s^{CZ} + \chi_s \geq 2n_s$ (resp. $\mu_s^{CZ} + \chi_s \geq n_s$) est satisfaite pour chaque sommet $s$ du graphe $A_C$ \`a distance impaire (resp. paire) de 
$s_0, \dots , s_b$, si $n_s \geq 0$ d\'esigne le nombre de points de la configuration $\underline{x}$ par lesquels passe
la composante cod\'ee par $s$. En sommant ces in\'egalit\'es sur tous les sommets du graphe $A_C$, on d\'eduit que la dimension totale attendue de la courbe $C$ se trouve
minor\'ee par $r + 2 r_X$ o\`u $r$ est le cardinal de $\underline{x} \cap L$ et $r_X$ le cardinal de $\underline{x} \setminus L$. Comme par hypoth\`ese cette dimension vaut
$r + 2 r_X$, les in\'egalit\'es pr\'ec\'edentes sont des \'egalit\'es. Il suit en particulier que toutes les courbes de l'\'etage $X \setminus L$ ont une dimension attendue paire ; elles ne peuvent 
par cons\'equent \^etre planes. Comme par ailleurs $L$ ne poss\`ede pas de g\'eod\'esique contractile, on vient de montrer que le graphe $A_C$ ne poss\`ede pas de feuille exception faite
\'eventuellement des sommets $s_0, \dots , s_b$. Remarquons \`a pr\'esent que 
chaque composante de la courbe \`a deux \'etages est asymptote \`a une r\'eunion de cylindres $J$-holomorphes sur les orbites de Reeb limites. Ces cylindres ont un module infini.
On en d\'eduit que les courbes compt\'ees par $N_d^g (\underline{x} , J)$, lorsque $J$ poss\`ede un cou extr\^emement long, poss\`edent \'egalement des anneaux de grands modules
dont les \^ames sont homotopes aux orbites de Reeb cod\'ees par les ar\^etes de l'arbre $A_C$. Or les courbes \`a bords compt\'ees par $N_d^g (\underline{x} , J)$ 
sont suppos\'ees avoir  une structure conforme fix\'ee. Tout anneau dont le module est sup\'erieur \`a une certaine quantit\'e donn\'ee par la structure conforme doit donc \^etre
contenu dans un disque. Par suite, lorsqu'on prive une telle courbe de la collection finie d'\^ames de nos anneaux de grands modules cod\'es par les ar\^etes de $A_C$, elle se trouve
disconnect\'ee en plusieurs composantes dont une au moins est un disque. Ce disque doit correspondre \`a une feuille du graphe $A_C$ distincte de $s_0, \dots , s_b$. Nous aboutissons
ainsi \`a une impossibilit\'e qui prouve que l'ensemble des courbes \`a deux \'etages sur lequel  nous avons fond\'e notre raisonnement est vide, ce qu'il fallait d\'emontrer.
$\square$

\begin{prop}
\label{propmintore}
Soit $(X, \omega , c_X)$ une vari\'et\'e symplectique r\'eelle ferm\'ee de dimension quatre dont le lieu r\'eel poss\`ede un tore lagrangien ou bien une surface hyperbolique lagrangienne $L$,
orientable ou non. On suppose que  $(X, \omega , c_X)$ ne poss\`ede pas de sph\`ere symplectique r\'eelle $S$ satisfaisant $\langle c_1 (X) , [S] \rangle > 1$ si $L$ est orientable et 
$\langle c_1 (X) , [S] \rangle > 0$ sinon. Soit $(d, g, b) \in H_2 (X , L ; \Z) \times \N \times \N^*$ 
et $K$ un compact de ${\cal M}_{g,b}$. Alors, pour toute structure presque-complexe g\'en\'erale ayant un cou suffisamment long au voisinage de $L$, la vari\'et\'e ne poss\`ede pas de membrane $J$-holomorphe homologue \`a $d$ \`a bord dans $L$ et conforme \`a un \'el\'ement de $K$.
\end{prop}

{\bf D\'emonstration de la Proposition \ref{propmintore} :}

On \'equipe \`a nouveau $L$ d'une m\'etrique \`a courbure constante et on allonge le cou d'une structure presque complexe g\'en\'erique au voisinage de $L$ jusqu'\`a briser la vari\'et\'e en deux 
morceaux. D'apr\`es le th\'eor\`eme de compacit\'e de th\'eorie symplectique des champs \cite{BEHWZ}, les membranes $J$-holomorphes homologues \`a $d$, de genre $g$ ayant $b$ composantes
de bord dans $L$ qui survivent \`a cette d\'eformation se brisent en courbes \`a deux \'etages cod\'ees par des graphes $A_C$ ayant $b$ sommets marqu\'es $s_1, \dots , s_b$ correspondant aux 
$b$ composantes de bord.  Les seules feuilles de ces arbres sont alors ces $b$ sommets $s_1, \dots , s_b$. En effet, ces feuilles coderaient sinon des  plans $J$-holomorphes asymptotes \`a
des orbites de Reeb du fibr\'e unitaire cotangent $S^* L$. Ces orbites de Reeb n'\'etant pas contractiles dans $T^* L$, les plans $J$-holomorphes doivent \^etre dans l'\'etage $X \setminus L$.
La r\'eunion d'un tel plan $P$, de son image par l'involution $c_X (P)$ et d'un cylindre $J$-holomorphe de $T^* L$ sur l'orbite de Reeb asymptote de $P$ fournit une sph\`ere \`a deux \'etages.
Cette sph\`ere se recolle en une sph\`ere symplectique $S$ de $(X, \omega)$ dont l'indice de Maslov vaut le double de l'obstruction \`a \'etendre la trivialisation canonique de $TX$ le long
de l'orbite de Reeb \`a $S$ tout entier. Lorsque $L$ est orientable, cette obstruction est nulle le long du cylindre de $T^* L$ et sup\'erieure \`a un le long de $P$ d'apr\`es la Proposition
\ref{propmaslov} et le Th\'eor\`eme $2.8$ de \cite{HWZ}. On en d\'eduit que l'indice de Maslov de $S$ serait sup\'erieur \`a quatre, ce qui contredit les hypoth\`eses. De la m\^eme mani\`ere lorsque $L$ 
est non-orientable, l'indice de Maslov de $S$ vaut la somme des indices de Conley-Zehnder de $P$, $c_X (P)$ et du cylindre. Ces derniers sont sup\'erieurs \`a leur caract\'eristique d'Euler 
puisqu'habitant des espaces
de modules de dimensions attendues positives, voir le Th\'eor\`eme $2.8$ de \cite{HWZ}. Par sommation, l'indice de Maslov de $S$ devrait \^etre sup\'erieur \`a deux ce qui contredit \`a nouveau
les hypoth\`eses. Les seules feuilles
des arbres codant les courbes \`a deux \'etages limites \'etant les $b$ sommets marqu\'es $s_1, \dots , s_b$, on d\'eduit comme dans la d\'emonstration de la Proposition  \ref{propmin1}
que peu avant la brisure de la vari\'et\'e, les membranes $J$-holomorphes homologues \`a $d$, de genre $g$ ayant $b$ composantes
de bord dans $L$ poss\`edent un anneau de grand module au moins, d'\^ame voisine d'une orbite de Reeb de $S^* L$. En particulier, elles n'appartiennent pas au compact $K$ de ${\cal M}_{g,b}$, 
ce qu'il fallait d\'emontrer. $\square$

\subsection{Formulaire}
\label{subsectformulaire}

D\'esignons par $L$ une vari\'et\'e compacte de dimension $n$ hom\'eomorphe \`a une sph\`ere, un tore ou un espace projectif r\'eel et munissons cette vari\'et\'e d'une m\'etrique \`a courbure constante.
Soit $C$ une courbe pseudo-holomorphe simple immerg\'ee d'\'energie de Hofer finie dans un remplissage symplectique du fibr\'e unitaire cotangent $(S^* L , \lambda)$ de $L$,  o\`u $\lambda$ 
d\'esigne la restriction de la forme de Liouville. 
La caract\'eristique d'Euler de $C$ vaut $\chi = 2-2g - v$ si l'on note $g$ son genre et $v$ le nombre de ses pointes.
Le flot de Reeb de $(S^* L , \lambda)$ trivialise le fibr\'e  normal de $C$ au voisinage de ses pointes ; nous noterons $\mu$ le double de l'obstruction \`a \'etendre cette trivialisation 
sur $C$ toute enti\`ere.

L'indice de Maslov $\mu$ est explicitement calcul\'e dans la Proposition \ref{propcotangent} lorsque $C$ est immerg\'ee dans le fibr\'e cotangent $T^* L$. Nous calculons sous cette m\^eme hypoth\`ese le 
nombre de points singuliers de $C$ dans le Lemme \ref{lemmepointsdoubles} tandis que la Proposition \ref{propmaslov} fournit une expression de la dimension de l'espace des d\'eformations de $C$
en fonction des quantit\'es $\chi$ et $\mu$. Les trois formules qui r\'esultent de ces calculs se montrant bien utiles, nous leurs consacrons ce paragraphe.

\begin{prop}
\label{propmaslov}
Soit $C$ une courbe pseudo-holomorphe simple 
d'\'energie de Hofer finie  dans un remplissage symplectique de dimension $2n$ du fibr\'e unitaire cotangent d'une sph\`ere, d'un tore ou d'un espace projectif r\'eel \`a courbure constante $L$.
Nous notons $g$ le genre de $C$, $v^-$ le nombre de ses pointes n\'egatives, $\mu$ son indice de Maslov et $\chi^- = 2-2g-v^-$. La dimension de l'espace des d\'eformations de $C$ vaut 
$\mu + (n-1)(2-2g)$ lorsque $L$ est une sph\`ere ou un espace projectif r\'eel et vaut $\mu + (n-1)\chi^-$ lorsque $L$ est un tore.
\end{prop}

{\bf D\'emonstration :}

Ces formules sont des cons\'equences de la formule d'indices  calcul\'ee par Fr\'ed\'eric Bourgeois dans sa th\`ese \cite{Bour}. La courbe $C$ converge en ses pointes vers
des orbites de Reeb qui appartiennent \`a des espaces de dimension $2(n-1)$ dans le premier cas et $n-1$ dans le second. Ces orbites contribuent donc \`a hauteur de
$2(n-1)v$ dans le premier cas et $(n-1)v$ dans le second \`a la dimension que l'on  calcule. Le reste de la contribution s'interpr\`ete comme l'indice de Fredholm de l'op\'erateur de 
Cauchy-Riemann associ\'e \`a $C$ et perturb\'e par un facteur $\mp \frac{d}{p} Id$ en ses pointes, ce qui le rend non-d\'eg\'en\'er\'e, voir la proposition $5.2$ de \cite{Bour}. 
Ce dernier est calcul\'e par le Th\'eor\`eme $2.8$ de \cite{HWZ} et vaut $\mu^{CZ} + (n-1)\chi$ o\`u $\mu^{CZ}$ d\'esigne l'indice de Conley-Zehnder normal total de $C$. L'indice de 
Conley-Zehnder normal total se d\'ecompose ici en la somme de l'indice de Maslov $\mu$ et des indices de Conley-Zehnder des op\'erateurs de Cauchy-Riemann perturb\'es en 
chaque pointe de $C$ et calcul\'es dans la trivialisation que l'on a fix\'e. Or ces indices de Conley-Zehnder des op\'erateurs de Cauchy-Riemann  perturb\'es
valent par d\'efinition $- (n-1)v$ dans le premier cas tandis qu'ils valent $0$ dans le second pour des pointes positives et $n-1$ pour des pointes n\'egatives. 
Ces r\'esultats sont \'etablis dans  le paragraphe $9.4$
de la th\`ese  \cite{Bour}. Signalons toutefois une d\'emonstration de ce dernier fait autre que celle propos\'ee par Fr\'ed\'eric Bourgeois. Lorsqu'on allonge la structure complexe de
$(\C P^1)^n$ au voisinage de $(\R P^1)^n$ jusqu'\`a briser la vari\'et\'e en deux morceaux, les fibres r\'eelles de $(\C P^1)^n \to (\C P^1)^{n-1}$ se brisent en un cylindre sur
un orbite simple de $T^* L$ et deux plans complexes conjugu\'es de $(\C P^1)^n \setminus L$. La dimension de l'espace des d\'eformations de chacun de ces morceaux vaut $n-1$ tandis que
les indices de Maslov de ces composantes sont tous nuls. Confrontons ce r\'esultat \`a ce qui pr\'ec\`ede.  La dimension de l'espace des d\'eformations du plan vaut 
$2(n-1)$ moins l'indice de l'op\'erateur de Cauchy-Riemann perturb\'{e} pour le plan, ce dernier vaut donc effectivement $n-1$. Elle vaut la moiti\'e de $2(n-1)$ moins le double de l'indice de 
l'op\'erateur perturb\'{e} pour le cylindre, puisque le cylindre poss\`ede deux pointes et se voit contraint d'\^etre pr\'eserv\'e par l'antipodation dans les fibres de $T^* L$. On en d\'eduit que l'indice 
de Conley-Zehnder des pointes  positives s'annule. $\square$

\begin{prop}
\label{propcotangent}
Soit $C$ une courbe pseudo-holomorphe simple d'\'energie de Hofer finie  dans le fibr\'e cotangent d'une sph\`ere, d'un tore ou d'un espace projectif r\'eel de dimension $n$
\`a courbure constante $L$.
Notons $\chi = 2-2g - v$ la caract\'eristique d'Euler de $C$ et $k$ la somme sur ses $v$ pointes des multiplicit\'es de ses orbites de Reeb limites. L'indice de Maslov $\mu$ de $C$
vaut $2(n-1)k - 2\chi$ lorsque $L$ est une sph\`ere, $(n-1)k - 2\chi$ lorsque $L$ est un espace projectif r\'eel et $- 2\chi$ lorsque $L$ est un tore.
\end{prop}

{\bf D\'emonstration :}

Consid\'erons le deux-cycle $C - c_L (C) - \sum_{i=1}^v \text{Cyl}_i$, o\`u $c_L$ est l'antipodation dans les fibres de $T^* L$ et $\text{Cyl}_i$ les cylindres sur les orbites de Reeb
limites de $C$. Ce deux-cycle se trouve renvers\'e par $c_L$ de sorte qu'il est homologue \`a z\'ero. La premi\`ere classe de Chern de $T^* L$ s'annule donc une fois \'evalu\'ee contre
ce cycle. Calculons cette derni\`ere comme l'obstruction \`a trivialiser le fibr\'e tangent en restriction \`a ce deux-cycle. La contribution de $C - c_L (C)$ vaut $2\chi + \mu$ o\`u $\chi$. La contribution 
d'un cylindre $\text{Cyl}_i$ vaut l'oppos\'e de son demi-indice de Maslov, soit $-k_i$ fois le demi-indice de Maslov du cylindre sur l'orbite simple
sous-jacente si $k_i$ d\'esigne la multiplicit\'e de l'orbite. 

Le demi-indice de Maslov d'un cylindre sur une orbite simple dans le cas d'une sph\`ere vaut le degr\'e du fibr\'e
normal d'une section plane de la quadrique ellipso\"{\i}de $Q^n$ puisque cette derni\`ere est obtenue en recollant deux plans de fibr\'es normaux triviaux de $Q^n \setminus L$
au cylindre en question. Ce dernier vaut donc $2n - 2$, d'o\`u la relation $2\chi + \mu - 2\sum_{i=1}^v (n-1)k_i = 0$. 

Le demi-indice de Maslov d'un cylindre sur une orbite simple dans le cas d'un espace projectif r\'eel vaut le degr\'e du fibr\'e
normal d'une droite dans l'espace projectif complexe de dimension $n$ puisque cette derni\`ere est obtenue en recollant deux plans de fibr\'es normaux triviaux de $\C P^n \setminus L$
au cylindre en question. Ce dernier vaut donc $n-1$, d'o\`u la relation $2\chi + \mu - \sum_{i=1}^v (n-1)k_i = 0$. 

Le demi-indice de Maslov d'un cylindre sur une orbite simple dans le cas d'un tore est trivial, d'o\`u la relation $2\chi + \mu  = 0$. $\square$

\begin{lemme}
\label{lemmepointsdoubles}
Soit $C$ une courbe pseudo-holomorphe d'\'energie de Hofer finie immerg\'ee dans le fibr\'e cotangent d'une sph\`ere de dimension deux ou d'un plan projectif r\'eel \`a courbure constante. 
Supposons que cette courbe  soit simple, rationnelle, r\'eelle et n'ayant que des points doubles transverses comme singularit\'es.
On note $v$ le nombre de paires de pointes complexes conjugu\'ees de $C$ et $k$ la multiplicit\'e totale des paires d'orbites de Reeb limites en ces pointes. Le nombre de points doubles de $C$ 
est major\'e par $k^2 - 2k + 1$  dans le cas d'une sph\`ere et $\frac{1}{2}(k^2 - 3k + 2)$ dans le cas d'un plan projectif r\'eel.
\end{lemme}

{\bf D\'emonstration :}

Consid\'erons le deux-cycle $C - \sum_{i=1}^v \text{Cyl}_i$, o\`u $\text{Cyl}_i$, $1 \leq i \leq v$, d\'esignent les cylindres sur les orbites de Reeb
limites de $C$. Ce deux-cycle se trouve renvers\'e par $c_L$ de sorte qu'il est homologue \`a z\'ero. Choisissons un cylindre $\text{Cyl}$ sur une orbite de Reeb distincte des limites de $C$.
L'indice d'intersection de $\text{Cyl}$ avec $C - \sum_{i=1}^v \text{Cyl}_i$ s'annule. Nous en d\'eduisons que l'indice d'intersection de $\text{Cyl}$ avec $C$ vaut $2k$ dans le cas de la sph\`ere de dimension
deux et $k$ dans celui du plan projectif r\'eel. Perturbons \`a pr\'esent $C$ en une courbe voisine $\widetilde{C}$ dont toutes les orbites de Reeb limites sont distinctes de celles de $C$.
L'indice d'intersection de $\widetilde{C}$ avec $C - \sum_{i=1}^v \text{Cyl}_i$ s'annule. On d\'eduit de ce qui pr\'ec\`ede que l'indice d'intersection de $\widetilde{C}$ avec $C$
se trouve major\'e par $2k^2$ dans le cas de la sph\`ere de dimension deux et $k^2$ dans celui du plan projectif r\'eel. Cet indice est par ailleurs minor\'e par deux fois le nombre de points doubles 
de $C$ auquel s'ajoute la moiti\'e de son
indice de Maslov  et le nombre de points d'intersection de $\widetilde{C}$ avec $C$ qui apparaissent au voisinage des pointes de $C$. Ces derniers sont au moins au nombre de $k_i -1$
au voisinage de chaque pointe convergeant vers une orbite de Reeb parcourue $k_i$ fois, ce qui d\'ecoule du Th\'eor\`eme $1.5$ de \cite{HWZ1}, soit $2(k-v)$ au total. L'indice de Maslov de $C$ 
est quant \`a lui donn\'e par la Proposition \ref{propcotangent}, il vaut $4k + 4v - 4$ dans le cas de la sph\`ere de 
dimension deux et $2k + 4v - 4$ dans celui du plan projectif r\'eel. Ainsi, l'indice d'intersection de $\widetilde{C}$ avec $C$ se trouve minor\'e par $4k-2$ (resp. $3k - 2$) plus deux fois le nombre 
de points doubles de $C$ si $L$ est une sph\`ere (resp. un plan projectif r\'eel). Le r\'esultat en d\'ecoule. $\square$

\begin{rem}
\label{reminters}
Nous avons \'etabli au cours de la d\'emonstration du Lemme \ref{lemmepointsdoubles} la majoration $\vert \widetilde{C} \circ C \vert \leq 2k^2$ ou $k^2$ selon que $L$ est une sph\`ere ou un plan 
projectif r\'eel. Cette majoration nous sera utile au \S \ref{sectcalculs}.
\end{rem}

\section{Congruences}
\label{sectcong}

\subsection{\'Enonc\'es des r\'esultats}

\'Etant donn\'ee une classe d'homologie $d \in H_2 (X ; \Z)$ d'une vari\'et\'e symplectique r\'eelle de dimension quatre $(X, \omega , c_X)$, nous noterons 
$g_d = \frac{1}{2} (d^2 - c_1 (X)d + 2)$ le genre lisse de $d$ et $c_d = c_1 (X)d -1$ le degr\'e attendu du polyn\^ome $\chi^d (T)$ d\'efini dans \cite{Wels1}.

\begin{theo}
\label{theocong1}
Soit $(X, \omega , c_X)$ une vari\'et\'e symplectique r\'eelle ferm\'ee de dimension quatre dont le lieu r\'eel poss\`ede une composante connexe $L$ hom\'eomorphe \`a une sph\`ere.
Soient $d \in H_2 (X ; \Z)$ et $r \in \N$. Lorsque $2r+1 < c_d$, la puissance $2^{\frac{1}{2} (c_d - 2r - 1)}$ divise $\chi^d_r (L)$. 

Supposons en outre la connexit\'e du lieu r\'eel de la vari\'et\'e $(X, \omega , c_X)$. Alors,

a)  Lorsque $2r-1 < c_d$ et lorsque de plus $g_d$ et $\frac{1}{2} (r+1)$ sont de m\^eme parit\'e, la puissance $2^{\frac{1}{2} (c_d - 2r + 1)}$ divise $\chi^d_r (L)$. 

b) Lorsque $2k < r+1 \leq \frac{1}{2} c_1 (X)d + 2$, la puissance $2^{\frac{1}{2} (c_d - 2r + 3)}$ divise $\chi^d_r (L)$, o\`u $k$ d\'esigne le maximum de l'ensemble
$\{ j \in  \N \, \vert \, j \neq g_d \mod(2) \text{ et } j \leq \vert d' \circ [L] \vert \text{ o\`u } $d'$ \text{ est effectif satisfaisant }\\ d' - c_X (d') = d  \}.$
\end{theo}
On entend ici par classe effective une classe d'homologie r\'ealisable par un deux-cycle pseudo-holomorphe sur son deux-squelette.\\

{\bf Exemple :}

Le Th\'eor\`eme \ref{theocong1} s'applique \`a l'ellipso\"{\i}de de dimension deux lorsque $d$ est un multiple positif, disons $\delta > 0$, d'une section plane r\'eelle. 
Dans ce cas, $c_d = 4\delta - 1$ et $g_d = \delta^2 - 2 \delta + 1 = \delta + 1 \mod (2)$. Par cons\'equent,  $2^{2 \delta - r - 1}$ divise $\chi^d_r (L)$ lorsque $r < 2 \delta - 1$,
$2^{2 \delta - r}$ divise $\chi^d_r (L)$ lorsque de plus $r = 2 \delta + 1 \mod (4)$ et $\chi^d_{2 \delta - 3} (L) = 0 \mod (16)$.

\begin{theo}
\label{theocong2}
Soient $(X, c_X)$ la quadrique ellipso\"{\i}de de dimension trois et $d$ un multiple positif, disons $\delta > 0$, d'une section hyperplane r\'eelle.
Lorsque $6r + 1 \leq 3 \delta$, la puissance $2^{\frac{3}{4} (\delta - 2r)}$ divise $\chi^d_r$. 
\end{theo}

\begin{theo}
\label{theocong3}
Soit $(X, \omega , c_X)$ une vari\'et\'e symplectomorphe au plan projectif complexe \'eclat\'e en six boules complexes conjugu\'ees au maximum. Soit $d \in H_2 (X ; \Z)$ 
une classe satisfaisant $c_d = c_1 (X)d -1 \geq 0$ et soient $r, r_X$ des entiers naturels satisfaisant la relation $r +  2r_X = c_d$. Lorsque $r+ 1 < r_X$, la puissance
$2^{r_X - r - 1}$ divise $\chi^d_r (L)$. Lorsque $r < r_X$ et lorsque de plus $r = \langle d , h \rangle + 1 \mod (4)$, o\`u $h$ est la classe d'une droite g\'en\'erique du plan, 
la puissance $2^{r_X - r }$ divise $\chi^d_r (L)$. Lorsqu'enfin 
la vari\'et\'e est le plan projectif complexe lui-m\^eme et $r + 1 < \langle d , h \rangle$, $\chi^d_{r} (L) = 0 \mod (64)$.
\end{theo}

{\bf Exemple :}

Le Th\'eor\`eme \ref{theocong3} s'applique au plan projectif complexe o\`u $d$ est un multiple positif, disons $\delta > 0$, d'une droite complexe. 
Dans ce cas, $8^{\frac{1}{2}(\delta  - r - 1)}$ divise $\chi^d_r$ lorsque $r+1 < \delta $, $2^{\frac{1}{2}(3\delta  - 3r - 1)}$ divise $\chi^d_r$ lorsque de plus $r = \delta  + 1 \mod (4)$
et $\chi^d_{\delta - 3} = 0 \mod (64)$.

\subsection{D\'emonstrations des Th\'eor\`emes \ref{theocong1},  \ref{theocong2} et \ref{theocong3}}

On suit la strat\'egie g\'en\'erale \'enonc\'ee au paragraphe \ref{subsubsectstrat} en \'equipant la sph\`ere ou le plan projectif r\'eel lagrangien d'une m\'etrique \`a courbure constante. 
Les courbes $J$-holomorphes rationnelles r\'eelles $C$ compt\'ees par l'invariant $\chi^d_r (L)$ qui survivent  \`a l'allongement du cou de $J$ jusqu'\`a la brisure de
la vari\'et\'e sont cod\'ees par des arbres $A_C$, voir la figure \ref{figarbre}. Ces derniers ont une racine $s_0$ qui code l'unique composante de la courbe \`a deux \'etages limite 
laiss\'ee invariante par l'involution $c_X$.
Le fibr\'e normal de chaque composante simple $C_s$ associ\'e \`a un sommet $s$
d'un arbre $A_C$ est canoniquement trivialis\'e le long des orbites de Reeb asymptotes par le flot de Reeb et l'on note $\mu_s$ l'obstruction \`a \'etendre cette trivialisation sur $C_s$ 
tout entier.
La dimension de l'espace des modules dans lequel habite $C_s$ s'exprime par la relation $\mu_s + 2(\text{dim}_\C X - 1)$ lorsque $s \neq s_0$
et $\frac{1}{2} \mu_s + \text{dim}_\C X - 1$ lorsque $s = s_0$ puisque la courbe $C_s$ est alors contrainte d'\^etre pr\'eserv\'ee par l'involution $c_X$, 
voir la Proposition \ref{propmaslov} de notre formulaire donn\'e au paragraphe \ref{subsectformulaire}. Notons
$s_0 , \dots , s_j$ les sommets de $A_C$ qui codent les composantes $C_s$ soumises \`a des conditions d'incidence et $\widetilde{A}_C$ le sous arbre de $A_C$ obtenu en ne 
retenant que les sommets $s_0 , \dots , s_j$ et les ar\^etes reliant ces sommets entre eux. \\

{\bf D\'emonstration du Th\'eor\`eme \ref{theocong1} :}

Nous allons commencer par minorer la contribution de chaque composante connexe
de $A_C \setminus \widetilde{A}_C$ \`a l'indice de Maslov total de la courbe $C$.
Lorsqu'une telle composante n'est pas connect\'ee \`a $s_0$, cette contribution est minor\'ee par $2(k' - 1)$ o\`u $k'$ d\'esigne la multiplicit\'e totale des ar\^etes reliant cette composante
connexe \`a $\widetilde{A}_C$. Ceci r\'esulte de l'in\'egalit\'e (\ref{equ1}) \'etablie au \S \ref{subsectoptdem}. Lorsqu'une telle composante est connect\'ee \`a $s_0$, la Proposition
\ref{propcotangent} fournit l'estimation $\sum_{s \in S_2'} \mu_s = 2\sum_{s \in S_2'} (k_s + v_s) - 4 \# S_2'$ de la contribution des sommets \`a distance paire de $s_0$, o\`u
$S_1'$ (resp. $S_2'$) d\'esigne l'ensemble des sommets \`a distance impaire (resp. paire) de $s_0$ de cette composante. 
Cette composante poss\`ede un unique sommet $s$ ayant la propri\'et\'e d'\^etre connect\'e \`a $s_0$. Notons $l^0$ le degr\'e du rev\^etement de la courbe $C_s$ cod\'ee par ce sommet et
$k^0$ la multiplicit\'e de l'ar\^ete qui le joint \`a $s_0$.
La minoration (\ref{equ3}) \'etablie au \S \ref{subsectoptdem} fournit, en reprenant les notations introduites dans
ce paragraphe, $\sum_{s \in S_1'} \mu_s \geq 2 \sum_{s \in S_1'} ( l_s  - l_s \tilde{v}_s  + v_s) - 4 \# S_1'$, soit $\sum_{s \in S_1'} \mu_s \geq 2 \sum_{s \in S_1'} ( l_s  - k_s  + v_s) - 
4 \# S_1' +2k^0 - 2l^0$. Nous en d\'eduisons 
apr\`es sommation $\sum_{s \in S_1'  \cup S_2'} \mu_s   \geq  2 \sum_{s \in S_1'} l_s - 2l^0 + 2k' +2v' - 2$, o\`u $v'$ (resp. $k'$) d\'esigne le nombre
d'ar\^etes (resp. leur multiplicit\'e totale) reliant cette composante \`a un sommet $s_1 , \dots , s_j$. L'indice de Maslov d'une telle composante se trouve donc finalement minor\'e par $2v'$ 
except\'e dans le cas
o\`u $v'$ est nul et cet indice vaut $-2$. Notons $c_{-2}$ le nombre de composantes de $A_C \setminus \widetilde{A}_C$ d'indice de Maslov total $-2$. De ces calculs r\'esulte que la contribution
de $A_C \setminus \widetilde{A}_C$  \`a l'indice de Maslov total de la courbe \`a deux \'etages cod\'ee par $A_C$ est minor\'ee par $2a - 2c_{-2}$ si $a+1$ d\'esigne le nombre de composantes 
connexes de $\widetilde{A}_C$.

La contribution des sommets $s_1 , \dots , s_j$ est quant \`a elle minor\'ee par $2 r_X - 2j$. La courbe r\'eelle que code le sommet $s_0$ se voit d'une part contrainte d'interpoler $r$ points 
de $L$ et d'autre part de converger en $c_{-2}$ paires complexes conjugu\'ees de ses pointes vers $c_{-2}$ paires complexes conjugu\'ees d'orbites de Reeb prescrites. En effet, 
les $c_{-2}$ sommets correspondant de $A_C \setminus \widetilde{A}_C$ adjacents \`a $s_0$ codent des courbes rigides. Enfin, chaque sommet $s_1 , \dots , s_j$ connect\'e \`a $s_0$ 
pr\'esente l'alternative suivante. Soit la minoration pr\'ec\'edente $2 r_X - 2j$ est atteinte pour ce sommet et la courbe correspondante, avec ses conditions 
d'incidences, est rigide ; ce qui ajoute donc une contrainte suppl\'ementaire pour une paire de pointes de $s_0$. Soit la minoration pr\'ec\'edente $2 r_X - 2j$ n'est pas atteinte pour ce sommet
et peut donc \^etre am\'elior\'ee de deux. Nous aboutissons dans tous les cas \`a la minoration $\sum_{s \in \widetilde{A}_C} \mu_s \geq 2 r_X - 2a + r - 1 +2c_{-2} $.  L'indice de Maslov total $\mu$ 
de la courbe \`a deux \'etages cod\'ee par $A_C$ se trouve ainsi minor\'e par $r + 2 r_X - 1 = c_d - 1$. Comme cet indice est par ailleurs major\'e par cette quantit\'e qui n'est autre que le degr\'e du 
fibr\'e normal d'une courbe rationnelle irr\'eductible immerg\'ee homologue \`a $d$, toutes nos in\'egalit\'es doivent \^etre \'egalit\'es. Nous en concluons que les composantes connexes
de $A_C \setminus \widetilde{A}_C$ qui ne sont pas connect\'ees \`a $s_0$ sont r\'eduites \`a un sommet codant un plan asymptote \`a une orbite de Reeb simple tandis que les composantes
connect\'ees \`a $s_0$  sont au nombre de $c_{-2}$ et leur indice de Maslov vaut $-2$. L'arbre $\widetilde{A}_C$ est en particulier connexe.

L'arbre $A_C$ vient avec une donn\'ee combinatoire suppl\'ementaire, une fonction qui associe \`a chaque sommet \`a distance impaire de $s_0$ les classes d'homologies relatives de la paire 
de courbes correspondantes cod\'ee par ce sommet ainsi que les paires de points complexes conjugu\'es de $\underline{x} $ que ces courbes contiennent. Notons $r_1 , \dots , r_j$ le nombre
de paires de points complexes conjugu\'es de $\underline{x} $ associ\'ees \`a $s_1 , \dots , s_j$ respectivement, de sorte que leur somme vaille $r_X$. Il y a $2^{r_i - 1}$ partitions d'un ensemble
de $r_i$ points complexes conjugu\'es en deux ensembles complexes conjugu\'es, soit ici $2^{r_X - j}$ partitions au total. Une fois attribu\'es \`a chaque courbe $C_s$ l'ensemble de points qu'elle
doit interpoler, certaines de ces courbes sont rigides et d'autres non. Notons $j^-$ le nombre de telles courbes rigides et $j^+ = j - j^-$. D'apr\`es ce qui pr\'ec\`ede, la courbe r\'eelle cod\'ee par $s_0$,
avec ses $j^- + c_{-2}$ paires d'asymptotes prescrites et ses $r$ points r\'eels \`a interpoler, est rigide. La dimension $2k_{s_0} + 2v_{s_0} -1$ donn\'ee par la Proposition \ref{propcotangent}
vaut donc en particulier $r + 2j^- + 2c_{-2}$. Par cons\'equent, les $j^+$ paires de courbes non-rigides pr\'ec\'edentes h\'eritent d'une contrainte suppl\'ementaire, elles ont une paire
d'asymptotes prescrites correspondant \`a une paire d'orbites de Reeb limites rest\'ees libres de la courbe $C_{s_0}$. Il y a deux bijections possibles entre une telle paire d'orbites de Reeb et
une telle paire de courbes non-rigides, soit $2^{j^+}$ bijections au total. Ainsi, le nombre de courbes \`a deux \'etages ayant une combinatoire donn\'ee par $A_C$ est divisible par $2^{r_X - j^-}$.
Or d'apr\`es le th\'eor\`eme de recollement en th\'eorie symplectique des champs \cite{Bour} et le Th\'eor\`eme \ref{theoinv}, la contribution \`a l'invariant $\chi^d_{r} (L)$ d'une courbe \`a deux 
\'etages ne d\'epend que de sa combinatoire, de sorte que $2^{r_X - j^-}$ divise $\chi^d_{r} (L)$. L'\'equation 
\begin{equation}
\label{relj1}
r + 2j^- + 2c_{-2} = 2k_{s_0} + 2v_{s_0} -1
\end{equation}
impose l'in\'egalit\'e $r+1 \geq 2k_{s_0}$.
On en d\'eduit $2j^- \leq 2v_{s_0} \leq 2k_{s_0} \leq r+1$ et le premier r\'esultat \'enonc\'e dans le Th\'eor\`eme \ref{theocong1}. 

Tous les arbres $A_C$ pour lesquels l'une des in\'egalit\'es $2j^- \leq 2v_{s_0} \leq 2k_{s_0} \leq r+1$ est stricte satisfont $2j^- \leq r-1$ et le deuxi\`eme \'enonc\'e du Th\'eor\`eme \ref{theocong1}
est imm\'ediat. Soit $A_C$ un arbre pour lequel $j^- = v_{s_0} = k_{s_0} = \frac{1}{2} (r+1)$. Toutes les courbes cod\'ees par les sommets de cet arbre sont simples puisque d'indices de Maslov 
positifs et ont des orbite de Reeb simples pour asymptotes. Consid\'erons la courbe $C'$ form\'ee de toutes les paires de courbes de l'\'etage $X \setminus L$
cod\'ees par $A_C$ et de paires de plans complexes conjugu\'es de $T^* L$ convergeant vers $\partial C'$. Une telle courbe \`a deux \'etages
se recolle en une courbe $J$-holomorphe r\'eductible homologue \`a $d$ ayant un nombre pair de composantes irr\'eductibles \'echang\'ees par $c_X$. Le nombre de points doubles d'une telle
courbe a la parit\'e de  $g_d + 1$, d'apr\`es la formule d'adjonction. Supposons le lieu r\'eel de $X$ connexe, ce nombre de points double est alors \'egalement de la m\^eme parit\'e que le nombre
de points d'intersection avec $L$, c'est-\`a-dire que le nombre de paires de plans de $T^* L$ que l'on a introduit ou encore le nombre total d'ar\^etes $v$ de l'arbre $A_C$. Or par hypoth\`ese,
$g_d$ et $\frac{1}{2} (r+1)$ sont de m\^eme parit\'e et $\frac{1}{2} (r+1) = v_{s_0}$. Par cons\'equent $v$ et $v_{s_0}$ ne sont pas de m\^eme parit\'e ce qui impose l'existence d'un sommet
de valence paire adjacent \`a $s_0$. La structure presque complexe de $X \setminus L$ \'etant g\'en\'erale, toutes les pointes de la courbe cod\'ee par ce sommet ont des asymptotes distinctes et
cette derni\`ere est rigide puisque $j^- = v_{s_0}$. Le r\'esultat $a)$ d\'ecoule \`a pr\'esent du fait qu'il y a un nombre pair de choix de la pointe de cette courbe \`a relier \`a $C_{s_0}$, ce qui
permet d'am\'eliorer d'une puissance de deux la divisibilit\'e du nombre de courbes cod\'ees par ces arbres $A_C$.

Enfin, tous les arbres $A_C$ pour lesquels $2k_{s_0} <  r+1$ ou $v_{s_0} + 1 < k_{s_0}$ satisfont $2j^- \leq r-3$, ce qui d\'ecoule de (\ref{relj1}). Or si $A_C$ est un arbre tel que $j^- = v_{s_0} = k_{s_0}$,
le raisonnement
que l'on vient de suivre fournit une paire de courbes $J$-holomorphes r\'eductibles complexes conjugu\'ees homologue \`a $d$. Notons $d'$ la classe d'homologie d'une telle courbe de sorte que
$d' - c_X (d') = d$. Cette courbe peut \^etre choisie de sorte que l'indice d'intersection $\vert d' \circ [L] \vert $ vaille la multiplicit\'e totale $k$ des ar\^etes de $A_C$. Comme $k$ et $g_d$ ne sont pas
de m\^eme parit\'e, $k$ fait partie de l'ensemble d\'efini en $b)$. L'hypoth\`ese implique \`a pr\'esent $2k_{s_0} \leq 2k < r+1$, d'o\`u le r\'esultat $b)$ dans ce cas. Si en revanche 
$j^- = v_{s_0} = k_{s_0} -1$, l'in\'egalit\'e $2j^- < r+1$ permet d'am\'eliorer d'une puissance de deux le premier \'enonc\'e du Th\'eor\`eme \ref{theocong1}. En outre, une ar\^ete adjacente \`a $s_0$
est de multiplicit\'e deux, de sorte qu'une courbe $C_s$ se trouve connect\'ee \`a $C_{s_0}$ par une orbite de Reeb double. Le th\'eor\`eme de recollement en th\'eorie symplectique des champs \cite{Bour}
garantit alors l'existence d'un nombre pair de courbes $J$-holomorphes convergeant vers une courbe \`a deux \'etages donn\'ee cod\'ee par $A_C$. Cette parit\'e provenant du param\`etre de recollement associ\'e \`a l'orbite double permet d'am\'eliorer le premier \'enonc\'e du Th\'eor\`eme \ref{theocong1} d'une puissance de deux suppl\'ementaire. D'o\`u le r\'esultat. $\square$\\

{\bf D\'emonstration du Th\'eor\`eme \ref{theocong2} :}

Le compl\'ementaire $X \setminus L$ est isomorphe au fibr\'e en droites de bidegr\'e $(1,1)$ sur la quadrique $\C P^1 \times \C P^1$. Les courbes rationnelles irr\'eductibles d'\'energie de Hofer finie
de ce compl\'ementaire sont donc d'indice de Maslov positif.
En effet, ce compl\'ementaire se compactifie en le fibr\'e en droites projectives $F$ obtenu \`a partir de la somme du fibr\'e trivial et du fibr\'e de bidegr\'e $(1,1)$ sur $\C P^1 \times \C P^1$. 
Une courbe d'\'energie de Hofer finie se compactifie en une courbe de $F$ dont la classe d'homologie s'\'ecrit $e + kf$, o\`u $e$ d\'esigne une classe d'homologie effective de la section nulle 
$\C P^1 \times \C P^1$ du fibr\'e et $f$ la classe d'une fibre. Or il suit de la formule d'adjonction que l'\'evaluation de la premi\`ere classe de Chern de $F$ sur $e$ est positive et vaut deux sur $f$.
Les indices de Maslov de ces courbes sont donc positifs et ne peuvent qu'augmenter par rev\^etements ramifi\'es. Il s'ensuit que toutes les courbes $C_s$ cod\'ees par les sommets $s$ 
de l'arbre sont des courbes simples pour peu que la configuration de points choisie soit suffisamment g\'en\'erale.
La dimension $\mu_s + 4$ des espaces de modules associ\'es est donc strictement positive et m\^eme sup\'erieure au quadruple (resp. au double) du nombre de points de la configuration que doit
interpoler $C_s$ si $s$ est \`a distance impaire (resp. paire) de $s_0$.
Par suite, la contribution totale \`a l'indice de Maslov des sommets $s \in S_1$ de $A_C$ \`a distance impaire de $s_0$, lesquels codent les courbes de $X  \setminus L$, se trouve minor\'ee
par $4 r_X - 4 \# S_1$. De m\^eme, la contribution totale \`a l'indice de Maslov des sommets $s \in S_2$ de $A_C$ \`a distance paire de $s_0$, lesquels codent les courbes de $T^* L$, 
se trouve minor\'ee par $2r - 4 \# S_2 + 2$ puisque $s_0$ code une courbe r\'eelle ayant $\frac{1}{2} \mu_{s_0} + 2$ degr\'es de libert\'e. En outre, chaque ar\^ete de l'arbre $A_C$ code 
une paire de pointes des courbes cod\'ees par les sommets adjacents et qui ont une m\^eme orbite de Reeb pour asymptote. 
Cette contrainte co\^ute quatre degr\'es de libert\'e suppl\'ementaires par ar\^ete, soit au total $4 \# S_1 + 4 \# S_2 - 4$ degr\'es, puisque le nombre d'ar\^etes d'un arbre diff\`ere du nombre de sommets
par un. Ceci nous permet finalement de minorer l'indice de Maslov de notre courbe \`a deux \'etages cod\'ee par $A_C$ par $2r + 4r_X - 2 = \langle c_1 (X) , d \rangle -2$. Cet indice est par ailleurs 
major\'e par le degr\'e $\langle c_1 (X) , d \rangle -2$
du fibr\'e normal de $C$, de sorte que toutes nos minorations sont des \'egalit\'es. En particulier, tous les sommets adjacents \`a $s_0$ sont des feuilles qui doivent coder des courbes interpolant
chacune au moins un point de la configuration, d'apr\`es la Proposition \ref{propcotangent}. 

Notons \`a pr\'esent $r_1 , \dots , r_j$ le nombre
de paires de points complexes conjugu\'es de $\underline{x} $ associ\'ees aux sommets $s_1 , \dots , s_j$ adjacents \`a $s_0$ respectivement, de sorte que leur somme vaille $r_X$. 
Il y a $2^{r_i - 1}$ partitions d'un ensemble
de $r_i$ points complexes conjugu\'es en deux ensembles complexes conjugu\'es, soit ici $2^{r_X - j}$ partitions au total. Une fois attribu\'es \`a chaque courbe $C_s$ l'ensemble de points qu'elle
doit interpoler, certaines de ces courbes sont rigides, d'autres conservent deux ou quatre degr\'es de libert\'e. Notons $j^-$ le nombre de telles courbes rigides, $j_1^+$ (resp. $j_2^+$) le nombre de
celles qui conservent deux (resp. quatre) degr\'es de libert\'e. D'apr\`es ce qui pr\'ec\`ede, la courbe r\'eelle cod\'ee par $s_0$,
avec ses $j^- + j_1^+$ paires d'asymptotes prescrites et ses $r$ points r\'eels \`a interpoler, est rigide. La dimension $4k_{s_0} + 2v_{s_0}$ donn\'ee par la Proposition \ref{propcotangent}
vaut donc en particulier $2r + 4j^- + 2j_1^+$. Par cons\'equent, les $j_2^+$ paires de courbes non-rigides pr\'ec\'edentes h\'eritent d'une contrainte suppl\'ementaire, elles doivent converger
en une paire de pointes complexes conjugu\'ees vers une paire prescrite d'orbites de Reeb limites de la courbe $C_{s_0}$. Il y a deux bijections possibles entre une telle paire d'orbites de Reeb et
une telle paire de courbes non-rigides, soit $2^{j_2^+}$ bijections au total. Ainsi, le nombre de courbes \`a deux \'etages ayant une combinatoire donn\'ee par $A_C$ est divisible par $2^{r_X - j^- - j_1^+}$.
Or d'apr\`es le th\'eor\`eme de recollement en th\'eorie symplectique des champs \cite{Bour} et le Th\'eor\`eme \ref{theoinv}, la contribution \`a l'invariant $\chi^d_{r}$ d'une courbe \`a deux 
\'etages ne d\'epend
que de sa combinatoire, de sorte que $2^{r_X - j^- - j_1^+}$ divise $\chi^d_{r}$. L'\'equation $2r + 4j^- + 2j_1^+ = 4k_{s_0} + 2v_{s_0}$ impose l'in\'egalit\'e $r \geq k_{s_0}$.
On en d\'eduit $j^- + j_1^+ \leq v_{s_0} \leq r$ de sorte que $2^{r_X - r}$ divise $\chi^d_{r}$. Or par hypoth\`ese, $2r + 4r_X = 3\delta$ puisque la premi\`ere classe de Chern
de la quadrique de dimension trois est Poincar\'e duale au triple de la section hyperplane. Le Th\'eor\`eme \ref{theocong2} en d\'ecoule. $\square$\\

{\bf D\'emonstration du Th\'eor\`eme \ref{theocong3} :}

D'apr\`es les hypoth\`eses que l'on a faites, le compl\'ementaire $X \setminus L$ est isomorphe au fibr\'e en droites de degr\'e quatre sur $\C P^1$ \'eclat\'e en six points complexes conjugu\'es
au maximum. Ce compl\'ementaire ne contient par cons\'equent pas de courbes rationnelles irr\'eductibles d'\'energie de Hofer finie et d'indice de Maslov n\'egatif ayant plus de deux pointes.
En effet, ce compl\'ementaire se compactifie en la surface r\'egl\'ee rationnelle de degr\'e quatre $\Sigma_4$ \'eclat\'ee en six points complexes conjugu\'es
au maximum. Une courbe d'\'energie de Hofer finie se compactifie en une courbe dont la classe d'homologie s'\'ecrit $ne + kf - \sum \alpha_i E_i$, o\`u $e$ d\'esigne la section nulle du fibr\'e, 
$f$ une fibre et $E_i$ les diviseurs exceptionnels des \'eclatements. L'irr\'eductibilit\'e de la courbe
force $0 \leq \alpha_i \leq n$ d\`es que $n \geq 1$ et la premi\`ere classe de Chern de $\Sigma_4$ est duale \`a $2e - 2f- \sum E_i$, de sorte que son \'evaluation $6n + 2k - \sum \alpha_i$ 
sur la courbe soit minor\'ee par $2k$ lorsque $n \geq 1$. L'indice de Maslov de telles courbes $C_s$ est donc positif, puisque minor\'e par $2k - 2 \chi (C_s)$. L'absence de courbes d'indice 
de Maslov $-2$ autres que planes a pour cons\'equence que pour tout
arbre $A_C$, les courbes $C_s$ cod\'ees par les sommets $s$ de l'arbre sont des courbes simples pour peu que la configuration de points choisie soit suffisamment g\'en\'erale.
La dimension $\mu_s + 2$ des espaces de modules associ\'es est donc positive et m\^eme sup\'erieure au double du nombre (resp. au nombre) de points de la configuration que doit
interpoler $C_s$ si $s$ est \`a distance impaire (resp. paire) de $s_0$.
Par suite, la contribution totale \`a l'indice de Maslov des sommets $s \in S_1$ de $A_C$ \`a distance impaire de $s_0$, lesquels codent les courbes de $X  \setminus L$, se trouve minor\'ee
par $2 r_X - 2 \# S_1$. De m\^eme, la contribution totale \`a l'indice de Maslov des sommets $s \in S_2$ de $A_C$ \`a distance paire de $s_0$, lesquels codent les courbes de $T^* L$, se trouve minor\'ee
par $r - 2 \# S_2 + 1$ puisque $s_0$ code une courbe r\'eelle ayant $\frac{1}{2} \mu_{s_0} + 1$ degr\'es de libert\'e. En outre, chaque ar\^ete de l'arbre $A_C$ code une paire de pointes des courbes cod\'ees par les sommets adjacents et qui ont une m\^eme orbite de Reeb pour asymptote. 
Cette contrainte co\^ute deux degr\'es de libert\'e suppl\'ementaires par ar\^ete, soit au total $2 \# S_1 + 2 \# S_2 - 2$ degr\'es, puisque le nombre d'ar\^etes d'un arbre diff\`ere du nombre de sommets
par un. Ceci nous permet finalement de minorer l'indice de Maslov de notre courbe \`a deux \'etages cod\'ee par $A_C$ par $r + 2r_X - 1 = c_d-1$. Cet indice est par ailleurs major\'e par le degr\'e $c_d-1$
du fibr\'e normal de $C$, de sorte que toutes nos minorations sont des \'egalit\'es. En particulier, tous les sommets \`a distance paire de $s_0$ autre que $s_0$ lui-m\^eme codent soit des cylindres 
convergeant vers des orbites de Reeb simplement rev\^etues, soit des plans convergeant vers des orbites doublement rev\^etues puisque ce sont d'apr\`es la Proposition \ref{propcotangent} les seules 
courbes de $T^* \R P^2$ rigides une fois leurs orbites de Reeb limites prescrites. 
Reprenons \`a ce stade la d\'emarche suivie dans le troisi\`eme paragraphe de la d\'emonstration 
du Th\'eor\`eme \ref{theocong1}. On note $r_1 , \dots , r_j$ le nombre
de paires de points complexes conjugu\'es de $\underline{x} $ associ\'ees \`a $s_1 , \dots , s_j$ respectivement, de sorte que leur somme vaille $r_X$. Il y a $2^{r_i - 1}$ partitions d'un ensemble
de $r_i$ points complexes conjugu\'es en deux ensembles complexes conjugu\'es, soit $2^{r_X - j}$ partitions au total. Une fois attribu\'es \`a chaque courbe $C_s$ l'ensemble de points qu'elle
doit interpoler, certaines de ces courbes sont rigides et d'autres non. Notons $j_1$ le nombre de telles courbes rigides non adjacentes \`a $C_{s_0}$ et $j^-$ le nombre de telles courbes rigides
adjacentes \`a $C_{s_0}$. Les $j^+$ courbes restantes sont adjacentes \`a $C_{s_0}$ et gardent deux degr\'es de libert\'e une fois interpol\'es les $r_k$ points qu'elles doivent interpoler ; c'est la condition
d'adjacence \`a $C_{s_0}$ qui les rigidifie. D'apr\`es ce qui pr\'ec\`ede, la courbe r\'eelle cod\'ee par $s_0$,
avec ses $j^- + c_{-2}$ paires d'asymptotes prescrites et ses $r$ points r\'eels \`a interpoler, est rigide, o\`u $c_{-2}$ d\'esigne \`a nouveau le nombre de courbes rigides adjacentes \`a $C_{s_0}$ et autres
que les $j^-$ courbes soumises \`a des conditions d'incidence. La dimension $k_{s_0} + 2v_{s_0} -1$ donn\'ee par la Proposition \ref{propcotangent}
vaut donc en particulier $r + 2j^- + 2c_{-2}$. Par cons\'equent, les $j^+$ paires de courbes non-rigides pr\'ec\'edentes h\'eritent d'une contrainte suppl\'ementaire, elles ont une paire
d'asymptotes prescrites correspondant \`a une paire d'orbites de Reeb limites rest\'ees libres de la courbe $C_{s_0}$. Il y a deux bijections possibles entre une telle paire d'orbites de Reeb et
une telle paire de courbes non-rigides, soit $2^{j^+}$ bijections au total. De m\^eme, les $j_1$ courbes rigides non adjacentes \`a $C_{s_0}$ sont cod\'ees par des sommets adjacents \`a au moins 
un sommet bivalent de l'arbre $A_C$ puisque ce dernier est connexe et que les autres sommets adjacents sont des feuilles. D'apr\`es ce qui pr\'ec\`ede, ce sommet bivalent code une paire de cylindres
complexes conjugu\'es de $T^* L$ reliant deux paires complexe d'orbites limites de deux paires complexes conjugu\'ees de courbes rigides de $X \setminus L$. Il y a deux fa\c{c}ons d'apparier ces
orbites, soit $2^{j_1}$ bijections au total. Ainsi, le nombre de courbes \`a deux \'etages ayant une combinatoire donn\'ee par $A_C$ est divisible par $2^{r_X - j^-}$.
Or d'apr\`es le th\'eor\`eme de recollement en th\'eorie symplectique des champs \cite{Bour} et le Th\'eor\`eme \ref{theoinv}, la contribution \`a l'invariant $\chi^d_{r} (L)$ d'une courbe \`a deux 
\'etages ne d\'epend que de sa combinatoire, de sorte que $2^{r_X - j^-}$ divise $\chi^d_{r} (L)$. Nous disposons cette fois-ci de la relation 
$r + 2j^- + 2c_{-2} = k_{s_0} + 2v_{s_0} -1$ qui impose l'in\'egalit\'e $r+1 \geq k_{s_0}$.
On en d\'eduit donc $j^- \leq v_{s_0} \leq k_{s_0} \leq r+1$ et le premier r\'esultat \'enonc\'e dans le Th\'eor\`eme \ref{theocong3}. 

Tous les arbres $A_C$ pour lesquels l'une des in\'egalit\'es $j^- \leq v_{s_0} \leq k_{s_0} \leq r+1$ est stricte satisfont $j^- \leq r$. Pour montrer le second r\'esultat \'enonc\'e dans le 
Th\'eor\`eme \ref{theocong3}, on peut donc se restreindre aux arbres satisfaisant $j^- = v_{s_0} = k_{s_0} = r+1$. En particulier, la courbe cod\'ee par $s_0$ n'a que des orbites simples
pour limites. Si un tel arbre poss\`ede une feuille \`a distance paire de $s_0$, on a vu qu'elle doit coder un plan asymptote \`a une orbite double. Le th\'eor\`eme de recollement en th\'eorie
symplectique des champs garantit alors qu'il y a deux fa\c{c}ons de recoller ce plan au restant de la courbe. Ce param\`etre de recollement permet donc \`a nouveau dans ce cas l\`a d'am\'eliorer
d'une puissance de deux le premier \'enonc\'e du Th\'eor\`eme \ref{theocong3}. On peut donc supposer que toutes les ar\^etes des arbres $A_C$ sont de multiplicit\'e un. 
Notons $C_1$ la r\'eunion des courbes cod\'ees par les sommets \`a distances impaires de $s_0$  et $\overline{C}_1$ sa compactifi\'ee dans $\Sigma_4$. La classe d'homologie de $\overline{C}_1$
s'\'ecrit $vf + g$, o\`u $f$ est la classe d'une fibre de $\Sigma_4$,  $v$ est le nombre d'ar\^etes de $A_C$ et $g \in H_2 (X \setminus L ; \Z)$. La classe d'homologie $d$ s'\'ecrit alors 
\begin{equation}
\label{relv}
v(f + (c_X)_* f ) + g +  (c_X)_* g = vh + g +  (c_X)_* g,
\end{equation}
 d'o\`u $\langle d , h \rangle = v \mod (4)$.
On d\'eduit donc des hypoth\`eses faites que $v = r-1 \mod (4)$, puis que le nombre de sommets \`a distance paire de $s_0$, $s_0$ exclu, est impair puisque $v_{s_0} = r+1$. 
Ceci force l'existence d'un sommet $s$
de valence paire parmi les sommets \`a distance impaire de $s_0$. Or, le nombre de courbes \`a deux \'etages cod\'ees par un tel arbre est pair. En effet, si $s$ est adjacent \`a $s_0$, il y a
parmi les pointes de $C_s$ un nombre pair de choix de celle reli\'ee \`a $s_0$. Si $s$ n'est pas adjacent \`a $s_0$, il est reli\'e \`a un nombre pair de sommets bivalents de $S_2$, eux-m\^emes reli\'es \`a un
nombre pair de sommets de $S_1$ de sorte que ces sommets bivalents ne font que connecter bijectivement ces derniers aux pointes de $C_s$. Le nombre de telles bijections \'etant pair,
nous pouvons \`a nouveau dans ce dernier cas d'am\'eliorer d'une puissance de deux le r\'esultat pr\'ec\'edent, ce qui d\'emontre le second \'enonc\'e du Th\'eor\`eme \ref{theocong3}. 

Enfin, dans le cas du plan projectif complexe, lorsque $r+1 = d-2$, on d\'eduit aussi de la relation (\ref{relv}) l'in\'egalit\'e $v \leq d$, in\'egalit\'e stricte d\`es que $\langle g , h \rangle \neq 0$. 
Or l'annulation $\langle g , h \rangle$ force les sommets adjacents \`a $s_0$ \`a \^etre des feuilles, ce qui est exclus par la majoration $k_{s_0} \leq r+1 < d$. Par suite,  
$v \leq d - 4$, $k_{s_0} \leq d-4 < r+1$ et l'\'equation $r + 2j^- + 2c_{-2} = k_{s_0} + 2v_{s_0} -1$ impose $j^- <  v_{s_0}$. Dans ce cas, on obtient donc $r_X - j^- =
d+1 - j^- \geq d + 2 - k_{s_0} \geq 6$, d'o\`u le r\'esultat. $\square$

\section{Calculs}
\label{sectcalculs}

\subsection{Invariants \'enum\'eratifs r\'eels de fibr\'es cotangents}

\subsubsection{Construction des espaces de modules}
\label{subsubmodules}

Soit $L$ une sph\`ere, un tore ou un espace projectif r\'eel de dimension $n=2$ ou $3$. Le fibr\'e cotangent
de $L$ est \'equip\'e de sa forme de Liouville $\lambda$ et de l'involution $c_L$ d\'efinie par $(q,p) \in T^* L \mapsto (q,-p) \in T^* L$. 
Cette derni\`ere satisfait $c_L^* \lambda = - \lambda$ de sorte que $(T^* L , d\lambda , c_L)$ est une vari\'et\'e symplectique r\'eelle.
Soit $g$ une m\'etrique \`a courbure constante sur $L$, $U^* L$  l'ensemble des couples $(q,p) \in T^* L$ tels que $g(p,p) \leq 1$ et $S^* L$ 
le bord de $U^* L$. La restriction de $\lambda$ \`a $S^* L$ est une forme de contact et l'on note $R_\lambda$
le champ de Reeb associ\'e. Le flot engendr\'e par $R_\lambda$ n'est autre que le flot g\'eod\'esique.
Notons ${\cal J}_\lambda$ l'espace des structures presque-complexes positives pour $d \lambda$ et asymptotiquement cylindriques
sur une structure $CR$ de $S^* L$.   Plus pr\'ecis\'ement, le champ radial de $ T^* L$ identifie le compl\'ementaire de la section nulle
avec la symplectisation $(\mathbb{R} \times S^* L , d(e^\rho \lambda))$ de $(S^* L , \lambda)$. On note ${\cal J}_\lambda$ 
l'espace des structures presque-complexes $J$ positives pour $d \lambda$, de classe $C^l$, $l \gg 1$,  qui satisfont
$J(\frac{\partial}{\partial \rho}) = R_\lambda$ et pr\'eservent le noyau de $\lambda$ pour  $\rho \gg 1$ et qui enfin sont invariantes par translation par $\rho$ au-del\`a
d'un certain rang $\rho_0$. Nous notons alors $\mathbb{R} {\cal J}_\lambda \subset {\cal J}_\lambda$ le sous-espace des structures
presque-complexes pour lesquelles $c_L$ est $J$-antiholomorphe. Ces espaces ${\cal J}_\lambda$ et $\Bbb{R} {\cal J}_\lambda$ 
sont tous deux des vari\'et\'es de Banach s\'eparable non-vides et contractiles.

Soit $S$ une sph\`ere de dimension deux orient\'ee et ${\cal J_S} $ l'espace des structures presque-complexes de classe $C^l$ sur $S$ qui sont 
compatibles avec son orientation. Soient $v_\C \in \N^*$ et $y_1, \dots , y_{v_\C}$ une collection de $v_\C$ points distincts sur $S$. Il suit de
\cite{HWZ1} et du Corollaire $5.1$ de \cite{Bour} qu'il existe $0 < d << 1$ tel que pour tous $J_S \in {\cal J_S} $, $J \in {\cal J_\lambda} $ et toute application $J$-holomorphe propre 
$u : S \setminus \{y_1, \dots , y_{v_\C} \} \to T^* L$ d'\'energie de Hofer finie, l'application $u$ a le comportement suivant au voisinage de chaque point $y_i$. Fixons un param\'etrage
local $J_S$-holomorphe de $S$ au voisinage de chaque point $y_i$, $1 \leq i \leq v_\C$ par l'anneau $\C_1 = \{ z \in \C \; \vert \; \vert z \vert \geq 1\}$, puis des coordonn\'ees cylindriques
$(s,t) \in \R^* \times [0 , 1] \mapsto e^{s + 2\pi it} \in \C_1$ de cet anneau. On en d\'eduit un param\'etrage de $S$ au voisinage de chaque point $y_i$ de la forme
$\phi_i : (s,t) \in \R^* \times [0 , 1]  \to S \setminus \{y_1, \dots , y_{v_\C} \}$, $1 \leq i \leq v_\C$. Pour $1 \leq i \leq v_\C$, notons $u \circ \phi_i = (\rho_i , \tilde{u}_i)$, 
o\`u $\rho_i : (s,t)  \in \R^+ \times [0 , 1]  \to \R$ et $\tilde{u}_i : (s,t)  \in \R^+ \times [0 , 1]  \to S^* L$. Alors, pour $1 \leq i \leq v_\C$, il existe $s_i \in \R$, $k_i \in \N^*$ et
des orbites $\gamma_i$ du flot de Reeb tels que les fonctions distances $\vert \rho_i (s,t) -(k_i A)s - s_i \vert$ et $d \big( \tilde{u} (s,t) - \gamma_i ((k_i A)t) \big)$
appartiennent \`a l'espace fonctionnel $L^{k,p}_d = \{ f : \R^+ \times [0 , 1] \to \R^+ \; \vert \; f(s,t) e^{ds} \in L^{k,p} (\R^+ \times [0 , 1] , \R) \}$ et ceci quel que soient $1 << k << l$ et $2 < p < + \infty$, o\`u
$A$ d\'esigne l'int\'egrale de $\lambda$ sur l'orbite de Reeb simple sous-jacente et $L^{k,p} (\R^+ \times [0 , 1] , \R)$ d\'esigne l'espace des fonctions ayant $k$ d\'eriv\'ees dans $L^p$. On note
$L^{k,p}_d (S \setminus \{y_1, \dots , y_{v_\C} \}  , T^*L)$ l'espace des fonctions propres $u : S \setminus \{y_1, \dots , y_{v_\C} \} \to T^* L$ ayant cette propri\'et\'e.
C'est une vari\'et\'e de Banach s\'eparable. Soient \`a pr\'esent $v_\C^- \in \N^*$, $k_1, \dots , k_{v_\C} \in \N^*$ et $\gamma_1, \dots , \gamma_{v_\C^-}$une collection d'orbites disjointes
du flot de Reeb, de sorte que lorsque $v_\C^-$ s'annule, $\Gamma = \{ \gamma_1, \dots , \gamma_{v^-_\C}\}$ soit vide. Soient $r_\C \in \N$, $z_1 , \dots , z_{r_\C} \in
S \setminus \{y_1, \dots , y_{v_\C} \} $ et $x_1 , \dots , x_{r_\C} \in T^* L$ une collection de points distincts, de sorte qu'\`a nouveau, lorsque $r_\C$ s'annule, $ \{ x_1, \dots , x_{r_\C}\}$ soit vide. 
Soit alors
$${\cal P} (\Gamma) = \{ (u, J_S , J) \in L^{k,p}_d (S \setminus \{y_1, \dots , y_{v_\C} \}  , T^*L) \times {\cal J_S} \times {\cal J_\lambda} \; \vert \; du + J \circ du \circ J_S = 0 , \; u (z_i) = x_i, $$
$$\lim_{z \to y_i^-} u(z) = k_i \gamma_i \text{ et } \lim_{s \to +\infty} \int (u \circ \phi_i)^* \lambda = k_i A \},$$
et ${\cal P}^* (\Gamma) \subset {\cal P} (\Gamma)$ le sous-espace des applications pseudo-holomorphes non-multiples.
Soit ${\cal D}iff^+ (S , \underline{z} , \underline{y})$ le groupe des diff\'eomorphismes de classe $C^{l+1}$ de $S$ qui pr\'eservent l'orientation et fixent $ \underline{z} , \underline{y}$.
Ce groupe agit sur ${\cal P}^* (\Gamma)$ par $(\phi , (u , J_S , J)) \in {\cal D}iff^+ (S , \underline{z} , \underline{y}) \times {\cal P}^* (\Gamma) \mapsto (u \circ \phi^{-1} ,
\phi^* J_S , J) \in {\cal P}^* (\Gamma) $,  o\`u $\phi^* J_S  = d \phi \circ J_S \circ d\phi^{-1}$. On note ${\cal M}_{r_\C}^{v_\C} (\Gamma , \underline{x})$ le quotient 
${\cal P}^* (\Gamma) / {\cal D}iff^+ (S , \underline{z} , \underline{y})$ et $\pi : {\cal M}_{r_\C}^{v_\C} (\Gamma , \underline{x}) \to {\cal J}_\lambda $ la projection induite par
$(u , J_S , J) \in {\cal P}^* (\Gamma) \to J \in {\cal J}_\lambda $

Fixons une m\'etrique $g_L$ sur $T^* L$ pr\'eserv\'ee par $c_L$ et invariante par translation par $\rho$ pour $\rho$ assez grand. Elle induit une connexion $\nabla$
sur $TT^* L$ et tous les fibr\'es associ\'es. Si $(u , J_S , J) \in {\cal P}^* (\Gamma)$, on note $D$ l'op\'erateur de Gromov $v \in L^{k,p}_d (S \setminus \{y_1, \dots , y_{v_\C} \}  , T^*L)
\mapsto \nabla v + J \circ \nabla v \circ J_S + \nabla_v J \circ du \circ J_S \in L^{k-1,p}_d (S \setminus \{y_1, \dots , y_{v_\C} \}  , \Lambda^{0,1} S \otimes u^*TT^*L)$.
Ce dernier induit un op\'erateur $\overline{D} : L^{k,p}_d (S \setminus \{y_1, \dots , y_{v_\C} \}  , u^*TT^*L) / du(L^{k,p}_d (S \setminus \{y_1, \dots , y_{v_\C} \}  , TS)
\to L^{k-1,p}_d (S \setminus \{y_1, \dots , y_{v_\C} \}  , \Lambda^{0,1} S \otimes u^*TT^*L) / du(L^{k-1,p}_d (S \setminus \{y_1, \dots , y_{v_\C} \}  , \Lambda^{0,1} S \otimes TS))$
(voir la formule $1.5.1$ de \cite{IShev} et le paragraphe $1.4$ de \cite{Wels1}). On note $H^0_{\overline{D}} ( S ; {\cal N}_{u , -\underline{z}})$ (resp. 
$H^1_{\overline{D}} ( S ; {\cal N}_{u , -\underline{z}})$) le noyau (resp. conoyau) de $\overline{D}$.

\begin{prop}
\label{propmod}
L'espace ${\cal M}_{r_\C}^{v_\C} (\Gamma , \underline{x})$ est une vari\'et\'e de Banach s\'eparable de classe $C^{l-k}$ La projection $\pi$ est Fredholm, de noyau (resp. conoyau) 
isomorphe \`a $H^0_{\overline{D}} ( S ; {\cal N}_{u , -\underline{z}})$ (resp. 
$H^1_{\overline{D}} ( S ; {\cal N}_{u , -\underline{z}})$). En particulier, l'indice de $\pi$ vaut $2(n-1)\sum_{i=1}^{v_\C} k_i + 2v_\C - 2 - 2(n-1)r_\C - 2(n-1)v_\C^-$ si $L$ est une sph\`ere,
$(n-1)\sum_{i=1}^{v_\C} k_i + 2v_\C - 2 - 2(n-1)r_\C - 2(n-1)v_\C^-$ si $L$ est un plan projectif r\'eel et $2v_\C + 2n - 6 - 2(n-1)r_\C - (n-1)v_\C^-$ si $L$ est un tore.
\end{prop}

{\bf D\'emonstration :}

La premi\`ere partie de la Proposition \ref{propmod} est classique. L'identification des noyau et conoyau de $\pi$ avec ceux de $ \overline{D}$ se d\'emontre comme le Th\'eor\`eme $2$ de \cite{IShev}.
Le caract\`ere Fredhom de $D$ d\'ecoule de la Proposition $5.2$ de \cite{Bour} (voir aussi le sixi\`eme paragraphe de \cite{HWZ}). Le calcul des indices d\'ecoule des Propositions
\ref{propmaslov} et \ref{propcotangent}. $\square$\\

Supposons \`a pr\'esent que $v_\C = 2 v_\R$, $v_\C^- = 2 v_\R^-$ et $r_\C = r + 2r_L$. Pour $1 \leq i \leq v_\R$ (resp. $1 \leq j \leq v_\R^-$ ), on note $y_{v_\R + i } = \overline{y}_i$
(resp. $\gamma_{v_\R^- + j } = \overline{\gamma}_j$).  De m\^eme, pour $1 \leq i \leq r_L$, on note $z_{r + r_L + i} = \overline{z}_{r+i}$ et $x_{r + r_L + i} = \overline{x}_{r+i}$.
On suppose cette fois-ci que $x_1 , \dots , x_r \in L \subset T^*L$ et que pour $1 \leq i \leq r_L$ (resp. $1 \leq j \leq v_\R^-$ ), $\overline{x}_{r+i} = c_L (x_{r+i})$ (resp. 
$ \overline{\gamma}_j = -c_L (\gamma_j )$). On note ${\cal D}iff (S , \underline{z} , \underline{y})$ le groupe des diff\'eomorphismes de classe $C^{l+1}$ de $S$ qui fixent
$\underline{z}$ et  $\underline{y}$ s'ils pr\'eservent l'orientation ou bien fixent $z_1 , \dots z_r$ et \'echangent $y_i$, $\overline{y}_i$ et $z_{r+j}$, $\overline{z}_{r+j}$ s'ils renversent
l'orientation, $1 \leq i \leq v_\R$, $1 \leq j \leq r_L$. Sous les hypoth\`eses que l'on vient de faire, cette extension d'indice deux de ${\cal D}iff^+ (S , \underline{z} , \underline{y})$ 
agit \'egalement sur ${\cal P}^* (\Gamma)$ par $(\phi , (u , J_S , J)) \in {\cal D}iff (S , \underline{z} , \underline{y}) \times {\cal P}^* (\Gamma) \mapsto (c_L \circ u \circ \phi^{-1} ,
\phi^* J_S , \overline{c}_L^* J) \in {\cal P}^* (\Gamma) $ lorsque $\phi \notin  {\cal D}iff^+ (S , \underline{z} , \underline{y})$, o\`u $\overline{c}_L^* J = -dc_L \circ J \circ dc_L$.
Par suite, le quotient ${\cal M}_{r_\C}^{v_\C} (\Gamma , \underline{x}) = {\cal P}^* (\Gamma) / {\cal D}iff^+ (S , \underline{z} , \underline{y})$ se trouve
\`a pr\'esent \'equip\'e d'une action de $\Z / 2\Z = {\cal D}iff (S , \underline{z} , \underline{y}) / {\cal D}iff^+ (S , \underline{z} , \underline{y})$, not\'ee $c_{\cal M}$. La projection
$\pi : ({\cal M}_{r_\C}^{v_\C} (\Gamma , \underline{x}) , c_{\cal M}) \to ({\cal J}_\lambda , \overline{c}_L^*) $ est alors  $\Z / 2\Z $-\'equivariante. On note 
$\R {\cal M}_{(r , r_L)}^{v_\R} (\Gamma , \underline{x})$
le lieu fixe de $c_{\cal M}$ et $\pi_\R : \R {\cal M}_{(r , r_L)}^{v_\R} (\Gamma , \underline{x}) \to \R {\cal J}_\lambda$ la projection induite par $\pi$. De la m\^eme mani\`ere que dans
\cite{Wels1}, les seuls \'el\'ements de ${\cal D}iff (S , \underline{z} , \underline{y}) $ qui peuvent avoir des points fixes dans ${\cal P}^* (\Gamma) $ sont d'ordre deux et renversent l'orientation de $S$
(voir le Lemme $1.3$ de \cite{Wels1}) et l'op\'erateur $D$ est ${\cal D}iff (S , \underline{z} , \underline{y}) $-\'equivariant (Lemme $1.5$ de \cite{Wels1}). Si 
$c_S \in {\cal D}iff (S , \underline{z} , \underline{y}) $ est un tel \'el\'ement d'ordre deux et $ (u , J_S , J) \in  {\cal P}^* (\Gamma) $ un point fixe de $c_S$, on note 
$L^{k,p}_d (S \setminus \{y_1, \dots , y_{v_\C} \}  , u^*TT^*L)_{+1}$, $L^{k,p}_d (S \setminus \{y_1, \dots , y_{v_\C} \}  , TS)_{+1}$ et
$ L^{k-1,p}_d (S \setminus \{y_1, \dots , y_{v_\C} \}  , \Lambda^{0,1} S \otimes u^*TT^*L)_{+1}$, $ L^{k-1,p}_d (S \setminus \{y_1, \dots , y_{v_\C} \}  , \Lambda^{0,1} S \otimes TS)_{+1}$
les espaces propres associ\'es aux valeurs propres $+1$ de l'action de $c_S$ sur
$L^{k,p}_d (S \setminus \{y_1, \dots , y_{v_\C} \}  , u^*TT^*L)$, $L^{k,p}_d (S \setminus \{y_1, \dots , y_{v_\C} \}  , TS)$ et
$ L^{k-1,p}_d (S \setminus \{y_1, \dots , y_{v_\C} \}  , \Lambda^{0,1} S \otimes u^*TT^*L)$, $ L^{k-1,p}_d (S \setminus \{y_1, \dots , y_{v_\C} \}  , \Lambda^{0,1} S \otimes TS)$
respectivement.  On note alors $\overline{D}_\R $ l'op\'erateur induit $L^{k,p}_d (S \setminus \{y_1, \dots , y_{v_\C} \}  , u^*TT^*L)_{+1} / du(L^{k,p}_d (S \setminus \{y_1, \dots , y_{v_\C} \}  , TS)_{+1})
\to L^{k-1,p}_d (S \setminus \{y_1, \dots , y_{v_\C} \}  , \Lambda^{0,1} S \otimes u^*TT^*L)_{+1} / du( L^{k-1,p}_d (S \setminus \{y_1, \dots , y_{v_\C} \}  , $ $\Lambda^{0,1} S \otimes TS)_{+1})$
et $H^0_{\overline{D}} ( S ; {\cal N}_{u , -\underline{z}})_{+1}$, $H^1_{\overline{D}} ( S ; {\cal N}_{u , -\underline{z}})_{+1}$ ses  noyau et conoyau.

\begin{prop}
\label{propmodreel}
L'espace $\R {\cal M}_{(r , r_L)}^{v_\R} (\Gamma , \underline{x})$ est une vari\'et\'e de Banach s\'eparable de classe $C^{l-k}$. La projection $\pi_\R$ est Fredholm, de noyau (resp. conoyau) 
isomorphe \`a $H^0_{\overline{D}} ( S ; {\cal N}_{u , -\underline{z}})_{+1}$ (resp. 
$H^1_{\overline{D}} ( S ; {\cal N}_{u , -\underline{z}})_{+1}$). En particulier, l'indice de $\pi_\R$ vaut $2(n-1)\sum_{i=1}^{v_\R} k_i + 2v_\R - 1 - (n-1)r - 2(n-1)r_L - 2(n-1)v_\R^-$ si $L$ est une sph\`ere,
$(n-1)\sum_{i=1}^{v_\R} k_i + 2v_\R - 1 - (n-1)r - 2(n-1)r_L - 2(n-1)v_\R^-$ si $L$ est un plan projectif r\'eel et $2v_\R + n - 3 - (n-1)r - 2(n-1)r_L - (n-1)v_\R^-$ si $L$ est un tore. $\square$
\end{prop}
La d\'emonstration de cette proposition est strictement analogue \`a celle de la Proposition $1.9$ de \cite{Wels1} et n'est pas reproduite ici.

\subsubsection{D\'efinition des invariants}
\label{subsubinvariants}

Nous allons compter les courbes $J$-holomorphes rationnelles r\'eelles point\'ees d'\'energie de Hofer finie proprement immerg\'ees dans $T^* L$ en fonction d'un signe $\pm 1$
de fa\c{c}on \`a obtenir un invariant associ\'e \`a $T^* L$. Rappelons que d'apr\`es le Th\'eor\`eme $1.2$ de \cite{HWZ1} et d'apr\`es \cite{Bour}, ces courbes rationnelles point\'ees
convergent en leurs pointes vers des orbites de Reeb parcourues un nombre entier de fois, que l'on appelle multiplicit\'e.
La dimension de l'espace des modules de telles courbes a \'et\'e calcul\'ee dans la Proposition \ref{propcotangent} et d\'epend du nombre de pointes
et des multiplicit\'es associ\'ees. Afin d'obtenir un nombre fini de courbes, nous allons soumettre ces courbes \`a quelques contraintes, soit en les
for\c{c}ant \`a converger vers des orbites de Reeb prescrites, soit en les for\c{c}ant \`a passer par des points de $L$ ou des paires de points
complexes conjugu\'ees de $T^* L \setminus L$.
Soit $e_i$, $i \geq 1$, la suite d'entiers partout nulle sauf au $i$-\`eme rang o\`u elle vaut un. Soient $\alpha = \sum_{i \in \Bbb{N}^*} \alpha_i e_i$ 
et $\beta = \sum_{i \in \Bbb{N}^*} \beta_i e_i$ deux suites d'entiers positifs qui s'annulent \`a partir d'un certain rang. Ces deux suites
codent respectivement le nombre de paires d'orbites de Reeb complexes conjugu\'ees limites prescrites et non prescrites de nos courbes, avec leur multiplicit\'es $i \in \Bbb{N}^*$.
Le nombre de pointes de nos courbes vaut donc $2v = 2 \sum_{i \in \Bbb{N}^*} (\alpha_i + \beta_i )$ et nous choisissons un ensemble 
$\Gamma$ de $\sum_{i \in \Bbb{N}^*} \alpha_i $ g\'eod\'esiques ferm\'ees disjointes de $L$ pour prescrire nos paires d'orbites de Reeb limites.
\`A pr\'esent, afin de fixer nos contraintes ponctuelles, soient  $r \in \Bbb{N}$ et $x_1 , \dots , x_r$  des points distincts de $L$. 
De m\^eme, soient $r_L \in \Bbb{N}$ et $\xi_1, \overline{\xi}_1 , \dots ,
\xi_{r_L}, \overline{\xi}_{r_L}$ des paires distinctes de points complexes conjugu\'es de $T^* L \setminus L$, c'est-\`a-dire satisfaisant $c_L (\xi_i) = \overline{\xi}_i$.
Nous supposons que
\begin{equation}
\label{dimsphere}
(n-1)r + 2(n-1)r_L + 2(n-1) \# \Gamma = 2 v + \epsilon (n-1) \sum_{i \in \Bbb{N}^*} i (\alpha_i + \beta_i ) + n - 3,
\end{equation}
o\`u $\epsilon = 2$ si $L$ est hom\'eomorphe \`a une sph\`ere et  $\epsilon = 1$ si $L$ est hom\'eomorphe \`a un espace projectif r\'eel, tandis que nous supposons
\begin{equation}
\label{dimtore}
(n-1)r + 2(n-1)r_L  = 2 v + n - 3 \text{ et } \alpha = 0
\end{equation}
si $L$ est hom\'eomorphe \`a un tore.

Alors, lorsque la structure presque-complexe $J \in \Bbb{R} {\cal J}_\lambda$ est g\'en\'erique, il n'y a qu'un nombre fini
de courbes $J$-holomorphes rationnelles r\'eelles d'\'energie de Hofer finie, proprement immerg\'ees dans $T^* L$
et ayant $2v$ pointes qui passent par $\underline{x}$, par chaque paire $\{ \xi_i ,
\overline{\xi}_i \}$ et qui convergent vers les orbites de Reeb relevant les \'el\'ements de $\Gamma$  ainsi que vers
$\beta_j$ autres paires d'orbites, $j \in \Bbb{N}^*$, chacune avec multiplicit\'e $j$ ou de classe d'homologie donn\'ee si $L$ est un tore.
En effet, si $L$ est un tore, il y a une infinit\'e de g\'eod\'esiques ferm\'ees primitives non homologues et la dimension (\ref{dimtore}) ne d\'epend pas du choix
des classes d'homologies de sorte qu'il y a une infinit\'e d'espaces de modules ayant la m\^eme dimension. Pour garantir la finitude, nous imposons les classes
 d'homologies des orbites de Reeb limites. 
Notons ${\cal R} (\alpha , \beta , \Gamma , \underline{x} , \underline{\xi} , J)$ cet ensemble fini de courbes,
la g\'en\'ericit\'e de $J$ garantit qu'elles sont toutes immerg\'ees. 
Si $L$ est de dimension deux, on d\'efinit alors comme dans \cite{Wels1} la masse $m(C)$ d'une telle courbe $C$ comme
le nombre fini de ses points doubles r\'eels isol\'es, c'est-\`a-dire de ses points doubles situ\'es sur $L$ et qui sont l'intersection
transverses de deux branches complexes conjugu\'ees. On pose 
$$F_{(r, r_L)} (\alpha , \beta , \Gamma , \underline{x} , \underline{\xi} , J) = \sum_{C \in {\cal R} (\alpha , \beta , \Gamma , \underline{x} , \underline{\xi} , J)} (-1)^{m(C)} \in \Bbb{Z}.$$
Si $L$ est de dimension trois, on l'\'equipe d'une structure spin. Ceci permet d'associer un \'etat spinoriel $\text{sp} (C)$ \`a chaque courbe 
$C \in {\cal R} (\alpha , \beta , \Gamma , \underline{x} , \underline{\xi} , J)$ comme expliqu\'e au paragraphe $5.2$ de \cite{Wels3}.
Cet \'etat  spinoriel est d\'efini comme suit. La lin\'earisation de l'\'equation de Cauchy-Riemann en $C$ fournit un op\'erateur surjectif de Cauchy-Riemann g\'en\'eralis\'e
d\'efini sur un espace de Banach de sections du fibr\'e normal de $C$ \`a valeurs dans un espace de Banach de formes de type $(0,1)$ \`a valeurs dans le
fibr\'e normal de $C$.  Si cet op\'erateur est $\C$-lin\'eaire, il induit une structure de fibr\'e vectoriel holomorphe sur le fibr\'e normal de $C$ qui se d\'ecompose comme
somme \'equilibr\'ee de fibr\'es en droites complexes. Cette d\'ecomposition fournit un rep\`ere mobile le long de la partie r\'eelle de $C$ qui permet de d\'efinir l'\'etat spinoriel de $C$
comme l'obstruction \`a relever ce rep\`ere \`a un rep\`ere du fibr\'e des spineurs, voir \cite{Wels2}, \cite{Wels3}. Si cet op\'erateur surjectif de Cauchy-Riemann g\'en\'eralis\'e n'est que
$\R$-lin\'eaire, il peut \^etre reli\'e \`a un op\'erateur $\C$-lin\'eaire par un chemin transverse \`a l'espace des op\'erateurs non-surjectifs. L'\'etat spinoriel de $C$ est
alors d\'efini comme \'etat spinoriel de l'op\'erateur $\C$-lin\'eaire corrig\'{e} par la parit\'e du nombre d'intersection du chemin choisi avec le mur des op\'erateurs non-surjectifs,
voir le \S $5.2$ de \cite{Wels3}. On pose alors
$$F_{(r, r_L)} (\alpha , \beta , \Gamma , \underline{x} , \underline{\xi} , J) = \sum_{C \in {\cal R} (\alpha , \beta , \Gamma , \underline{x} , \underline{\xi} , J)} \text{sp} (C) \in \Bbb{Z}.$$

\begin{theo}
\label{theoinv}
Soit $L$ une sph\`ere, un tore ou un espace projectif r\'eel de dimension $n=2$ ou $3$ muni d'une m\'etrique \`a courbure constante. Soient
$\alpha , \beta$ deux suites d'entiers positifs qui s'annulent \`a partir d'un certain rang. On choisit comme ci-dessus
un ensemble $\Gamma$ de g\'eod\'esiques ferm\'ees et des ensembles $\underline{x}$, $\underline{\xi}$
de $r$ et $r_L$ points dans $L$ et $T^* L \setminus L$ respectivement de sorte que ces nombres satisfassent (\ref{dimtore}) dans le cas du tore et (\ref{dimsphere}) sinon. 
Lorsque $n=3$, on suppose $r \neq 0$ et lorsque de plus $L \in \{ S^3 , \R P^3 \}$, on suppose que $J$ est invariante par le flot de Reeb pour $\rho \gg 1$. Alors, l'entier $F_{(r, r_L)} (\alpha , \beta , \Gamma , \underline{x} , \underline{\xi} , J)$ d\'efini ci-dessus ne d\'epend ni du choix 
des contraintes $\Gamma , \underline{x} , \underline{\xi}$, ni du choix g\'en\'erique de la structure presque-complexe $J \in \Bbb{R} {\cal J}_\lambda$.
\end{theo}

Nous n'utiliserons ce Th\'eor\`eme \ref{theoinv} que dans le cas de la sph\`ere et du plan projectif r\'eel mais incluons toutefois le cas du tore ou de l'espace projectif de dimension trois
puisque la d\'emonstration est analogue. Remarquons que dans le cas de la sph\`ere ou de l'espace projectif r\'eel, $F_{(r, r_L)} (\alpha , \beta)$ n'est autre qu'un invariant r\'eel
analogue \`a celui d\'efini dans \cite{Wels1}, \cite{Wels2} relatif \`a la quadrique imaginaire pure si $L \in \{ \R P^2 , \R P^3 \}$ ou bien relatif \`a une section hyperplane r\'eelle de
la quadrique disjointe de l'ellipso\"{\i}de si  $L \in \{ S^2 , S^3 \}$. L'existence d'un tel invariant relatif a \'et\'e ind\'ependamment observ\'ee par Cheol-Hyun Cho dans \cite{Cho}.
Le lien entre la th\'eorie des invariants relatifs et le point de vue de la th\'eorie symplectique des champs est d\'evelopp\'e dans \cite{Katz}.\\

{\bf D\'emonstration :}

On consid\`ere l'espace des modules $\R {\cal M}^v_{r, r_L} (\Gamma ,  \underline{x} , \underline{\xi})$ des sph\`eres $J$-holomorphes r\'eelles d'\'energie de Hofer finie 
proprement immerg\'ees dans $T^* L$ ayant $2v$ pointes qui passent par $\underline{x}$, par chaque paire $\{ \xi_i ,\overline{\xi}_i \}$ et qui convergent vers les 
orbites de Reeb relevant les \'el\'ements de $\Gamma$  ainsi que vers $\beta_j$ autres paires d'orbites, $j \in \Bbb{N}^*$, chacune avec multiplicit\'e $j$, voir le \S \ref{subsubmodules}.
Soient $J_0$, $J_1 \in \R {\cal J}_\lambda$ deux structures presque complexes g\'en\'eriques de sorte que $F_{(r, r_L)} (\alpha , \beta , \Gamma , \underline{x} , \underline{\xi} , J_0)$ et 
$F_{(r, r_L)} (\alpha , \beta , \Gamma , \underline{x} , \underline{\xi} , J_1)$ soient bien d\'efinis. Soient $\gamma : [0,1] \to \Bbb{R} {\cal J}_\lambda$
une homotopie g\'en\'erique reliant $J_0$ \`a $J_1$, $\R {\cal M}_\gamma = \R {\cal M}^v_{r, r_L} (\Gamma ,  \underline{x} , \underline{\xi}) \times_\gamma [0,1] $ et
$\pi_\gamma : \R {\cal M}_\gamma  \to [0,1] $ la projection associ\'ee. 
Supposons pour commencer que $n=2$.
Les seuls points \`a \'etudier sont l'absence de compacit\'e de $\R {\cal M}_\gamma$ et les points critiques de
$\pi_\gamma$. En effet, en dehors de ce nombre fini de valeurs de $[0,1]$, les seules autres valeurs $t$ o\`u $F_{(r, r_L)} (\alpha , \beta , \Gamma , \underline{x} , \underline{\xi} , J_t)$
n'est pas d\'efini correspondent \`a des courbes ayant un point triple r\'eel ordinaire ou un point de tangence non-d\'eg\'en\'er\'e, et comme dans \cite{Wels1}, l'invariance de 
$F_{(r, r_L)} (\alpha , \beta , \Gamma , \underline{x} , \underline{\xi} , J_t)$ au passage de ces valeurs se v\'erifie facilement. Les points critiques de $\pi_\gamma $ correspondent aux courbes
ayant un unique  point de rebroussement de premi\`ere esp\`ece ordinaire. En effet, d'apr\`es la Proposition \ref{propmodreel}, en un tel point critique l'espace 
$H^0_{\overline{D}} ( S ; {\cal N}_{u , -\underline{z}})_{+1}$ ne s'annule pas. Or par d\'efinition, l'indice total des z\'eros d'une section du fibr\'e normal \`a une courbe immerg\'ee
vaut l'indice de Maslov de cette courbe. Une telle section doit par ailleurs s'annuler en $\underline{x} , \underline{\xi}$ et en les $v^-$ pointes prescrites, ce qui ne se peut pas. 
Par suite, le g\'en\'erateur de $H^0_{\overline{D}} ( S ; {\cal N}_{u , -\underline{z}})_{+1}$ est forc\'ement de torsion, de sorte que la courbe n'est pas immerg\'ee. La g\'en\'ericit\'e de $\gamma$
assure l'unicit\'e du point de rebroussement et son caract\`ere ordinaire. Ces courbes sont des points critiques non-d\'eg\'en\'er\'es de $\pi_\gamma$ qui correspondent \`a l'apparition ou la disparition de
deux courbes dont la masse diff\`ere de un. Nous ne reproduisons pas la d\'emonstration de ces deux faits ici puisqu'elle est strictement analogue \`a celle de \cite{Wels1}. L'invariance de
$F_{(r, r_L)} (\alpha , \beta , \Gamma , \underline{x} , \underline{\xi} , J)$  au passage de telles valeurs critiques $t$ de $\pi_\gamma$  en d\'ecoule.
D'apr\`es le th\'eor\`eme de compacit\'e en th\'eorie symplectique des champs \cite{BEHWZ}, l'absence de compacit\'e de $\R {\cal M}_\gamma $ peut provenir de trois ph\'enom\`enes, \`a savoir la d\'{e}g\'{e}n\'{e}rescence d'une suite d'\'el\'ements de $\R {\cal M}_\gamma $ vers une courbe multiple, r\'eductible ou \`a plusieurs \'etages. Supposons pour commencer qu'une telle suite d\'eg\'en\`ere 
vers un rev\^etement $l$-uple d'une courbe $C'$ de $T^*L$, $l \geq 2$. Alors, le nombre de pointes de $C'$ est inf\'erieur \`a celui de $C$ et la somme des multiplicit\'es associ\'ees est inf\'erieure au 
$l^{\text{\`eme}}$ de celle de $C$.
Par suite, dans le cas d'une sph\`ere, l'espace des modules contenant $C'$ vient avec une projection Fredholm sur $\Bbb{R} {\cal J}_\lambda $ dont l'indice est major\'e par celui de $C$ qui est nul moins
le double de la somme des multiplicit\'e des pointes. De telles courbes multiples ne peuvent appara\^{\i}tre en codimension un. Elles le peuvent dans le cas du plan projectif uniquement lorsque
$C$ est un rev\^etement double d'un cylindre $C'$ sur des orbites de Reeb simples, ramifi\'e en les pointes. Dans ce cas, $r$ vaut un ou trois et la courbe $C'$ a une partie r\'eelle connexe, sans point 
double et non triviale dans $H_1 (\R P^2 ; \Z/2\Z)$. Le compl\'ementaire de ces parties r\'eelles est donc toujours connexe par arc et par suite le compl\'ementaire dans $\Bbb{R} {\cal J}_\lambda $ 
des structures presque-complexes pour lesquelles une courbe multiple satisfait nos conditions d'incidence est lui aussi connexe par arc. Dans le cas du tore enfin, l'espace des modules
contenant $C'$ vient avec une projection Fredholm sur $\Bbb{R} {\cal J}_\lambda $ dont l'indice est major\'e par celui de $C$ qui est nul moins deux sauf si le nombre de pointes de  $C'$
est le m\^eme que celui de $C$. Notant $v'$ ce nombre de pointes, la formule de Riemann-Hurwitz impose que l'indice total des points de ramification situ\'es au dessus de ces pointes
vaille $(l-1)v'$. Cet indice de ramification \'etant major\'e par $2l-2$, cela force $C'$ et $C$ \`a \^etre des cylindres et  $F_{(r, r_L)} (\alpha , \beta , \Gamma , \underline{x} , \underline{\xi} , J) = 1$.
Dans tous les cas, l'\'eventuelle d\'eg\'en\'erescence vers des courbes multiples ne fait pas obstacle \`a l'invariance de $F_{(r, r_L)} (\alpha , \beta , \Gamma , \underline{x} , \underline{\xi} , J)$.
Supposons \`a pr\'esent qu'une telle suite d'\'el\'ements de $\R {\cal M}_\gamma $ converge vers une courbe $J_{t_0}$-holomorphe r\'eductible $C_{t_0}$. La g\'en\'ericit\'e de $\gamma$ impose alors 
que cette courbe r\'eductible
poss\`ede deux composantes irr\'eductibles, toutes deux r\'eelles et que ses  points singuliers soient des points doubles ordinaires. De plus, d'apr\`es ce que l'on vient de voir, ces composantes 
doivent \^etre toutes deux simples lorsque $L$ est une sph\`ere ou un plan projectif r\'eel, mais peuvent \^etre des rev\^etements de cylindres dans le cas du tore. Excluons ce dernier cas
pour commencer. Il se produit alors le m\^eme  ph\'enom\`ene que dans  les vari\'et\'es ferm\'ees, voir \cite{Wels1}, \`a savoir que pour toute structure presque-complexe $J_t$ proche de $J_{t_0}$,
\`a l'exclusion de $J_{t_0}$, et pour chaque point d'intersection r\'eel $p$ entre les deux composantes irr\'eductibles de $C_{t_0}$, il y a exactement une courbe $J_t$-holomorphe dans 
$T^* L$ satisfaisant nos conditions d'incidence. En effet, s'il y en avait deux, elles s'intersecteraient en chaque point de notre configuration $\underline{x} , \underline{\xi}$, en deux points
au voisinage de chaque point double de $C_{t_0}$ autre que $p$ et en $2 i $ (resp. $2i - 2$) points au voisinage de chaque pointe convergeant vers une orbite de Reeb prescrite
(resp. non prescrite) de multiplicit\'e $i$, ce qui d\'ecoule du Th\'eor\`eme $1.5$ de \cite{HWZ1}. D'apr\`es le Lemme \ref{lemmepointsdoubles} et (\ref{dimsphere}), cela ferait au total
$ r + 2r_L + 2k^2 - 4k + 2 +  2k -2v + 2\# \Gamma  =   r + 2r_L  + 2\# \Gamma + 2k^2 - 2k - 2v + 2 = 2k^2 + 1$ (resp. 
$ r + 2r_L + k^2 - 3k + 2 +  2k -2v + 2\# \Gamma  =   r + 2r_L  + 2\# \Gamma + k^2 - k - 2v + 2 = k^2 + 1$) dans le cas de la sph\`ere (resp. du plan projectif r\'eel), ce qui ne se peut pas
d'apr\`es la Remarque \ref{reminters}. La contradiction \`a laquelle nous venons d'aboutir provient du fait que le nombre de points
de notre configuration est strictement sup\'erieur \`a l'indice de Maslov de la courbe $C_t$ calcul\'e dans la Proposition \ref{propcotangent}. Ceci vaut \'egalement lorsque $L$ est un tore,
de sorte que nous aboutissons \`a la m\^eme conclusion pour peu qu'aucune des composantes de la courbe r\'eductible ne soit multiple. Ainsi, dans tous ces cas, pour toute structure presque-complexe 
$J_t$ proche de $J_{t_0}$,
\`a l'exclusion de $J_{t_0}$, et pour chaque point d'intersection r\'eel $p$ entre les deux composantes irr\'eductibles de $C_{t_0}$, il y a au plus une courbe $J_t$-holomorphe dans 
$T^* L$ satisfaisant nos conditions d'incidence. La d\'emonstration du fait qu'il y a au moins une  courbe $J_t$-holomorphe satisfaisant nos conditions d'incidence est la m\^eme que
celle de la Proposition $2.14$ de \cite{Wels1} et nous ne reproduisons pas ici cet argument local. Si une des deux composantes de la courbe est multiple, on applique le cas pr\'ec\'edent
au recoll\'e
d'un voisinage de la composante simple avec un rev\^etement d'un voisinage du cylindre multiple, lequel recoll\'e se projette sur un voisinage $C_{t_0}$ dans $T^* L$, pour aboutir \`a la m\^eme
conclusion. On proc\`ede de m\^eme si les deux composantes sont multiples de sorte que dans tous les cas, l'\'eventuelle d\'eg\'en\'erescence vers 
des courbes r\'eductibles ne fait pas obstacle \`a l'invariance de 
$F_{(r, r_L)} (\alpha , \beta , \Gamma , \underline{x} , \underline{\xi} , J)$. Il reste \`a \'etudier la possibilit\'e qu'une suite d'\'el\'ements de $\R {\cal M}_\gamma $ d\'eg\'en\`ere vers une courbe
\`a plusieurs \'etages. Cela ne se produit pas lorsque le chemin $\gamma $ est choisi de fa\c{c}on suffisamment g\'en\'erale. En effet, soit $D_X$ une composante non-triviale d'un \'etage
$\R \times S^* L$. Le nombre d'orbites de Reeb positives limites de $D_X$ compt\'ees avec multiplicit\'e moins le nombre d'orbite de Reeb n\'egatives compt\'ees avec multiplicit\'e vaut au moins
deux (resp. quatre) si $L$ est une sph\`ere (resp. un plan projectif r\'eel), ce qui d\'ecoule de la positivit\'e de l'aire $\int_{D_X} d\lambda$ et de l'isomorphisme $H_1 (S^* L ; \Z) \cong \Z / 2 \Z$
(resp. $H_1 (S^* L ; \Z) \cong \Z / 4 \Z$). De m\^eme, si $L$ est un tore, le nombre de pointes positives de $D_X$ est strictement plus grand que un.  La courbe \`a plusieurs \'etages que l'on 
consid\`ere poss\`ede une composante r\'eelle dans l'\'etage $T^* L$ que l'on note $D_L$. Soit $m$ le nombre de telles composantes
non-triviales $D_X$ de $\R \times S^* L$ adjacentes \`a $D_L$. D'apr\`es les Propositions \ref{propmaslov} et \ref{propcotangent}, le nombre d'asymptotes non-prescrites de $D_L$ est major\'e 
par $v-  \# \Gamma + m$. La dimension
virtuelle de l'espace des modules contenant $D_L$ est donc major\'ee par $2 v + \epsilon  k - 4m  - 1 - r - 2r_L - 2 \# \Gamma + 2m$ (resp. $2(v-m) - 1 - r - 2r_L $)
si $L$ est une sph\`ere ou un plan projectif r\'eel (resp. si $L$ est un tore). Dans tous les cas cette dimension est inf\'erieure \`a $-2$, ce qu'il fallait d\'emontrer. Ceci ach\`eve la d\'emonstration du
Th\'eor\`eme \ref{theoinv} dans le cas o\`u $n=2$. Lorsque $n=3$, le passage par un point critique ou la d\'{e}g\'{e}n\'{e}rescence vers une courbe multiple ou r\'eductible se traite 
\`a nouveau de la m\^eme
mani\`ere que dans le cas absolu \cite{Wels2}, \cite{Wels3}. Le seul ph\'enom\`ene nouveau \`a exclure est la d\'{e}g\'{e}n\'{e}rescence vers une courbe \`a plusieurs \'etages.
Lorsque $L$ est un tore, le nombre de pointes positives de chaque composante $D_X$ est \`a nouveau strictement plus grand que un ce qui force la dimension
virtuelle de l'espace des modules contenant $D_L$ \`a \^etre major\'ee par $-2$. Lorsque $L \in \{ S^3 , \R P^3 \}$, on a suppos\'e que $J$ est invariante par le flot de Reeb pour 
$\rho \gg 1$, de sorte que les espaces de modules contenant chaque composante $D_X$ de $S^* L \times \R$ sont munis d'une action de $\C^*$. Ceci force
la dimension
virtuelle de l'espace des modules contenant $D_L$ \`a chuter de deux de sorte que cette d\'{e}g\'{e}n\'{e}rescence en une courbe \`a plusieurs \'etages ne peut se produire en codimension un.
$\square$

\begin{rem}
Lorsque $L \in \{ S^3 , \R P^3 \}$, le Th\'eor\`eme \ref{theoinv} utilise une hypoth\`ese qui n'appara\^{\i}t pas en dimension deux, \`a savoir que la structure presque-complexe $J$ est invariante 
par le flot de Reeb pour $\rho \gg 1$. Cette hypoth\`ese semble n\'ecessaire en g\'en\'eral pour la raison suivante. Munissons la symplectisation $\R \times S^* L$ d'une structure presque-complexe
asymptotiquement cylindrique g\'en\'erique $J$. Pour tout $k$ strictement positif, cette symplectisation poss\`ede un cylindre $J$-holomorphe convergent positivement vers une orbite de Reeb 
rev\^etue $k+1$ fois (resp. $k+2$ fois)
et convergeant n\'egativement vers une orbite de Reeb rev\^etue $k$ fois si $L$ est une sph\`ere (resp. un espace projectif r\'eel) de dimension trois. La dimension attendue d'un tel cylindre vaut huit
d'apr\`es la Proposition \ref{propcotangent}. Fixons l'orbite de Reeb positive d'un tel cylindre. Sans l'hypoth\`ese d'invariance de $J$ par le flot de Reeb, l'orbite de Reeb n\'egative appartient
\`a un espace de dimension trois d'orbites tandis qu'avec cette condition, elle n'appartient qu'\`a un espace de dimension deux, puisque l'espace des modules de tels cylindres est muni d'une action
de $\R$ par translation dans le premier cas et de $\C^*$ dans le second, lesquelles actions pr\'eservent l'\'evaluation de l'orbite n\'egative dans son espace de dimension quatre d'orbites de Reeb.
Par cons\'equent, au-dessus d'un chemin g\'en\'erique de structures presque-complexes asymptotiquement cylindriques de $T^* L$, on ne peut \'eviter que les courbes $J$-holomorphes rigides
d'\'energie de Hofer finie ayant une orbite  de Reeb prescrite de multiplicit\'e plus grande que deux ou trois, selon que $L = S^3$ ou $\R P^3$, se brisent en courbes \`a deux \'etages dont l'\'etage
sup\'erieur poss\`ede un tel cylindre non-trivial ainsi que des cylindres triviaux. Le nombre de telles courbes n'est donc pas invariant. Il l'est si l'on se restreint aux structures invariantes par le
flot de Reeb \`a l'infini.
\end{rem}

\subsubsection{Quelques calculs}
\label{subsubsectcalculs}

L'entier $F_{(r, r_L)} (\alpha , \beta , \Gamma , \underline{x} , \underline{\xi} , J)$ \'etant ind\'ependant de $\Gamma , \underline{x} , \underline{\xi} , J$ d'apr\`es le Th\'eor\`eme \ref{theoinv},
nous le noterons $F_{(r, r_L)} (\alpha , \beta)$.
Afin d'all\'eger encore cette notation, nous noterons cet entier $F (\alpha , \beta)$ lorsque $r_L = 0$,  puisque la valeur de $r$ est alors d\'efinie sans ambigu\"{\i}t\'e par les calculs
de dimensions (\ref{dimsphere}) et  (\ref{dimtore}).

\begin{lemme}
\label{lemmecalc1}
Si $L$ est hom\'eomorphe \`a une sph\`ere de dimension deux et $r_L = 0$, on a
$F (e_1 , 0) = F (0 , e_1) = 1$, $F (e_2 , 0) = 2$, $F (0 , e_2) = 8$, $F (2e_1 , 0) = 2$,
$F (e_1 , e_1) = 4$ et $F (0 , 2e_1) = 6$.
\end{lemme}

{\bf D\'emonstration :}

D'apr\`es le Lemme \ref{lemmepointsdoubles}, les cylindres asymptotes \`a des orbites de Reeb simples sont plong\'es et il ressort de la Remarque \ref{reminters}
que deux tels cylindres 
s'intersectent en deux points au maximum. Par suite, $F_{(3,0)} (0 , e_1) = 1$ et $F_{(1,0)} (e_1 , 0) = F_{(1,1)} (0 , e_1) = 1$. De m\^eme, $F_{(1,1)} (e_2 , 0) = 0$ puisqu'en faisant tendre la paire
de points complexes conjugu\'es vers l'infini les courbes devraient converger vers des courbes \`a deux \'etages non-triviales et n'ayant qu'une orbite double comme limite positive. De telles
courbes \`a deux \'etages n'existent pas. Pour comparer $F_{(1,1)} (e_2 , 0) $ \`a $F_{(3,0)} (e_2 , 0)$, on proc\`ede comme au \S $3$ de \cite{Wels1} en faisant tendre la
paire de points complexes conjugu\'es vers un point $y$ de $L$. Il n'y a comme pr\'ec\'edemment qu'un seul cylindre asymptote \`a une orbite double, passant par un point $x$ de $L$
et ayant un point double en $y$. On montre de la m\^eme mani\`ere que le Th\'eor\`eme $3.2$ de \cite{Wels1} la relation $F_{(3,0)} (e_2 , 0) - F_{(1,1)} (e_2 , 0) = 2$, d'o\`u 
on d\'eduit $F_{(3,0)} (e_2 , 0)  = 2$. On \'etablit de m\^eme la relation $F_{(5,0)}  (0 , e_2) = 2F_{(3 + \times ,0)}  (0 , e_2) + F_{(3,1)}  (0 , e_2) $ o\`u $F_{(3 + \times ,0)}  (0 , e_2) $ d\'esigne
le nombre alg\'ebrique de cylindres $J_L$-holomorphes r\'eels de $T^* L$ asymptotes \`a une orbite de Reeb double non prescrite passant par trois points de $L$ et ayant son point double  impos\'e
en un quatri\`eme point, voir le $\S 3.1$ de \cite{Wels1}. On a $F_{(3,1)}  (0 , e_2) = 2 F_{(3,0)}  (e_2 , 0) = 4$. De m\^eme, $ F_{(3 + \times ,0)}  (0 , e_2) = F_{(1 + \times ,1)}  (0 , e_2) =
2 F_{(1 + \times ,0)}  (e_2, 0) = 2$. D'o\`u finalement $F_{(5,0)}  (0 , e_2) = 4 + 4 =8$. Enfin,  $F_{(3,0)} (2e_1 , 0) = 2F_{(1 + \times ,0)} (2e_1 , 0) + F_{(1,1)} (2e_1 , 0) = 2 + 0 = 2 $,
$F_{(5,0)} (e_1 , e_1) = 2 F_{(3 + \times ,0)} (e_1 , e_1) + F_{(3 , 1)} (e_1 , e_1) = 2 + F_{(3,0)} (2e_1 , 0) = 4$, $F_{(7,0)} (0, 2e_1) = 2F_{(5 + \times ,0)} (0, 2e_1) + F_{(5,1)} (0, 2e_1) = 
2 + F_{(5,0)} (e_1 , e_1) = 6$. $\square$

\begin{lemme}
\label{lemmecalc2}
Si $L$ est hom\'eomorphe \`a un plan projectif r\'eel  et $r_L = 0$, on a
$F (e_1 , 0) = F (0 , e_1) = F (e_2 , 0) = F (2e_1 , 0) = F (e_1 , e_1)  = F (0 , 2e_1) = 1$
et $F (0 , e_2 ) = 4$.
\end{lemme}

{\bf D\'emonstration :}

D'apr\`es le Lemme \ref{lemmepointsdoubles}, les cylindres asymptotes \`a des orbites de Reeb simples sont plong\'es et d'apr\`es la Remarque \ref{reminters}, deux tels cylindres 
s'intersectent en un point au maximum. Par suite, $F_{(2,0)} (0 , e_1) = F_{(0,0)} (e_1 , 0) = 1$. De m\^eme,  des cylindres asymptotes \`a des orbites de Reeb doubles ou des sph\`eres
ayant quatre pointes asymptotes \`a des orbites de Reeb simples s'intersectent en quatre points au maximum, de sorte que 
$F_{(1,0)} (e_2 , 0) = F_{(1,0)} (2e_1 , 0) = F_{(3,0)} (e_1 , e_1)  = F _{(5,0)} (0 , 2e_1) = 1$. Enfin, $F_{(1,1)}  (0 , e_2 ) = 2F_{(1,0)}  (e_2 , 0 ) = 2$. Lorsque l'on fait converger la paire
de points complexes conjugu\'es vers un point $y$ de $L$, les deux courbes r\'eelles compt\'ees par $F_{(1,1)}  (0 , e_2 )$ convergent vers deux courbes r\'eelles ayant une tangente prescrite
en $y$. En effet, ces derni\`eres ne peuvent converger vers des courbes ayant un point double en $y$ d'apr\`es le  Lemme \ref{lemmepointsdoubles} et le rev\^etement  double du cylindre
sur une orbite simple passant par $y$ ne peut se d\'eformer en un cylindre interpolant une paire de points complexes conjugu\'es. Par contre, ce cylindre double se d\'eforme en  au moins
un cylindre asymptote \`a une paire d'orbites doubles et passant par trois points de $L$, de sorte que $F (0 , e_2 ) \geq 4$. Or, le nombre de cylindres complexes de $T^* L$ asymptotes
\`a deux orbites doubles et passant par trois points vaut quatre, ce que l'on obtient en faisant tendre deux points vers l'infini. D'o\`u le r\'esultat. $\square$

\begin{lemme}
\label{lemmecalc3}
Si $L$ est hom\'eomorphe \`a un plan projectif r\'eel  et $r_L = 0$, on a 
$F (e_3 , 0) = 2$, $F (0, e_3) = 12$, $F (e_1 + e_2 , 0) = 2$, $F (e_1 , e_2) = 8$,
$F (e_2 , e_1) = 4$, $F (0 , e_1 + e_2) = 24$, $F (3e_1 , 0) = 2$, $F (2e_1 , e_1) = 4$,
 $F (e_1 , 2e_1) = 6$ et $F (0 , 3e_1) = 8$.
\end{lemme}

{\bf D\'emonstration :}

On proc\`ede comme dans la d\'emonstration du Lemme \ref{lemmecalc1}. On obtient avec les m\^emes notations $F_{(2,0)}  (e_3 , 0) = 2 F_{( \times ,0)}  (e_3 , 0) + 
F_{(0, 1)}  (e_3 , 0) = 2 + 0 = 2$. 
De m\^eme, $F_{(4,0)} (0, e_3) =  2 F_{( 2 + \times ,0)}  (0 , e_3 ) + F_{(2, 1)}  (0 , e_3) = 2 F_{( 2 + \times ,0)}  (0 , e_3 ) + 3F_{(2, 0)}  ( e_3 , 0)  = 2 F_{( 2 + \times ,0)}  (0 , e_3 ) + 6$
et $F_{( 2 + \times ,0)}  (0 , e_3 ) = F_{( \times ,1)}  (0 , e_3 ) = 3 F_{( \times ,0)}  (e_3 , 0 ) = 3$, de sorte que $F_{(4,0)} (0, e_3) = 12$. De m\^eme, $F_{(2,0)}  (e_1 + e_2 , 0) =
2 F_{(\times,0)}  (e_1 + e_2 , 0) + F_{(0 , 1)}  (e_1 + e_2 , 0) = 2 + 0 = 2$ ;  $F_{(4,0)}  (e_1, e_2 ) =
2 F_{(2 + \times,0)}  (e_1 , e_2) + F_{(2 , 1)}  (e_1 , e_2 ) = 2 F_{( \times,1)}  (e_1 , e_2) + 2 F_{(2 , 0)}  (e_1 + e_2 , 0 ) =  4 F_{( \times,0)}  (e_1 + e_2 , 0) + 4 = 8$ ;
$F_{(4,0)}  (e_2, e_1 ) = 2 F_{(2 + \times,0)}  (e_2 , e_1) + F_{(2 , 1)}  (e_2 , e_1) =  2 + F_{(2 , 0)}  (e_1 + e_2 , 0 ) = 4$ ; $F_{(6,0)} (0 , e_1 + e_2) =  2 F_{(4 + \times,0)}  (0 , e_1 + e_2) 
+ F_{(4, 1)} (0 , e_1 + e_2) = 2 F_{(2 + \times,1)}  (0 , e_1 + e_2) + F_{(4, 0)} (e_1 , e_2) + 2F_{(4, 0)} (e_2 , e_1) = 2 F_{(2 + \times,0)}  (e_1 , e_2) + 4 F_{(2 + \times,0)}  (e_2 , e_1) +
8 + 8 = 4 + 4 + 16 = 24$. Enfin, $F_{(2,0)} (3e_1 , 0) = 2 F_{(\times,0)} (3e_1 , 0) + F_{(0,1)} (3e_1 , 0) = 2 + 0 = 2$ ; $F_{(4,0)} (2e_1 , e_1) = 2 F_{(2 + \times,0)} (2e_1 , e_1) + F_{(2,1)} (2e_1 , e_1) 
= 2 + F_{(2,0)} (3e_1 , 0) = 4$ ; $F_{(6,0)} (e_1 , 2e_1) = 2 F_{(4 + \times,0)} (e_1 , 2e_1) + F_{(4,1)} (e_1 , 2e_1) = 2 + F_{(4,0)} (2e_1 , e_1) = 6$ ; 
$F_{(8,0)} (0 , 3e_1) = 2 F_{(6 + \times,0)} (0 , 3e_1) + F_{(6,1)} (0 , 3e_1) = 2 + F_{(6,0)} (e_1 , 2e_1) = 8$. $\square$

\subsection{Calculs dans le plan projectif complexe}

\subsubsection{Arbres projectifs}

Soient $r, r_X$ et $d$ trois entiers naturels satisfaisant la relation $r + 2r_X = 3d - 1$.

\begin{defi}
Un arbre projectif est un arbre fini connexe dont toutes les ar\^etes sont \'etiquet\'ees par des entiers strictement positifs. De plus, un tel arbre poss\`ede une racine $s_0$ et tous 
les sommets \`a distance paire
de $s_0$ sont soit monovalents connect\'es \`a une ar\^ete double, soit bivalents connect\'es \`a deux ar\^etes simples.
\end{defi}

On note ${\cal B}_r$ l'ensemble des arbres projectifs qui satisfont
$$\sum_{a \in {\cal A}(s_0)} k(a) - 1 \leq r \leq  \sum_{a \in {\cal A}(s_0)} k(a) - 1 + 2v (s_0) \text{ et } r = \sum_{a \in {\cal A}(s_0)} k(a) - 1 \mod (2),$$
o\`u $v(s)$ d\'esigne la valence d'un sommet $s$, $ {\cal A}(s)$ l'ensemble des ar\^etes adjacentes \`a $s$ et $k (a)$ d\'esigne la multiplicit\'e de l'ar\^ete $a$.
Notons \'egalement $k_s$ la somme des multiplicit\'es des ar\^etes adjacentes au sommet $s$ et $k$ la multiplicit\'e totale de toutes les ar\^etes de l'arbre.
On pose $r_L (s_0) = \frac{1}{2} \big(    \sum_{a \in {\cal A}(s_0)} k(a) - 1 + 2v (s_0) - r     \big)$, de sorte que $0 \leq r_L (s_0) \leq v (s_0)$. Enfin, on note $S_1$ (resp. $S_2$) l'ensemble
des sommets \`a distance impaire (resp. paire) de $s_0$.

\begin{defi}
Un arbre projectif d\'ecor\'e est un arbre projectif $A \in {\cal B}_r$ \'{e}quip\'{e} d'une partition $S_1^+ \sqcup S_1^-$ de l'ensemble des sommets adjacents \`a $s_0$ telle que
$\# S_1^- = r_L (s_0)$ et $\# S_1^+ = v(s_0) - r_L (s_0)$. Cet arbre est de plus \'equip\'e des fonctions :
\begin{itemize}
\item $f_A : S_1 \to {\cal P} (\{1, \dots , r_X \})$ satisfaisant $f_A (s) \cap f_A (s') = \emptyset$ d\`es que $s \neq s'$ et $\cup_{s \in S_1} f_A (s) = \{1, \dots , r_X \}$.
\item  $g_A : S_1 \to \N$ telle que $k + 4\sum_{s \in S_1} g_A (s) = d$ et pour tout $s \in S_1^+$ (resp. $s \in S_1^-$ ), $6g_A (s) + k_s + v(s) - 1 = \# f_A (s) + 1$
(resp. $6g_A (s) + k_s + v(s) - 1 = \# f_A (s)$).
\end{itemize}
\end{defi}

On note ${\cal B}_r^d$ l'ensemble fini des arbres projectifs d\'ecor\'es. Soit $A \in {\cal B}_r^d$, on pose $m_1^- (A) = \prod_{s \in S_1^-} \#{\{ a \in  {\cal A}(s)} \; \vert \; k(a) = k(ss_0)\}$,
o\`u $k(ss_0)$ d\'esigne la multiplicit\'e de l'ar\^ete reliant $s$ \`a $s_0$. On note de m\^eme $m_1^+ (A)$ le nombre d'injections 
$\phi : \{s \in S_1^+ \; \vert \; f_A (s) \neq \emptyset \} \to   {\cal A}^+ (s_0)$ satisfaisant
$k (\phi (s)) = k(ss_0)$ pour tous les sommets $s \in S_1^+$, o\`u ${\cal A}^+ (s_0)$ d\'esigne l'ensemble des ar\^etes reliant $s_0$ \`a un sommet de $S_1^+$. Notons enfin $S_2^b$ l'ensemble des sommets bivalents de $S_2 \setminus \{ s_0 \}$. L'arbre $A$ priv\'e de ces sommets bivalents
n'est pas connexe. On note $m_2 (A)$ le nombre de fa\c{c}ons de reconnecter $A \setminus S_2^b$ de mani\`ere \`a obtenir un arbre isomorphe \`a $A$.  Posons
$$\text{mult} (A) = 2^{\sum_{s \in S_1^+} \# f_A (s) + \sum_{s \in S_1 \setminus S_1^+} \max (\# f_A (s) - 1 , 0) + \#S_2^b} m_1^+ (A) m_1^- (A) m_2 (A) \prod_a k(a).$$
C'est la {\it multiplicit\'e} de l'arbre $A \in {\cal B}_r^d$.

\begin{figure}[htb]
\label{figarbresprojectifs}
\begin{center}
\includegraphics{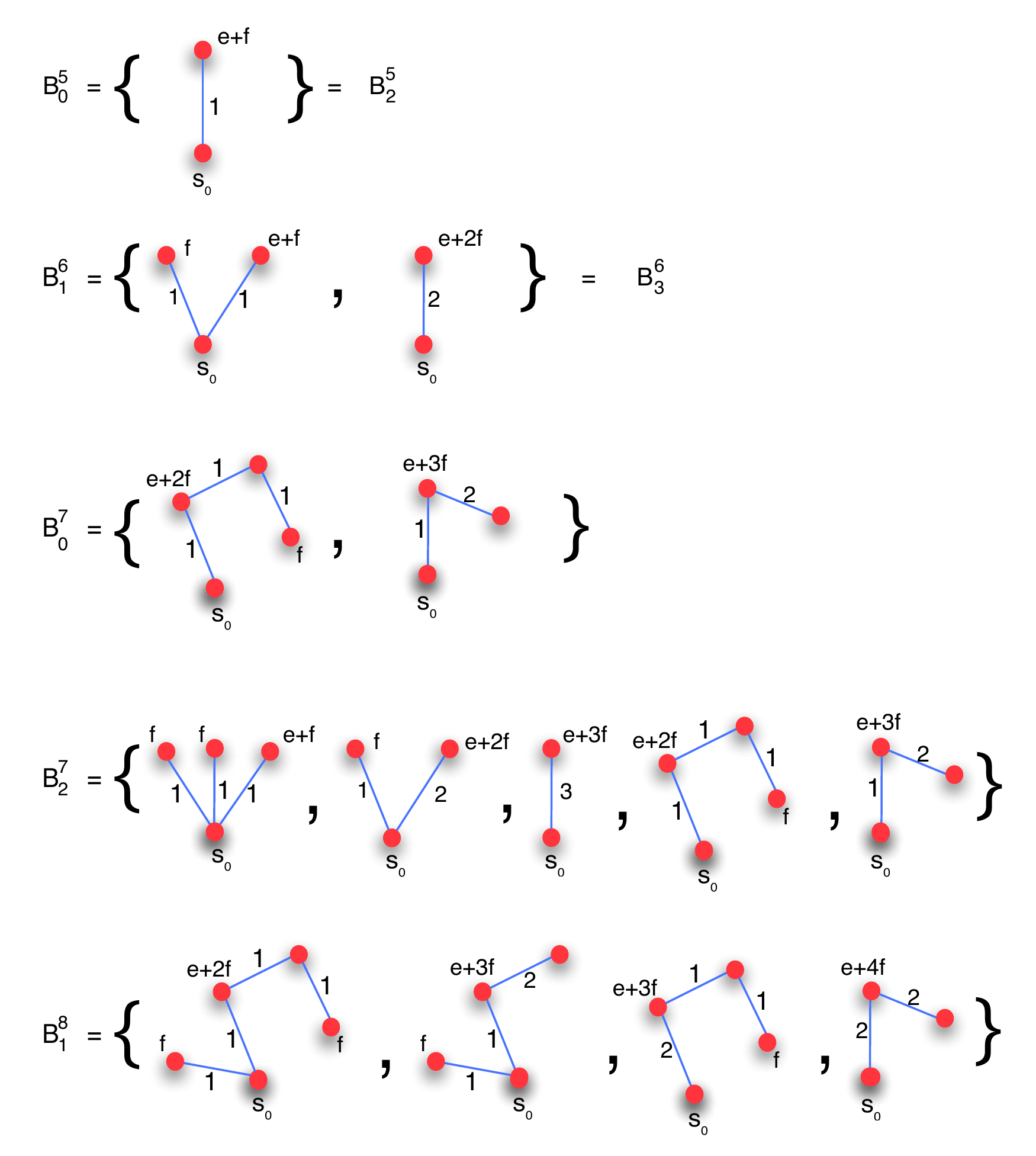}
\end{center}
\caption{Arbres projectifs d\'ecor\'es}
\end{figure}

\subsubsection{Relation avec l'invariant relatif}

Soit $\Sigma_4$ la surface rationnelle r\'egl\'ee de degr\'e quatre, $e$ la classe d'une section holomorphe d'autointersection quatre et $f$ la classe d'une fibre. \'Etant donn\'es $a,b \in \N$
et $\alpha, \beta$ des suites d'entiers positifs, on note $N_4^{ae+bf} (\alpha, \beta)$ le nombre de courbes rationnelles de $\Sigma_4$, homologues \`a $ae+bf$, ayant $\alpha_i + \beta_i$
points de tangence d'ordre $i$ avec la section exceptionnelle de $\Sigma_4$ parmi lesquels $\alpha_i $ sont prescrits et qui passent par le nombre ad\'equat de points fix\'es.

\begin{theo}
\label{theocalcproj}
Soit $(X, \omega , c_X)$ une vari\'et\'e symplectomorphe au plan projectif complexe et $r, r_X, d \in \N$ satisfaisant la relation $r + 2r_X = 3d - 1$. Alors,
$$\chi^d_r = \sum_{A \in {\cal B}_r^d} (-1)^{\# S_2 + 1} \text{mult} (A) F_{(r,0)} (\alpha_A^-, \beta_A^+) \prod_{s \in S_1 \setminus S_1^+} N_4^{g_A (s) e + k(s) f} (0, \beta_A)
\prod_{s \in S_1^+} N_4^{g_A (s) e + k(s) f} (e_{k(ss_0)} , \beta_A^0),$$
o\`u $(\alpha_A^-)_i$ (resp. $(\beta_A^+)_i$) vaut le  nombre d'ar\^etes de multiplicit\'e $i$ reliant $S_1^-$ (resp. $S_1^+$) \`a $s_0$, $( \beta_A)_i$ (resp. $( \beta_A^0)_i$)
vaut le nombre d'ar\^etes de multiplicit\'e $i$  adjacentes \`a $s$ (resp. moins un si l'ar\^ete reliant $s$ \`a $s_0$ est de multiplicit\'e $i$) et $F_{(r,0)} (\alpha_A^-, \beta_A^+)$ d\'esigne
l'invariant d\'efini dans le fibr\'e cotangent du plan projectif r\'eel au \S \ref{subsubinvariants}.
\end{theo}

\begin{rem}
Un algorithme permettant le calcul de $\chi^d_{3d-1}$ a d\'ej\`a \'et\'e propos\'e par G. Mikhalkin dans \cite{Mikh} et \'etendu par E. Shustin dans \cite{Shu} aux autres invariants $\chi^d_r$,
$0 \leq r \leq 3d-1$.
Par ailleurs, I. Itenberg, V. Kharlamov et E. Shustin \cite{IKS2} ont adapt\'e la formule de r\'ecurrence  \cite{Gath} au cas r\'eel et ainsi d\'eduit une formule calculant l'invariant $\chi^d_{3d-1}$ 
en fonction d'invariants relatifs tropicaux dans les surfaces de Del Pezzo toriques r\'eelles. J'avais introduit des invariants relatifs r\'eels par rapport \`a une courbe r\'eelle ayant une
partie r\'eelle non-vide dans \cite{Wels5}. Les formules apparaissant dans les Th\'eor\`emes \ref{theocalcproj}, \ref{theocal2spher} et \ref{theocal3spher} calculent $\chi^d_r$ en fonction 
d'invariants relatifs \`a un diviseur r\'eel ayant une partie r\'eelle vide ; il serait int\'eressant de calculer ces invariants $F_{(r, r_L)} (\alpha , \beta )$ et de d\'eterminer exactement 
dans quelles situations l'expression de $\chi^d_r$ en fonction de tels invariants relatifs \`a un diviseur r\'eel de partie r\'eelle vide permet le calcul effectif de $\chi^d_r$.
Lorsque $r \leq 2$, ces invariants $F_{(r, r_L)} (\alpha , \beta )$ sont calcul\'es au \S \ref{subsubsectcalculs} tandis que $N_4^{ae+bf} (\alpha, \beta)$ est calcul\'e dans \cite{Vak}, de sorte que
le Th\'eor\`eme \ref{theocalcproj} permet le calcul effectif de $\chi^d_r$, voir le Corollaire \ref{corcalcproj} pour les premi\`eres valeurs de cet invariant.
\end{rem}

{\bf D\'emonstration du Th\'eor\`eme \ref{theocalcproj} :}

On poursuit la strat\'egie g\'en\'erale \'enonc\'ee au $\S$ \ref{subsubsectstrat} en allongeant le cou d'une structure presque complexe g\'en\'erique jusqu'\`a briser la vari\'et\'e en deux morceaux
$T^* \R P^2$ et $\C P^2 \setminus \R P^2$. Les courbes $J$-holomorphes rationnelles r\'eelles compt\'ees par $\chi^d_r$ se brisent en courbes \`a deux \'etages interpolant  $r$ points de $\R P^2$
et $r_X$ paires de points complexes conjugu\'es de $\C P^2 \setminus \R P^2$. Il est apparu au cours de la d\'emonstration du Th\'eor\`eme \ref{theocong3} que ces courbes \`a deux \'etages
sont cod\'ees par les arbres projectifs d\'ecor\'es $A \in {\cal B}_r^d$. Il s'agit donc de d\'enombrer les  courbes \`a deux \'etages qui sont cod\'ees par un arbre donn\'e $A \in {\cal B}_r^d$, puis
de d\'enombrer les  courbes $J$-holomorphes rationnelles r\'eelles compt\'ees par $\chi^d_r$ qui d\'eg\'en\`erent sur une courbe \`a deux \'etages donn\'ee. Le nombre de fa\c{c}ons de r\'epartir les
points complexes conjugu\'es parmi les composantes de  la courbe \`a deux \'etages qui se trouvent dans $\C P^2 \setminus \R P^2$ a \'et\'e calcul\'e dans la d\'emonstration du Th\'eor\`eme
\ref{theocong3} et vaut $2^{\sum_{s \in S_1^+} \# f_A (s) + \sum_{s \in S_1 \setminus S_1^+} \max (\# f_A (s) - 1 , 0)}$. Les composantes cod\'ees par $S_1 \setminus S_1^+$ sont rigides
avec leurs conditions d'incidence, il y en a $\prod_{s \in S_1 \setminus S_1^+} N_4^{g_A (s) e + k(s) f} (0, \beta_A)$. Puis, il y a  $m_1^- (A)$ fa\c{c}ons de choisir les orbites de Reeb prescrites
de la courbe cod\'ee par $s_0$. Le nombre de courbes r\'eelles cod\'ees par $s_0$ satisfaisant nos conditions d'incidence et compt\'ees avec signe vaut $F_{(r,0)} (\alpha_A^-, \beta_A^+)$.
Il y a alors $m_1^+ (A)$  fa\c{c}ons de choisir la mani\`ere de connecter les courbes cod\'ees par $S_1^+$ aux  orbites de Reeb rest\'ees libres de la courbe cod\'ee par $s_0$. Il y a enfin
$2^{ \#S_2^b} m_2(A)$ fa\c{c}ons de connecter ces composantes entre elles par des paires de cylindres complexes conjugu\'es de $T^* \R P^2$ cod\'es par les sommets bivalents de $S_2^b$.
Ceci fournit $2^{\sum_{s \in S_1^+} \# f_A (s) + \sum_{s \in S_1 \setminus S_1^+} \max (\# f_A (s) - 1 , 0) + \#S_2^b} m_1^+  m_1^- m_2 F_{(r,0)} (\alpha_A^-, \beta_A^+)$ $\prod_{s \in S_1 
\setminus S_1^+} N_4^{g_A (s) e + k(s) f} (0, \beta_A) \prod_{s \in S_1^+} N_4^{g_A (s) e + k(s) f} (k(ss_0) , \beta_A^0)$ courbes cod\'ees par un arbre donn\'e $A \in {\cal B}_r^d$. 
Or, d'apr\`es le th\'eor\`eme de
recollement de th\'eorie symplectique des champs \cite{Bour}, il y a $ \prod_a k(a)$ courbes $J$-holomorphes rationnelles r\'eelles compt\'ees par $\chi^d_r$ qui d\'eg\'en\`erent sur une courbe \`a deux 
\'etages donn\'ee. Le r\'esultat d\'ecoule \`a pr\'esent du fait que chaque courbe cod\'ee par $S_2 \setminus \{ s_0 \}$ intersecte $\R P^2$ en un point et contribue donc \`a la masse des courbes
$J$-holomorphes rationnelles r\'eelles en question, d'o\`u le signe $(-1)^{\# S_2 + 1}$. $\square$

\begin{cor}
\label{corcalcproj}
Soit $(X, \omega , c_X)$ une vari\'et\'e symplectomorphe au plan projectif complexe. Alors,
$\chi^4 (T) = o(T^2)$, $\chi^5 (T) = 64 + 64T^2 + o(T^3)$,  $\chi^6 (T) = 1024T + 1536T^3 + o(T^4)$,
$\chi^7 (T) = -14336  + 11776T^2 + o(T^3)$ et $\chi^8 (T) = -280576T + o(T^2)$.
\end{cor}

\begin{rem}
Les valeurs $\chi^4_1 = 0$, $\chi^5_0 = \chi^5_2 = 64$ ont d\'ej\`a \'et\'e obtenues dans  \cite{IKS1} \`a l'aide d'un ordinateur et de l'algorithme \cite{Shu}. 
\end{rem}

{\bf D\'emonstration du Corollaire \ref{corcalcproj} :}

L'annulation de $\chi^4_1$ tient au fait que l'ensemble d'arbres $B^4_1$ est vide. Les arbres intervenant dans la d\'emonstration de ce Corollaire \ref{corcalcproj} sont
repr\'esent\'es dans la Figure \ref{figarbresprojectifs}. Lorsque $d=5$ et $r \leq 2$, un seul arbre d\'ecor\'e intervient. Le Th\'eor\`eme \ref{theocalcproj} fournit
$\chi^5_0 = 2^6 F (e_1 , 0) N^{e+f} (0 , e_1) = 64$ et $\chi^5_2 = 2^6 F (0 , e_1) N^{e+f} (e_1 , 0) = 64$. Lorsque $d=6$ et $r=1$, deux arbres projectifs d\'ecor\'es interviennent
qui sont repr\'esent\'es par la Figure \ref{figarbresprojectifs}.
La contribution du premier vaut $C_8^1 2^6 F (2e_1 , 0) N^{e+f} (0 , e_1) =  512$ et celle du second donn\'e par cette figure vaut $2^7 2 F (e_2 , 0) N^{e+2f} (0 , e_2) = 2^9 
N^{e+2f} (e_2 , 0) = 512$, voir le Lemme \ref{lemmecalc2}, de sorte
que $\chi^6_1 = 1024$. Lorsque $d=6$ et $r=3$, les arbres intervenant sont les m\^emes. Toutefois, la fonction $f_A$ du premier arbre peut affecter soit six, soit sept paires de points complexes
conjugu\'es au sommet $e+f$. Dans le premier cas, l'arbre contribue \`a hauteur de $C_7^1 2^6 F (e_1 , e_1) N^f  (0 , e_1) N^{e+f} (e_1 , 0) =  448$ ; dans le second, il contribue
\`a hauteur de $2^6 F (e_1 , e_1) N^f  (e_1 , 0) N^{e+f} (0 , e_1) = 64$, soit une contribution totale de $512$. La contribution du second arbre vaut
$2^7 2 F (0 , e_2) N^{e+2f} (e_2 , 0) = 2^{10}$, de sorte que $\chi^6_3 = 1536$. Lorsque $d=7$ et $r=0$, deux arbres contribuent. Le premier donn\'e par la Figure \ref{figarbresprojectifs}
contribue \`a hauteur de $- C_{10}^1 2^9 2 F (e_1 , 0) N^{e+2f} (0 , 2e_1) = -10240$ tandis que le second contribue \`a hauteur de $- 2^9 2 F (e_1 , 0) N^{e+3f} (0 , e_1 + e_2)$. Or
$N^{e+3f} (0 , e_1 + e_2) = N^{e+3f} (e_1 , e_2) + 2N^{e+3f} (e_2 , e_1) = 4 N^{e+3f} (e_1 + e_2 , 0) = 4$, de sorte que finalement $\chi^7_0 = -10240 - 4096 = -14336$. Lorsque $d=7$ et $r=2$, 
cinq arbres contribuent. La contribution du premier arbre donn\'e par la Figure \ref{figarbresprojectifs} vaut $C_9^2 2^6 F (3e_1 , 0) N^{e+f} (0 , e_1) = 4608$. La contribution du deuxi\`eme 
arbre vaut $C_9^1 2^7 2 F (e_1 + e_2 , 0) N^{e+2f} (0 , e_2) = 
9 2^9 2 N^{e+2f} (e_2 , 0) = 9216$. La contribution du troisi\`eme arbre vaut $2^8 3 F (e_3 , 0) N^{e+3f} (0 , e_3) = 9 2^9 N^{e+3f} (e_3 , 0) = 4608$.  La contribution du quatri\`eme 
arbre vaut 
$-C_9^1 2^9 F (0 , e_1) N^{e+2f} (e_1 , e_1) N^f (0 , e_1) = -4608$. Enfin, la contribution du cinqui\`eme arbre vaut $-2^9 2 F (0 ,e_1) N^{e+3f} (e_1 , e_2) = -2^{11} N^{e+3f} (e_1 + e_2 , 0) = -2048$.
On en d\'eduit $\chi^7_2 = 4608 + 9216 + 4608 - 4608 - 2048 = 11776$. Lorsque $d=8$ et $r=1$, quatre arbres projectifs d\'ecor\'es contribuent. La contribution du premier arbre donn\'e
par la Figure \ref{figarbresprojectifs} vaut $- 2 C_{11}^2 2^9 2 F (2e_1 , 0) 
N^{e+2f} (0 , 2e_1) N^f  (0 , e_1)^2 = -112640$.  La contribution du deuxi\`eme arbre vaut $-C_{11}^1 2^9 2 F (2e_1 , 0) N^{e+3f} (0 , e_1 + e_2) N^f (0 , e_1) = -45056$.  La contribution du troisi\`eme arbre vaut 
$-C_{11}^1 2^{10} 2 F (e_2 , 0) N^{e+3f} (0 , e_1 + e_2) N^f  (0 , e_1) = -90112$. Enfin, la contribution du quatri\`eme arbre vaut $-2 2^{10} 2^2 F (e_2 , 0) N^{e+4f} (0 , 2e_2) = -2^{15} N^{e+4f} (2e_2 , 0)
= -32768$. On en conclut $\chi^8_1 = -112640 - 45056 - 90112 - 32768 = -280576$. $\square$

\subsection{Calculs dans l'ellipso\"{\i}de de dimension deux}

\subsubsection{Arbres deux-sph\'eriques}

Soient $r, r_X$ et $d$ trois entiers naturels satisfaisant la relation $r + 2r_X = 4d - 1$, laquelle impose que $r$ soit impair.

\begin{defi}
Un arbre deux-sph\'erique est un arbre fini connexe dont toutes les ar\^etes sont \'etiquet\'ees par des entiers strictement positifs. De plus, un tel arbre poss\`ede une racine $s_0$ et tous les sommets \`a 
distance paire de $s_0$ sont monovalents connect\'es \`a une ar\^ete simple.
\end{defi}

En particulier, la distance d'un sommet \`a $s_0$ est major\'ee par deux et seules les ar\^etes connect\'ees \`a $s_0$ peuvent avoir une multiplicit\'e non-triviale.
On note ${\cal A}_r$ l'ensemble des arbres deux-sph\'eriques qui satisfont
$$2\sum_{a \in {\cal A}(s_0)} k(a) - 1 \leq r \leq  2\sum_{a \in {\cal A}(s_0)} k(a) - 1 + 2v (s_0),$$
o\`u $v(s)$ d\'esigne la valence d'un sommet $s$, $ {\cal A}(s)$ l'ensemble des ar\^etes adjacentes \`a $s$ et $k (a)$ d\'esigne la multiplicit\'e de l'ar\^ete $a$.
Notons \'egalement $k_s$ la somme des multiplicit\'es des ar\^etes adjacentes au sommet $s$ et $k$ la multiplicit\'e totale de toutes les ar\^etes de l'arbre.
On pose $r_L (s_0) = \frac{1}{2} \big(    2\sum_{a \in {\cal A}(s_0)} k(a) - 1 + 2v (s_0) - r     \big)$, de sorte que $0 \leq r_L (s_0) \leq v (s_0)$. Enfin, on note $S_1$ (resp. $S_2$) l'ensemble
des sommets \`a distance impaire (resp. paire) de $s_0$.

\begin{defi}
Un arbre deux-sph\'erique d\'ecor\'e est un arbre deux-sph\'erique $A \in {\cal A}_r$ \'{e}quip\'{e} d'une partition $S_1^+ \sqcup S_1^-$ de l'ensemble des sommets adjacents \`a $s_0$ telle que
$\# S_1^- = r_L (s_0)$ et $\# S_1^+ = v(s_0) - r_L (s_0)$. Cet arbre est de plus \'equip\'e des fonctions :
\begin{itemize}
\item $f_A : S_1 \to {\cal P} (\{1, \dots , r_X \})$ satisfaisant $f_A (s) \cap f_A (s') = \emptyset$ d\`es que $s \neq s'$ et $\cup_{s \in S_1} f_A (s) = \{1, \dots , r_X \}$.
\item  $g_A : S_1 \to \N$ telle que $k + 2\sum_{s \in S_1} g_A (s) = d$ et pour tout $s \in S_1^+$ (resp. $s \in S_1^-$ ), $4g_A (s) + k_s + v(s) - 1 = \# f_A (s) + 1$
(resp. $4g_A (s) + k_s + v(s) - 1 = \# f_A (s)$).
\end{itemize}
\end{defi}

On note ${\cal A}_r^d$ l'ensemble des arbres deux-sph\'eriques d\'ecor\'es, c'est un ensemble fini. Soit $A \in {\cal A}_r^d$, on pose 
$m_1^- (A) = \prod_{s \in S_1^-} \#{\{ a \in  {\cal A}(s)} \; \vert \; k(a) = k(ss_0)\}$, o\`u $k(ss_0)$ d\'esigne la multiplicit\'e de l'ar\^ete reliant $s$ \`a $s_0$, de sorte que chaque terme du 
produit vaille un ou $v(s)$. On note de m\^eme $m_1^+ (A)$ le nombre d'injections $\phi : \{s \in S_1^+ \; \vert \; f_A (s) \neq \emptyset \} \to  {\cal A}^+ (s_0)$ satisfaisant
$k (\phi (s)) = k(ss_0)$ pour tous les sommets $s \in S_1^+$.   On pose alors 
$$\text{mult} (A) = 2^{\sum_{s \in S_1^+} \# f_A (s) + \sum_{s \in S_1 \setminus S_1^+} \max (\# f_A (s) - 1 , 0) + \#S_2 -1} m_1^+ (A) m_1^- (A) \prod_a k(a),$$
c'est la {\it multiplicit\'e} de l'arbre $A \in {\cal A}_r^d$.

\begin{figure}[htb]
\label{figarbres2spheriques}
\begin{center}
\includegraphics{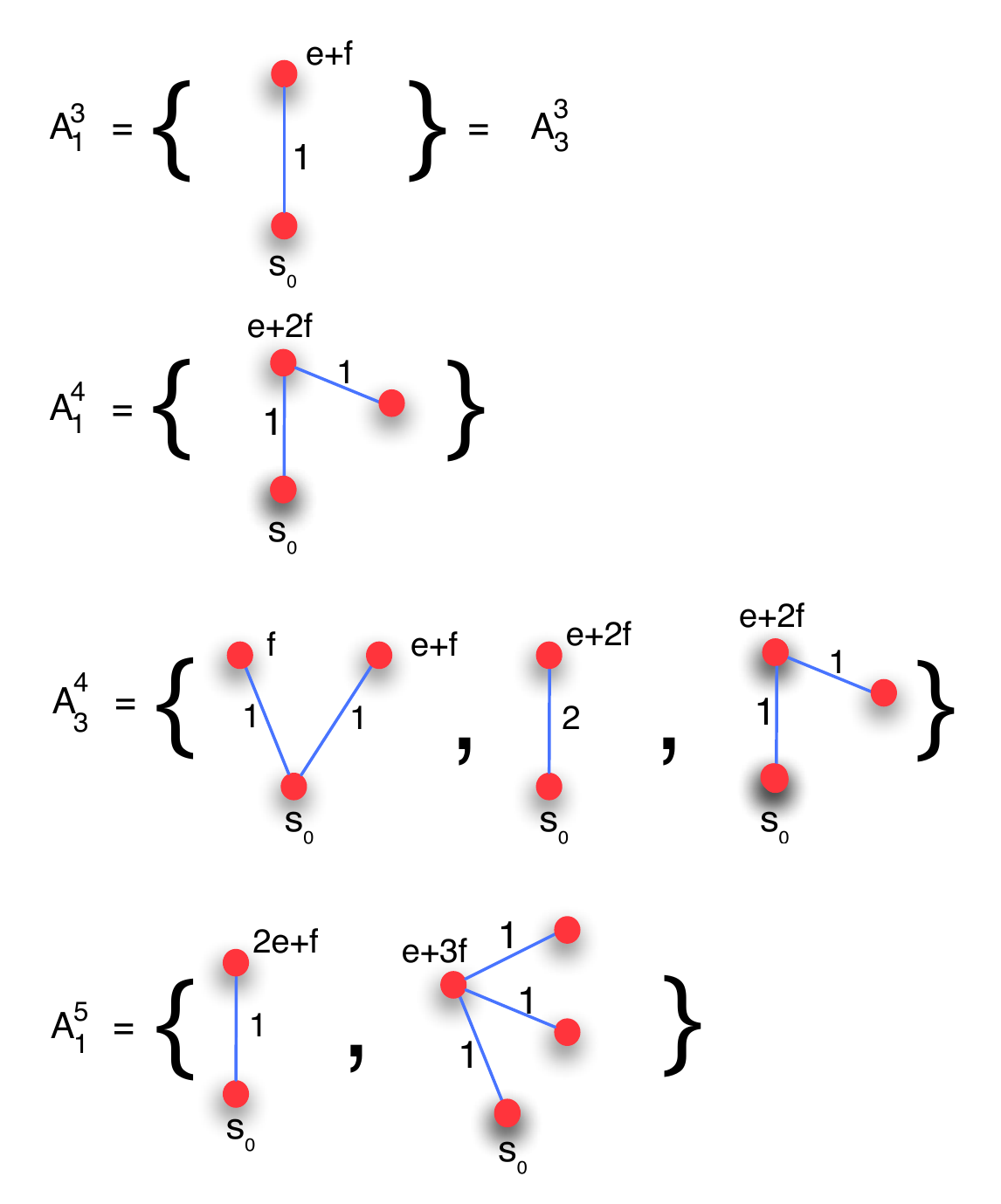}
\end{center}
\caption{Arbres deux-sph\'eriques d\'ecor\'es}
\end{figure}

\subsubsection{Relation avec l'invariant relatif}

Soit $\Sigma_2$ la surface rationnelle r\'egl\'ee de degr\'e deux, $e$ la classe d'une section holomorphe d'autointersection deux et $f$ la classe d'une fibre. \'Etant donn\'es $a,b \in \N$
et $\alpha, \beta$ des suites d'entiers positifs, on note $N^{ae+bf}_2 (\alpha, \beta)$ le nombre de courbes rationnelles de $\Sigma_2$, homologues \`a $ae+bf$, ayant $\alpha_i + \beta_i$
points de tangence d'ordre $i$ avec la section exceptionnelle de $\Sigma_2$ parmi lesquels $\alpha_i $ sont prescrits et qui passent par le nombre ad\'equat de points fix\'es.

\begin{theo}
\label{theocal2spher}
Soit $(X, \omega , c_X)$ une vari\'et\'e symplectomorphe \`a la quadrique ellipso\"{\i}de de dimension deux, $r, r_X, d \in \N$ satisfaisant la relation $r + 2r_X = 4d - 1$
et $h$ la classe d'une section plane r\'eelle de bidegr\'e $(1,1)$. Alors,
$$\chi^{dh}_r = \sum_{A \in {\cal A}_r^d} (-1)^{\# S_2 + 1} \text{mult} (A) F_{(r,0)} (\alpha_A^-, \beta_A^+) \prod_{s \in S_1 \setminus S_1^+} N_2^{g_A (s) e + k(s) f} (0, \beta_A)
\prod_{s \in S_1^+} N_2^{g_A (s) e + k(s) f} (e_{k(ss_0)} , \beta_A^0),$$
o\`u $(\alpha_A^-)_i$ (resp. $(\beta_A^+)_i$) vaut le  nombre d'ar\^etes de multiplicit\'e $i$ reliant $S_1^-$ (resp. $S_1^+$) \`a $s_0$, $( \beta_A)_i$ (resp. $( \beta_A^0)_i$)
vaut le nombre d'ar\^etes de multiplicit\'e $i$  adjacentes \`a $s$ (resp. moins un si l'ar\^ete reliant $s$ \`a $s_0$ est de multiplicit\'e $i$) et $F_{(r,0)} (\alpha_A^-, \beta_A^+)$ d\'esigne
l'invariant d\'efini dans le fibr\'e cotangent de la sph\`ere de dimension deux au \S \ref{subsubinvariants}.
\end{theo}

\begin{rem}
Un algorithme permettant le calcul de $\chi^d_r$, $1 \leq r \leq 4d-1$, est propos\'e par  E. Shustin dans \cite{Shu1}. Remarquons que cet invariant $\chi^d_r$ peut se d\'efinir
purement en termes de fractions rationnelles complexes. Lorsque $r = 4d-1$ par exemple, il compte alg\'ebriquement le nombre de fractions rationnelles $u = P/Q$, $P, Q \in \C [X]$ de degr\'es $d$,
modulo reparam\'etrage par les homographies r\'eelles de $PGL_2 (\R)$, telles que l'image $u(\R P^1)$ interpole un ensemble donn\'e g\'en\'erique de $4d-1$ points de la sph\`ere
de Riemann. Le signe en fonction duquel il convient de compter ces fractions rationnelles $u$ est pair si $u$ poss\`ede un nombre pair de points critiques dans chaque
h\'emisph\`ere $\C P^1 \setminus \R P^1$ et impair sinon. Il serait int\'eressant d'\'etudier cet invariant \`a l'aide de la th\'eorie des fractions rationnelles (un probl\`eme que j'avais propos\'e au MSRI
au printemps $2004$ lors d'une s\'eance de probl\`emes ouverts).
\end{rem}

{\bf D\'emonstration du Th\'eor\`eme \ref{theocal2spher} :}

On poursuit la strat\'egie g\'en\'erale \'enonc\'ee au $\S$ \ref{subsubsectstrat} en allongeant le cou d'une structure presque complexe g\'en\'erique jusqu'\`a briser la vari\'et\'e en deux morceaux
$T^* S^2$ et $X \setminus S^2$. Les courbes $J$-holomorphes rationnelles r\'eelles compt\'ees par $\chi^d_r$ se brisent en courbes \`a deux \'etages interpolant  $r$ 
points de $S^2$ et $r_X$ paires de points complexes conjugu\'es de $X \setminus S^2$. Il est apparu au cours de la d\'emonstration du Th\'eor\`eme \ref{theocong1} que 
ces courbes \`a deux \'etages
sont cod\'ees par les arbres deux-sph\'eriques d\'ecor\'es $A \in {\cal A}_r^d$. Il s'agit donc de d\'enombrer les  courbes \`a deux \'etages qui sont cod\'ees par un arbre donn\'e $A \in {\cal A}_r^d$, puis
de d\'enombrer les  courbes $J$-holomorphes rationnelles r\'eelles compt\'ees par $\chi^d_r$ qui d\'eg\'en\`erent sur une courbe \`a deux \'etages donn\'ee. Le nombre de fa\c{c}ons de r\'epartir les
points complexes conjugu\'es parmi les composantes de  la courbe \`a deux \'etages qui se trouvent dans $X \setminus S^2$ a \'et\'e calcul\'e dans la d\'emonstration du Th\'eor\`eme
\ref{theocong1} et vaut $2^{\sum_{s \in S_1^+} \# f_A (s) + \sum_{s \in S_1 \setminus S_1^+} \max (\# f_A (s) - 1 , 0)}$. Les composantes cod\'ees par $S_1 \setminus S_1^+$ sont rigides
avec leurs conditions d'incidence, il y en a $\prod_{s \in S_1 \setminus S_1^+} N_2^{g_A (s) e + k(s) f} (0, \beta_A)$. Puis, il y a  $m_1^- (A)$ fa\c{c}ons de choisir les orbites de Reeb prescrites
de la courbe cod\'ee par $s_0$. Le nombre de courbes r\'eelles cod\'ees par $s_0$ satisfaisant nos conditions d'incidence et compt\'ees avec signe vaut $F_{(r,0)} (\alpha_A^-, \beta_A^+)$.
Il y a alors $m_1^+ (A)$  fa\c{c}ons de choisir la mani\`ere de connecter les courbes cod\'ees par $S_1^+$ aux  orbites de Reeb rest\'ees libres de la courbe cod\'ee par $s_0$. 
Enfin, chaque sommet de $S_2$ autre que $s_0$ code un plan $J$-holomorphe de $T^* S^2$ asymptote \`a une orbite de Reeb simple. Il y a deux tels plans pour une structure presque
complexe g\'en\'erique $J$ de $T^* S^2$ qui sont les deux relev\'es du plan de $T^* \R P^2$ asymptote \`a une orbite de Reeb double et se compactifient en les deux droites de la quadrique
complexe passant par un point donn\'e. Il y a donc $2^{ \#S_2 -1}$ fa\c{c}ons de choisir les plans cod\'es par les \'el\'ements de $S_2 \setminus \{ s_0 \}$. Ceci fournit $2^{\sum_{s \in S_1^+} \# f_A (s) + 
\sum_{s \in S_1 \setminus S_1^+} \max (\# f_A (s) - 1 , 0) + \#S_2 -1} m_1^+ (A)  m_1^- (A) F_{(r,0)} (\alpha_A^-, \beta_A^+)$ $\prod_{s \in S_1 \setminus S_1^+} N_2^{g_A (s) e + k(s) f} (0, \beta_A) 
\prod_{s \in S_1^+} N_2^{g_A (s) e + k(s) f} (e_{k(ss_0)}  , \beta_A^0)$ courbes cod\'ees par un arbre donn\'e $A \in {\cal A}_r^d$. Or, d'apr\`es le th\'eor\`eme de
recollement de th\'eorie symplectique des champs \cite{Bour}, il y a $ \prod_a k(a)$ courbes $J$-holomorphes rationnelles r\'eelles compt\'ees par $\chi^d_r$ qui d\'eg\'en\`erent sur une courbe \`a deux 
\'etages donn\'ee. Le r\'esultat d\'ecoule \`a pr\'esent du fait que chaque plan cod\'e par $S_2 \setminus \{ s_0 \}$ intersecte $S^2$ en un point et contribue donc \`a la masse des courbes
$J$-holomorphes rationnelles r\'eelles en question, d'o\`u le signe $(-1)^{\# S_2 + 1}$. $\square$

\begin{cor}
\label{corcalc2spher}
Soit $(X, \omega , c_X)$ une vari\'et\'e symplectomorphe \`a la quadrique ellipso\"{\i}de de dimension deux. On note $h$ la classe d'une section plane r\'eelle de bidegr\'e $(1,1)$. Alors,
$\chi^{2h} (T) = 2T^3 + 4T^5 + 6T^7$, $\chi^{3h} (T) = 16T + 16T^2 + o(T^3)$,  $\chi^{4h} (T) = -256T + 320T^3 + o(T^4)$ et $\chi^{5h} (T) = 26880T + o(T^2)$.
\end{cor}

\begin{rem}
Les valeurs $\chi^{4h}_1$, $\chi^{4h}_3$ et $\chi^{5h}_1$ \'enonc\'ees dans la Proposition $2.3$ de \cite{Wels4} sont incorrectes et corrig\'ees ici, note \cite{Wels4} qui fut d'ailleurs soumise
en l'\'etat en janvier et non d\'ecembre $2006$.
\end{rem}

{\bf D\'emonstration du Corollaire \ref{corcalc2spher} :}

Le calcul de $\chi^{2h} (T)$ d\'ecoule imm\'ediatement du Lemme \ref{lemmecalc1}, puisque les seules composantes de $X \setminus L$ apparaissant sont des fibres.
Les arbres intervenant dans la d\'emonstration de ce Corollaire \ref{corcalc2spher} sont
repr\'esent\'es dans la Figure \ref{figarbres2spheriques}. Lorsque $d=3$ et $r \leq 3$, un seul arbre d\'ecor\'e intervient. Le Th\'eor\`eme \ref{theocalcproj} fournit
$\chi^{3h}_1 = 2^4 F (e_1 , 0) N^{e+f} (0 , e_1) = 16$ et $\chi^{3h}_3 = 2^4 F (0 , e_1) N^{e+f} (e_1 , 0) = 16$. Lorsque $d=4$ et $r=1$, un seul arbre deux-sph\'erique d\'ecor\'e $A$ intervient
dans le calcul de $\chi^{dh}_r$. On obtient $\chi^{4h}_1 = - 2^7 2 F (e_1 , 0) N^{e+2f} (0 , 2e_1) = -256$, puisque $m_1^- (A) = 2$. L'ensemble ${\cal A}_3^{4h}$ contient quant \`a lui
trois arbres deux-sph\'eriques d\'ecor\'es. La contribution du premier arbre donn\'e par la Figure \ref{figarbres2spheriques} vaut $C_6^1 2^4 F (2e_1 , 0) N^{f} (0 , e_1) N^{e+f} (0 , e_1) =  192$ 
puisqu'il  y a $C_6^1$ fa\c{c}ons de choisir la fonction $f_A$ ; celle du second vaut $2^5 2 F (e_2 , 0) N^{e+2f} (0 , e_2) = 2^8 N^{e+2f} (e_2 , 0) = 256$ et celle du troisi\`eme
$- 2^7 F (0 , e_1) N^{e+2f} (e_1 , e_1) = -128$, de sorte que $\chi^{4h}_3 = 192 + 256 - 128 = 320$. L'ensemble ${\cal A}_1^{5h}$ contient 
deux arbres deux-sph\'eriques d\'ecor\'es. La contribution du premier arbre donn\'e par la Figure \ref{figarbres2spheriques} vaut $2^8 F (e_1 , 0) N^{2e+f} (0 , e_1) = 2^8 93$, puisque d'apr\`es
le Th\'eor\`eme $6.8$ de \cite{Vak}, $N^{2e+f} (0 , e_1) = 93$.  La contribution du deuxi\`eme arbre vaut $2^{10} 3 F (e_1 , 0) N^{e+3f} (0 , 3e_1) = 3072$, de sorte que $\chi^{5h}_1 =
23808 + 3072 = 26880$. $\square$

\subsection{Calculs dans l'ellipso\"{\i}de de dimension trois}
\label{subsec3spher}

\subsubsection{Arbres trois-sph\'eriques}

Soient $r, r_X$ et $d$ trois entiers naturels satisfaisant la relation $2r + 4r_X = 3d$, laquelle impose que $d$ soit pair.

\begin{defi}
Un arbre trois-sph\'erique est un arbre connexe fini dont toutes les ar\^etes sont \'etiquet\'ees par des entiers strictement positifs. De plus, un tel arbre poss\`ede une racine $s_0$ et tous les sommets qui
lui sont adjacents sont monovalents.
\end{defi}

En particulier, la distance maximale d'un sommet \`a $s_0$ vaut un.
On note ${\cal C}_r$ l'ensemble des arbres trois-sph\'eriques qui satisfont
$$4\sum_{a \in {\cal A}(s_0)} k(a) - 2v(s_0) \leq 2r \leq  4\sum_{a \in {\cal A}(s_0)} k(a) + 2v (s_0) \text{ et } r = v (s_0)  \mod (2),$$
o\`u $v(s)$ d\'esigne la valence d'un sommet $s$, $ {\cal A}(s)$ l'ensemble des ar\^etes adjacentes \`a $s$ et $k (a)$ d\'esigne la multiplicit\'e de l'ar\^ete $a$.
Notons \'egalement $k_s$ la somme des multiplicit\'es des ar\^etes adjacentes au sommet $s$ et $k$ la multiplicit\'e totale de toutes les ar\^etes de l'arbre.
On pose $r_L (s_0) = \frac{1}{4} \big(    4\sum_{a \in {\cal A}(s_0)} k(a) + 2v (s_0) - 2r     \big)$, de sorte que $0 \leq r_L (s_0) \leq v (s_0)$. Enfin, on note $S_1$ l'ensemble
des sommets adjacents \`a $s_0$.

\begin{defi}
Un arbre trois-sph\'erique d\'ecor\'e est un arbre trois-sph\'erique $A \in {\cal C}_r$ \'{e}quip\'{e} d'une partition $S_1^+ \sqcup S_1^-$ de l'ensemble des sommets adjacents \`a $s_0$ telle que
$\# S_1^- = r_L (s_0)$ et $\# S_1^+ = v(s_0) - r_L (s_0)$. Cet arbre est de plus \'equip\'e des fonctions :
\begin{itemize}
\item $f_A : S_1 \to {\cal P} (\{1, \dots , r_X \})$ satisfaisant $f_A (s) \cap f_A (s') = \emptyset$ d\`es que $s \neq s'$ et $\cup_{s \in S_1} f_A (s) = \{1, \dots , r_X \}$.
\item  $g_A : S_1 \to \N$ telle que $2k + 2\sum_{s \in S_1} g_A (s) = d$ et pour tout $s \in S_1^+$ (resp. $s \in S_1^-$ ), $3g_A (s) + k_s + 1 = 2\# f_A (s) + 2$
(resp. $3g_A (s) + k_s + 1 = 2\# f_A (s)$).
\end{itemize}
\end{defi}

On note ${\cal C}_r^d$ l'ensemble des arbres trois-sph\'eriques d\'ecor\'es, c'est un ensemble fini. Soit $A \in {\cal C}_r^d$, on pose 
$m_1^+ (A)$ le nombre d'injections $\phi : \{s \in S_1^+ \; \vert \; f_A (s) \neq \emptyset \} \to  {\cal A}^+ (s_0)$ satisfaisant
$k (\phi (s)) = k(ss_0)$ pour tous les sommets $s \in S_1^+$.   On pose alors 
$$\text{mult} (A) = 2^{\sum_{s \in S_1^+} \# f_A (s) + \sum_{s \in S_1^-} \max (\# f_A (s) - 1 , 0)} m_1^+ (A) \prod_a k(a),$$
c'est la {\it multiplicit\'e} de l'arbre trois-sph\'erique $A \in {\cal C}_r^d$.

\subsubsection{Relation avec l'invariant relatif}

Soit $Y$ la vari\'et\'e r\'egl\'ee $P({\cal O}_Q (1,1) \oplus {\cal O}_Q)$ sur la quadrique de dimension deux $Q$ et $f$ la classe d'une fibre de $Y$. \'Etant donn\'es $a,b,c  \in \N$
et $\alpha, \beta$ des suites d'entiers positifs, on note $N^{(a,b)+cf}_3 (\alpha, \beta)$ le nombre de courbes rationnelles de $Y$, homologues \`a $(a,b)+cf$, ayant $\alpha_i + \beta_i$
points de tangence d'ordre $i$ avec la section exceptionnelle $P({\cal O}_Q)$ de $Y$ parmi lesquels $\alpha_i $ sont prescrits et qui passent par le nombre ad\'equat de points fix\'es.

\begin{theo}
\label{theocal3spher}
Soient $(X, \omega , c_X)$ une vari\'et\'e symplectomorphe \`a la quadrique ellipso\"{\i}de de dimension trois et $r, r_X, d \in \N$ satisfaisant la relation $2r + 4r_X = 3d$. Alors,
$$\chi^{d}_r = \sum_{A \in {\cal C}_r^d} \text{mult} (A) F_{(r,0)} (\alpha_A^-, \beta_A^+) \prod_{s \in S_1^-} \sum_{a+b = g_A (s)} N_3^{(a,b) + k_s f} (0, e_{k_s}) 
\prod_{s \in S_1^+} \sum_{a+b = g_A (s)} N_3^{(a,b)+ k_s f} (e_{k_s} , 0),$$
o\`u $(\alpha_A^-)_i$ (resp. $(\beta_A^+)_i$) vaut le  nombre d'ar\^etes de multiplicit\'e $i$ reliant $S_1^-$ (resp. $S_1^+$) \`a $s_0$ et $F_{(r,0)} (\alpha_A^-, \beta_A^+)$ d\'esigne
l'invariant d\'efini dans le fibr\'e cotangent de la sph\`ere de dimension trois au \S \ref{subsubinvariants}.
\end{theo}

{\bf D\'emonstration :}

On poursuit la strat\'egie g\'en\'erale \'enonc\'ee au $\S$ \ref{subsubsectstrat} en allongeant le cou d'une structure presque complexe g\'en\'erique jusqu'\`a briser la vari\'et\'e en deux morceaux
$T^* S^3$ et $X \setminus S^3$. Les courbes $J$-holomorphes rationnelles r\'eelles compt\'ees par $\chi^d_r$ se brisent en courbes \`a deux \'etages interpolant  $r$ 
points de $S^3$ et $r_X$ paires de points complexes conjugu\'es de $X \setminus S^3$. Il est apparu au cours de la d\'emonstration du Th\'eor\`eme \ref{theocong2} que 
ces courbes \`a deux \'etages sont cod\'ees par les arbres trois-sph\'eriques d\'ecor\'es $A \in {\cal C}_r^d$. Il s'agit donc de d\'enombrer les  courbes \`a deux \'etages qui sont cod\'ees 
par un arbre donn\'e $A \in {\cal C}_r^d$, puis
de d\'enombrer les  courbes $J$-holomorphes rationnelles r\'eelles compt\'ees par $\chi^d_r$ qui d\'eg\'en\`erent sur une courbe \`a deux \'etages donn\'ee. Le nombre de fa\c{c}ons de r\'epartir les
points complexes conjugu\'es parmi les composantes de  la courbe \`a deux \'etages qui se trouvent dans $X \setminus S^3$ a \'et\'e calcul\'e dans la d\'emonstration du Th\'eor\`eme
\ref{theocong2} et vaut $2^{\sum_{s \in S_1^+} \# f_A (s) + \sum_{s \in S_1^-} \max (\# f_A (s) - 1 , 0)}$. Les composantes cod\'ees par $S_1^-$ sont rigides
avec leurs conditions d'incidence, il y en a $\prod_{s \in S_1^-} \sum_{a+b = g_A (s)} N_3^{(a,b) + k_s f} (0, k_s) $.  Le nombre de courbes r\'eelles cod\'ees par $s_0$ 
satisfaisant nos conditions d'incidence et compt\'ees avec signe vaut $F_{(r,0)} (\alpha_A^-, \beta_A^+)$.
Il y a alors $m_1^+ (A)$  fa\c{c}ons de choisir la mani\`ere de connecter les courbes cod\'ees par $S_1^+$ aux  orbites de Reeb rest\'ees libres de la courbe cod\'ee par $s_0$.  
Ceci fournit $2^{\sum_{s \in S_1^+} \# f_A (s) + \sum_{s \in S_1^-} \max (\# f_A (s) - 1 , 0)} m_1^+ (A) F_{(r,0)} (\alpha_A^-, \beta_A^+) \prod_{s \in S_1^-} \sum_{a+b = g_A (s)} N_3^{(a,b) + k_s f} (0, k_s)$\\
$\prod_{s \in S_1^+} \sum_{a+b = g_A (s)} N_3^{(a,b)+ k_s f} (k_s , 0) $ courbes cod\'ees par un arbre donn\'e $A \in {\cal C}_r^d$. Or, d'apr\`es le th\'eor\`eme de
recollement de th\'eorie symplectique des champs \cite{Bour}, il y a $ \prod_a k(a)$ courbes $J$-holomorphes rationnelles r\'eelles compt\'ees par $\chi^d_r$ qui d\'eg\'en\`erent sur 
une courbe \`a deux \'etages donn\'ee. Le r\'esultat en d\'ecoule. $\square$

\begin{cor}
\label{corcalc3spher}
Soit $(X, \omega , c_X)$ une vari\'et\'e symplectomorphe \`a la quadrique ellipso\"{\i}de de dimension trois. Alors,
$\chi^{2}_1 = - 1$, $\chi^{6}_1 = 0$ et $\chi^{10}_1 = -896$.
\end{cor}

\begin{rem}
L'invariant $\chi^{4k}_0$, $k \in \N$, n'est pas d\'efini. Toutefois, d'apr\`es le Th\'eor\`eme \ref{theomin2}, lorsque les points complexes conjugu\'es sont suffisamment proches
d'une section hyperplane r\'eelle disjointe du lieu r\'eel de $X$, il n'y a aucune courbe rationnelle r\'eelle ayant une partie r\'eelle non-vide qui satisfait nos conditions d'incidence.
Remarquons aussi que l'invariant d\'efini dans \cite{Wels2}, \cite{Wels3} a \'et\'e \'etendu aux vari\'et\'es symplectiques r\'eelles de dimension six non fortement semi-positives dans 
\cite{Sol} (voir aussi \cite{Cho}) et calcul\'e dans le cas des quintiques r\'eelles de $\C P^4$ dans \cite{PSW}.
\end{rem}

{\bf D\'emonstration du Corollaire \ref{corcalc3spher} :}

Il n'y a qu'une seule section plane de $X$ qui passe par trois points. De plus, l'\'etat spinoriel de cette conique vaut $-1$, de sorte que $\chi^{2}_1 = - 1$.
En effet, l'\'etat spinoriel d'une conique dans une quadrique de dimension deux vaut $-1$ et le fibr\'e normal d'une section hyperplane r\'eelle de $X$ est trivial en restriction 
\`a sa partie r\'eelle. L'annulation de $\chi^{6}_1$ tient au fait qu'il n'y a pas de courbe rationnelle de bidegr\'e $(a,b)$ qui passe par quatre points dans la quadrique de
dimension deux. Enfin, l'ensemble $ {\cal C}_1^{10}$ ne contient qu'un seul arbre qui n'a que deux sommets, est de multiplicit\'e $2^6$, pour lequel $S_1^+$ est vide et 
$F (\alpha_A^-, \beta_A^+) = -1$
d'apr\`es ce qui pr\'ec\`ede. Par sept point de la quadrique $Q$ de dimension deux, il passe douze courbes rationnelles de bidegr\'e $(2,2)$, une courbe de bidegr\'e $(3,1)$ et
une de bidegr\'e $(1,3)$. La restriction du r\'eglage $Y \to Q$ \`a ces courbes est isomorphe \`a la surface r\'egl\'ee rationnelle de degr\'e quatre $\Sigma_4$. Par sept points de $\Sigma_4$,
il ne passe qu'une seule courbe rationnelle homologue \`a $e+f$. On en d\'eduit que $N_3^{(3,1) + f} (0, e_1) = N_3^{(1, 3) + f} (0, e_1) = 1$ et $N_3^{(2,2) + f} (0, e_1) = 12$.
D'o\`u $\chi^{10}_1 = -896$. $\square$

\addcontentsline{toc}{part}{\hspace*{\indentation}Bibliographie}


\bibliographystyle{abbrv}

\vspace{0.7cm}
\noindent 
Unit\'e de math\'ematiques pures et appliqu\'ees de l'\'Ecole normale sup\'erieure de Lyon,\\
 CNRS - Universit\'e de Lyon.

\end{document}